\documentclass[11pt,final,
titlepage]{article}

\usepackage{tikz}
\usepackage{epsfig}
\usepackage{amsfonts}

\usepackage{amssymb}
\usepackage{amsmath}
\usepackage{amsthm}
\usepackage{latexsym}
\usepackage{graphicx}
\usepackage{bm}
\usepackage{tabularx}
\usepackage{booktabs}
\usepackage{enumitem}

\usepackage[titletoc]{appendix}
\usepackage[numbers,sort&compress]{natbib}
\usepackage{caption}
\usepackage{subcaption}

\usepackage{mathtools}
\usepackage{multirow}
\usepackage{authblk}
\newcommand{\field}[1]{\mathbb{#1}}
\def\A{\field{A}}                                
\def\X{\field{X}}                                

\def\Z{\field{Z}}

\catcode`\@=11

\usepackage[margin=1in]{geometry}

\usepackage{setspace}
\numberwithin{equation}{section}
\numberwithin{table}{section}
\numberwithin{figure}{section}

\theoremstyle{plain}
\newtheorem{theorem}{Theorem}[section]
\newtheorem{lemma}[theorem]{Lemma}
\newtheorem{corollary}[theorem]{Corollary}
\newtheorem{proposition}[theorem]{Proposition}

\theoremstyle{definition}
\newtheorem{definition}{Definition}[section]

\newtheorem{example}[definition]{Example}
\newtheorem{condition}[definition]{Condition}

\theoremstyle{remark}
\newtheorem{remark}{Remark}

\graphicspath{{./}{Figs/}}

\usepackage[colorlinks=true,linkcolor=magenta,citecolor=blue,pagebackref]{hyperref}
\usepackage[noabbrev,capitalize]{cleveref}

\title{
Balancing Independent and Collaborative Service}
\author[1]{Shuwen Lu}
\author[2] {Mark E. Lewis}
\author[3]{Jamol Pender}

\affil[1]{\small Department of Systems Engineering, Cornell University, Ithaca, NY 14850, USA {\tt\small sl3243@cornell.edu}}%
\affil[2,3]{\small School of Operations Research and Information Engineering,
Cornell University, Ithaca, NY 14850, USA
{\tt\small mark.lewis@cornell.edu, jjp274@cornell.edu}}

\begin{document}
\maketitle

\begin{abstract}
    We study a two-type server queueing system where flexible Type-I servers, upon their initial interaction with jobs, decide in real time whether to process them independently or in collaboration with dedicated Type-II servers. Independent processing begins immediately, as does collaborative service if a Type-II server is available.  Otherwise, the job and its paired Type-I server wait in queue for collaboration.  Type-I servers are non-preemptive and cannot engage with new jobs until their current job is completed.

    We provide a complete characterization of the structural properties of the optimal policy for the clearing system. In particular, an optimal control is shown to follow a threshold structure based on the number of jobs in the queue before a Type-I first interaction and on the number of jobs in either independent or collaborative service. 

    We propose simple threshold heuristics, based on linear approximations, for real-time decision-making. In much of the parameter and state spaces, we establish theoretical bounds that compare the thresholds proposed by our heuristics to those of optimal policies and identify parameter configurations where these bounds are attained. Outside of these regions, the optimal thresholds are infinite. Numerical experiments further demonstrate the accuracy and robustness of our heuristics, particularly when the initial queue length is high. Our proposed heuristics achieve costs within 0.5\% of the optimal policy on average and significantly outperform benchmark policies that exhibit extreme sensitivity to system parameters, sometimes incurring costs exceeding 100\% of the optimal.
\end{abstract}
\section{Introduction} \label{sec:intro}
We consider a controlled queueing system equipped with two types of servers: flexible servers and dedicated servers. Flexible servers can make real-time, sequential decisions about whether to perform independent service or collaborate with dedicated servers. Collaboration may enhance the speed and/or quality of service but can introduce delays due to the limited availability of dedicated servers. 
In manufacturing, for example, engineers or technicians may either complete a task independently or operate (or supervise) a specialized machine. While machine-assisted operations offer higher efficiency and convenience, both the worker and the task may have to wait for a scarce machine to become available. Similar dynamics arise in service systems such as call centers, ticketing platforms, and in-person queues, where agents may resolve issues directly or escalate for a collaborative resolution involving both frontline staff, higher-level experts, and the customer in a comprehensive three-way interaction. 
Although such collaboration often improves service quality and customer satisfaction, it also results in additional waiting if specialist participation is not instantly accessible.

Such queueing structures with decisions are ubiquitous in industry and everyday life. 
In healthcare delivery, for instance, frontline providers (e.g., nurse practitioners) may independently manage low-acuity patients but can also initiate collaborative consultations with physicians, often via teleconferencing, to improve diagnostic accuracy and patient confidence \cite{pelone2017interprofessional, pappas2019diagnosis}. In practice, such decisions are often based on individual experience rather than objective policies that account for system-wide optimality, underscoring the need for mathematical models that guide effective decision-making.

Importantly, the decision-makers in these systems need not be human. A prominent example arises in cloud computing, where computational tasks are allocated between general-purpose CPUs and specialized GPU accelerators, forming a queueing system with two distinct types of servers. More specifically, large-scale matrix computations typically proceed sequentially through preprocessing, the core multiplication or solve, and post-processing or refinement. The core computation usually dominates runtime, while the other stages are comparatively lightweight \cite{islam2025improving, elafrou2017performance}. 
This aligns with our queueing abstraction, which treats upstream and downstream processing times as negligible, focusing on the detailed efforts in between. 
Existing libraries either market themselves as GPU-only solvers that offload main kernels to GPUs \cite{naumov2015amgx}, or provide hybrid CPU–GPU execution guided by heuristics or manual configuration \cite{balay2024petsc, heroux2005overview, nvidia2024cublas}. These approaches highlight the opportunity for a principled CPU–GPU collaborative framework—particularly when unstructured sparsity or GPU overhead reduces acceleration gains—where the central problem is to determine when collaboration should be triggered, given time and quality trade-offs.

In these scenarios, where flexible servers are capable of completing jobs independently, the objective is to identify a control policy that is impartial to the preferences of individual servers or jobs and maximizes the \textit{overall quality of service} for all jobs. This leads to several key questions:
\begin{itemize}
    \item How can we quantitatively define the measure `service quality' to account for the time each job spends in queues and at service stations, as well as the broader impact on other jobs?
    \item Under what conditions should servers prefer collaborative service over independent service?
    \item If an optimal policy is computationally infeasible, can we develop simple, effective, and robust heuristic rules for real-time decision-making?
\end{itemize}

To address these challenges, we model the sequential decision-making scenario using a Markov Decision Process (MDP) model. Each job incurs different costs per unit of time spent in queues or service stations, enabling us to measure service quality through the expected total cost charged by the system. We refer to these costs as \textit{holding costs} to align with the literature. Maximizing overall service quality thus translates to minimizing the total expected cost. We also examine a clearing system model that clears existing jobs, as in practice, it is of interest to investigate systems with daily job completion requirements (in manufacturing) or how a system recovers to its normative state after a surge of arrivals \cite{minh2011performance}. This approach results in an MDP clearing system model with holding costs, where the goal is to determine the control policy that minimizes the total expected cost of returning to an empty state.

To answer the second question above, we completely characterize the structural properties of optimal policies across the entire parameter space. We demonstrate that the optimal policy follows a threshold structure in the number of jobs in the queue, with thresholds that may extend from zero to infinity. Additionally, we analyze how the optimal policy evolves as the number of jobs at service stations varies, revealing that the policy exhibits \textit{diagonal monotonicity} with respect to job counts at the two service stations.

In real-time decision-making, obtaining an optimal policy may be impractical, particularly when flexible servers, who are often people, are not programmable or when computational constraints arise. For example, if the initial queue length is large, solving for the optimal policy can become computationally expensive or even intractable, especially when the number of servers is also large due to the rapid growth of the problem size. 
To address this challenge, we develop threshold heuristic policies that are computationally efficient and suitable for real-time application. 
In much of the parameter and state spaces, we derive theoretical bounds comparing the thresholds suggested by our heuristics relative to those of an optimal policy and provide parameter configurations where these bounds are achieved. Outside of these regions, the optimal thresholds are infinite.

Finally, numerical analysis demonstrates the strong performance of our proposed heuristics, especially when the initial queue is long, achieving near-optimal accuracy (within 1\% of the optimal policy) and exhibiting robustness across parameter variations. In contrast, existing benchmark policies are highly sensitive to parameter changes and can incur costs exceeding 100\% of the optimal.

\subsection{Operational implications}
Based on our model of real-time decision-making between immediate service and a higher-quality but potentially delayed alternative, we outline managerial implications or insights for practical applications:

\begin{itemize}
    \item \textbf{Model Generality.} 
    Our framework applies to a wide range of settings where resources in the second option are limited (but reusable). In service systems, this scenario occurs when there are expert staff or higher-level specialists with restricted availability. In manufacturing and supply chains, the same can refer to specialized equipment with limited supply or facilities with capacity constraints that limit the number of simultaneous jobs.
    
    Furthermore, the concept of service quality extends beyond traditional measures of product excellence. It can also encompass factors such as customer satisfaction, operational efficiency through automation, or convenience from computer-assisted processes.
    \item \textbf{Job Completion and System Recovery.}
    Our model emphasizes how systems recover from demand surges and ensure timely job completion, especially when the initial load is high. This insight is particularly relevant for businesses dealing with fluctuating demand, helping them anticipate bottlenecks and allocate resources accordingly. 
    \item \textbf{System-Wide Optimality Over Individual Preferences.} 
    Our model prioritizes overall system performance over individual job or server preferences. This principle helps managers assess whether immediate demands from individual jobs or servers align with long-term system optimality.
    \item \textbf{Threshold-Based Decision Policies.}
    We establish that the optimal control strategy follows a threshold-based policy across the whole parameter space (with infinite thresholds under certain conditions). Our findings indicate when it is more efficient to serve jobs independently versus collaboratively. Additionally, our results suggest that as more jobs wait for limited resources rather than receiving immediate service, it naturally discourages further queueing. This insight may help managers dynamically allocate resources across multiple locations.
    \item \textbf{Practical and Scalable Heuristic Implementation.}
    We propose adaptive threshold heuristics that apply across all parameter configurations, enabling real-time decision-making. These heuristics are easy to compute --- even with a basic calculator --- making them particularly practical and scalable when exact solutions are infeasible or learning-based approaches are too computationally expensive.
    \item \textbf{Heuristic Accuracy and Robustness.}
    We establish theoretical performance guarantees by comparing our heuristically estimated thresholds to the optimal ones (when finite). Numerical analysis further demonstrates their robustness and near-optimal accuracy, particularly in high initial queueing scenarios.
    \item \textbf{Actionable Insights for Decision-Makers and Stakeholders.}
    Our findings provide decision-makers with guidance on adjusting policies in response to parameter changes and implementing heuristics effectively. Additionally, we offer recommendations to managers and stakeholders on when investing in the expansion of limited resources --- such as hiring more expert staff or acquiring specialized equipment --- would be most beneficial. For instance, there are certain parameter regions where collaborative service is consistently required within certain initial queue lengths, making additional resources essential.
\end{itemize}

\subsection{Organization of paper}
The remainder of the paper is organized as follows: Related literature is covered in Section \ref{sec:lit-rev}. Section \ref{sec:description} formally introduces the model, with Section \ref{sec:main-results} presenting the key results on decision-making and Section \ref{sec:support-results} providing their supporting results. 
Section \ref{sec:heuristics} details the heuristic design and its theoretical justification. Numerical experiments are conducted in Section \ref{sec:numerical}. Finally, the paper concludes in Section \ref{sec:conclusions} with a summary of the findings and future research directions.

\section{Literature Review}\label{sec:lit-rev}
To the best of our knowledge, despite its widespread applicability, no prior research has addressed control policies for the type of queueing system studied in this work.
In this section, we highlight the unique queuing and policy structures of our work, while situating them within the broader literature.

Existing studies on queueing systems involving both flexible and dedicated servers (with potential collaboration on a single job) typically focus on scheduling policies for the flexible servers. In tandem queues, where a single class of arrivals progresses through multiple service stages, scheduling policies (or server allocation strategies) determine which stage(s) the flexible servers should prioritize \cite{farrar1993optimal, wu2006dynamic, wu2008heuristics, pandelis2008optimal1, andradottir2007dynamic, pandelis2008optimal2, zayas2016dynamic, zayas2019policies, papachristos2018optimal, yarmand2015maximizing, zou2018asymptotically}.
In contrast, single-stage systems with multi-class jobs allocate flexible servers across different job types \cite{bhulai2003queueing, bhulai2012optimal, ansari2019optimal}. A canonical example is the N-networks, dedicated servers focus exclusively on Type-I jobs, while cross-trained flexible servers can assist with either Type-I or Type-II jobs \cite{bell2001dynamic, ahn2004optimal, down2010n, tezcan2010dynamic, mathew2021two}. 
When both job classes are supported by dedicated servers and flexible servers are available to assist either, the structure corresponds to the M-model, as analyzed in \cite{chen2024managing}.
Our model diverges from the above in both queueing structure and control dynamics. Flexible (Type-I) servers not only choose between two service types --- independent or collaborative service with a dedicated (Type-II) server --- but also remain engaged for the full duration of the selected mode.

There is extensive literature examining queueing systems in which decisions involve selecting between service types or operating modes \cite{heyman1968optimal, borthakur1987poisson, lippman1975applying, george2001dynamic, badian2021optimal, rykov2004optimal, kim2011managing, efrosinin2016optimal, andradottir2021optimizing, yu2023optimal}. 
Within single-server queues, early work considers systems where the server operates in multiple modes (e.g., working, idling, or being turned off) with policies optimizing transitions among these modes \cite{heyman1968optimal, borthakur1987poisson}. 
Subsequent work extends to controlling service rates dynamically in response to system load \cite{lippman1975applying, george2001dynamic, badian2021optimal}. 
Other studies include allocating faster or more reliable servers based on current queue states \cite{rykov2004optimal, kim2011managing, efrosinin2016optimal}.
By contrast, our setting features two server types with distinct roles and resource asymmetry.
As outlined in Section \ref{sec:intro}, real-world examples include senior engineers with broader technical knowledge compared to junior technicians, or capital-intensive assembly machines with limited availability. In such cases, opting for collaborative service may require the flexible server (and its assigned job) to wait for an available dedicated server. This delay not only prevents the flexible server from initiating a new job but also propagates congestion throughout the system. 
The resulting interdependence introduces additional complexity, further distinguishing our work from prior work. 
Unlike our focus on selecting between two operational modes within the queueing system, a related but distinct line of research considers time allocation between queueing tasks and other responsibilities. For example, \citet{andradottir2021optimizing} analyze how an attending physician should allocate time between supervising residents and managing personal responsibilities in a two-stage queueing model, with extensions to multiple supervisors (attending physicians) in \cite{yu2023optimal}.

Clearing system models are widely used to analyze control policies aimed at minimizing the total expected cost required to process all jobs in a system. In the context of two-station tandem queues, \citet{farrar1993optimal} establishes the existence of an optimal transition-monotone policy for a system with one dedicated and one reconfigurable server. \citet{wu2006dynamic} extend this result, demonstrating that similar structural properties hold even when multiple dedicated and reconfigurable servers are present. In a system with two flexible servers and no dedicated ones, \citet{ahn1999optimal} derive necessary and sufficient conditions under which it is optimal to allocate both servers to the same station. Related work on parallel queues includes \citet{ahn2004optimal}, who characterize the optimal scheduling policy in a system with one dedicated and one flexible server, assuming no external arrivals.

In our study, we establish the existence of a threshold-structured optimal policy by partitioning the parameter and state space. The direction of monotonicity (non-increasing or non-decreasing) and the values of the thresholds (one, other finite values, or infinity) depend on specific parameter relationships. Based on these structural insights, we design efficient and accurate heuristics of the same threshold type to approximate the optimal decisions. 
Despite their simplicity, threshold policies describe many optimal solutions --- including several in the aforementioned literature --- and have thus become some of the most extensively studied and widely implemented policies in queueing control \cite{bhulai2003queueing, badian2021optimal, pang2015logarithmic, yu2023optimal, iravani1997two}.
For example, \citet{bhulai2003queueing} examine a system with two job types served by a common pool of servers with a waiting time constraint on the first type and show that the optimal allocation policy follows a threshold structure on the number of available servers, assuming equal service rates for the two job types. In an M/M/1 queue with a removable server capable of dynamically choosing from finitely many service rates, \citet{badian2021optimal} identify a threshold policy for turning the server on, with the optimal service rate characterized by a series of monotone non-decreasing thresholds relative to the number of jobs in the system. In tandem queueing systems, \citet{yu2023optimal} derive an optimal threshold policy with respect to the abandonment cost. Even in cases where threshold policies are not optimal, they often serve as heuristics with demonstrated performance. For instance, \citet{pang2015logarithmic} propose a logarithmic safety staffing rule combined with a threshold control policy so that the server utilization is always close to one. Similarly, \citet{iravani1997two} develop a``Triple-Threshold" policy to approximate optimal behavior in the first stage with high accuracy. Finally, in the N-system with one flexible and one dedicated server, \citet{bell2001dynamic} demonstrate that a threshold control policy inspired by the Brownian control problem is asymptotically optimal in the heavy traffic limit.
\section{Problem Description and Main Results}\label{sec:description}
\begin{figure}[htbp]
\centering
    \includegraphics[scale=0.6]{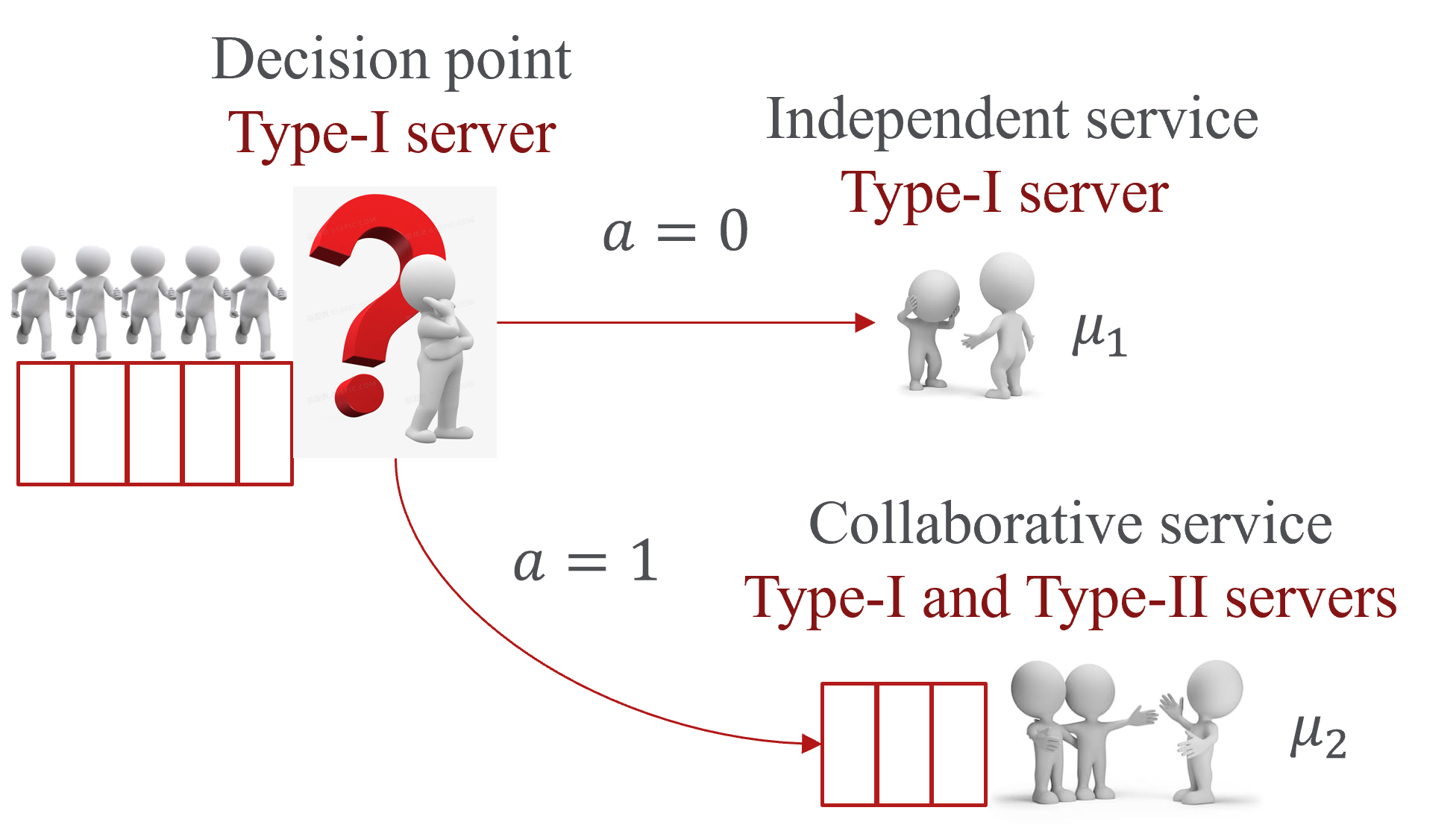}
\caption{System with independent and collaborative services.} \label{Fig:decision}
\end{figure}

We consider a clearing system queueing model equipped with two types of servers, as illustrated in Figure \ref{Fig:decision}. In this system, there are no external arrivals; instead, a fixed number of jobs are initially present and must be serviced before leaving the system. Upon completing a service with its previous job, a Type-I server immediately engages with a new job (if any) and decides whether to perform the service independently or collaboratively with a dedicated Type-II server.
The services are assumed to be non-preemptive and operate in a first-come, first-serve manner. Importantly, once a Type-I server begins serving a job, it must remain engaged until the job is completed, and cannot start another service until the current one concludes. 

Let the number of Type-I and Type-II servers be denoted as $C_1$ and $C_2$, respectively, where $C_2 < C_1$ for now. (The case where $C_2 \geq C_1$ is addressed as a corollary in our discussion.) All jobs require an exponential amount of work with rate 1. For clarity, we refer to the independent service station as Station 1, with a service rate $\mu_1$, and the collaborative service station as Station 2, with a service rate $\mu_2$. The state space of the stochastic process modeling this system is defined as follows:
\begin{align*}
    \X & := \Big\{(i,k,\ell)\in \Z_+^3\Big{|} i = 0, k+\ell < C_1 \text{ or } i \geq 0, k+\ell = C_1\Big\}, 
\end{align*}
where
\begin{itemize}
    \item $i$ is the number of jobs waiting in queue before the decision point,
    \item $k$ is the number of jobs currently receiving Station 1 service, and
    \item $\ell$ is the number of jobs at Station 2 (including those
        in service).
\end{itemize}
When a Type-I server completes a service at either Station 1 or 2, the control action is to decide what kind of service should be provided for the following job if one exists. Hence the set of states where a decision is made is 
\begin{align*}
    \X_{D} := \Big\{x = (i,k,\ell) \in \X | i \geq 1 \Big\}. 
\end{align*} 
In that sense, the action set is
\begin{align}
    \A(i,k,\ell) 
    & = \begin{cases}
        \left\{\big((1,a),(2,a)\big) | a \in \{0, 1\}\right\} & \text{if $x \in \X_D$,}\\
        \{0\} & \text{otherwise.}\\
    \end{cases}\label{def:action_set}
\end{align}

For each allowable action, note that $a \in \{0, 1\}$ is a binary variable indicating the server’s decision: $a = 0$ corresponds to providing service at Station 1, while $a = 1$ corresponds to serving at Station 2. The value $n=1$ (resp. $n=2$) at the first position indicates that a service at Station 1 (resp. 2) finishes before any at Station 2 (resp. 1). For example, the action $\big((1,0),(2,1)\big) \in \A(i,k,\ell)$ chooses to go to Station 1 if a Station 1 service finishes before any Station 2 services and chooses to go to Station 2 otherwise.

Let $\{Q_{0}(t), t\geq 0\}$, $\{Q_1(t), t\geq 0\}$ and $\{Q_2(t), t\geq 0\}$ denote the stochastic processes describing the number of jobs in queue before the decision point, at stations 1 and 2, respectively. As a measure of service quality, we assume the system is charged $h_0$, $h_1$, and $h_2$ per job, per unit of time for each.
For any control policy $\pi$, the total cost $T^{\pi}$ given initial state $x=(i,k,\ell) \in \X$ is
\begin{align*}
  T^{\pi}(x) &  = \int_{0}^{\infty} \Big(h_0 Q^{\pi}_{0}(t) +
    h_1 Q^{\pi}_1(t) + h_2 Q^{\pi}_2(t)\Big) dt,
\end{align*}
where $\big(Q^{\pi}_{0}(t), Q^{\pi}_1(0), Q^{\pi}_2(0)\big) =(i,k,\ell) = x$, and the dependence on the control policy has been added to each stochastic process.

Starting with any finite number of jobs initially, $T^{\pi} < \infty$ almost surely under any non-idling policy. We seek a control policy that minimizes the expectation of $T^{\pi}$ for any (fixed) initial state $x=(i,k,\ell)$, i.e.,
\begin{align*}
  \mathbb{E}\big[T^{\pi}(x)\big] &  = \mathbb{E}_{s} \Bigg[\int_{0}^{\infty}\Big(h_0 Q^{\pi}_{0}(t) +
    h_1 Q^{\pi}_1(t) + h_2 Q^{\pi}_2(t)\Big) dt \Bigg].
\end{align*}

Define $v(i,k,\ell)$, often referred to as the \emph{value function}, to be the optimal expected total cost incurred starting at state $(i,k,\ell) \in \X$ until the system clears all the jobs (so that $v(0,0,0) = 0$). Define $d(k,\ell) = k\mu_1 + \min\{\ell,C_2\}\mu_2$ and note that $d(k, \ell)$ is the overall system service rate at the state $(i,k,\ell)$. By the \emph{Principle of Optimality} (see
Section 4.3 in \citet{puterman2014markov}), for non-zero states, $v(i,k,\ell)$ satisfies the optimality equations as follows (c.f. Theorem 7.3.3 in \citet{puterman2014markov} for the optimality equations of discrete-time MDP).
\begin{enumerate}
    \item If $i = 0$, there are no decisions to make.
    \begin{align}
        \begin{split}
            v(0,k,\ell) & = \frac{k h_1 + \ell h_2 }{d(k,\ell)}
                + \bigg[ \frac{k \mu_1}{d(k,\ell)} v(0, k-1, \ell)  \\
                & \quad + \frac{\min\{\ell, C_2\}\mu_2}{d(k,\ell)} v(0, k, \ell-1) \bigg].
        \end{split}
    \end{align}
    \item If $i \geq 1$,
    \begin{align}
        \begin{split}
            v(i,k,\ell)
                &=\frac{i h_0+ k h_1 + \ell h_2 }{d(k,\ell)}
                + \bigg[ \frac{k \mu_1}{d(k,\ell)} \overbrace{\min\{v(i-1, k, \ell), v(i-1, k-1, \ell+1)\}}^{\text{where to work after a Station 1 service completion}} \label{eq:opt-clearing}\\
                & \quad + \frac{\min\{\ell, C_2\}\mu_2}{d(k,\ell)} \underbrace{\min\{v(i-1, k+1, \ell-1), v(i-1, k, \ell)\}}_{\text{where to work after a Station 2 service completion}} \bigg]
        \end{split}
    \end{align}
\end{enumerate}
For the ease of exposition, define the following sets of states
\begin{align*}
    \X_{diff} & := \Big\{x = (i,k,\ell) \in \X | k \geq 1 \Big\},\\ 
    \tilde{\X}_{diff}  &:= \Big\{x = (i,k,\ell) \in \X | k+\ell=C_1, k \geq 1 \Big\}. 
\end{align*} 
For any $x \in \X_{diff}$, define 
\begin{align}
    D(i,k,\ell) := v(i,k,\ell) - v(i,k-1,\ell+1), \label{def:D}
\end{align} 
where $v$ is a value function. Note that from the optimality equations \eqref{eq:opt-clearing}, for state $x=(i,k,\ell) \in \X_D$, an optimal decision of where to serve next depends on the sign of the difference of this type, regardless of whether a job from Station 1 or Station 2 is completed first. In particular, replacing $i$ with $i-1$, $D(i-1,k,\ell)$ is the difference considered when $n=1$, i.e., a service at Station 1 is completed first. Similarly, the same applies by replacing $(i,k,\ell)$ with $(i-1,k+1,\ell-1)$ in $D(i,k,\ell)$ when $n=2$. Consequently, if a threshold policy with threshold $T_D(k)$ (resp. $\tilde{T}_D(\ell)$) is defined by the first time $D(i,k,\ell)$ changes sign, say at $i = N(k)$ (resp.  $\tilde{N}(\ell)$), then $T_D(k) = N(k)+1$ (resp. $\tilde{T}_D(\ell) = \tilde{N}(\ell) + 1$) if $n=1$, and $T_D(k) = N(k+1)+1$ (resp. $\tilde{T}_D(\ell) = \tilde{N}(\ell-1) + 1$) if $n=2$. 

The set $\X_{diff}$ represents the states where $k\geq1$ so that $(i,k-1,\ell+1)$ remains in the state space. Furthermore, 
\begin{align*}
    \tilde{\X}_{diff} & = \X_{diff} \setminus \Big\{x = (0,k,\ell) \in \X | k+\ell < C_1, k \geq 1 \Big\} \subset \X_{diff},
\end{align*}
is the set of states where $k \geq 1$ and all the Type-I servers are busy. 

\subsection{Main results}\label{sec:main-results}
In this subsection, we present our main results regarding optimal policy structure. 
The condition $\frac{h_1}{\mu_1} > \frac{h_2}{\mu_2}$ implies that the expected cost to complete a single collaborative service is higher than independent service and likewise for the other direction for a single independent service versus a collaborative one. Recall $C_2 < C_1$ unless otherwise stated.

\begin{definition} \label{def:threshold_policy}
    We say a policy is of \textbf{collaborative (independent) threshold type with threshold $N$} if it chooses collaborative (independent) service if the number of jobs in queue is below $N$ and seeks independent (collaborative) service otherwise.
\end{definition}

\begin{theorem}\label{thm:policy1_i}
    Suppose $\frac{h_1}{\mu_1} \geq \frac{h_2}{\mu_2}$ and consider any $x =(i,k,\ell ) \in X_D$, the following holds.
    \begin{enumerate}
        \item If $\mu_1 \leq \mu_2$,
        \begin{enumerate}
            \item And a job at the independent station (Station 1) is completed first,
            \begin{enumerate}
                \item \label{state1.1:larger_mu2_smaller_ell}
                If $\ell < C_2$ (there is no Station 2 queue), an optimal policy always chooses beginning collaborative service next.
                \item \label{state1.1:larger_mu2_larger_ell}
                If $\ell \geq C_2$ (there is a Station 2 queue), there exists an optimal control policy of collaborative threshold type with threshold $N_2(\ell)+1$, where $N_2(\ell)$ is a finite threshold depending on $\ell$ only (or $k$ only since $k$ and $\ell$ satisfy $k+\ell = C_2$).
                \item If $\ell \geq C_2$ and in addition $\frac{h_1}{\mu_1} \leq \frac{\ell+1}{C_2}\frac{h_2}{\mu_2}$, regardless of the number of jobs in queue, an optimal policy always chooses independent service.  \label{state1.1:larger_mu2_larger_ell_action0}
            \end{enumerate}
            \item And a job at the collaborative station (Station 2) is completed first,
            \begin{enumerate}
                \item \label{state1.2:larger_mu2_smaller_ell}
                If $\ell \leq C_2$ (there is no Station 2 queue), an optimal policy always chooses the collaborative service.
                \item \label{state1.2:larger_mu2_larger_ell}
                If $\ell > C_2$ (there is a Station 2 queue), there exists an optimal control policy that is of collaborative threshold type with finite threshold $N_2(\ell-1)+1$ depending on $\ell$ only (or $k$ only).
                \item \label{state1.2:larger_mu2_larger_ell_action0}
                If $\ell > C_2$, and in addition $\frac{h_1}{\mu_1} \leq \frac{\ell}{C_2}\frac{h_2}{\mu_2}$, an optimal policy always chooses independent service.
            \end{enumerate}
        \end{enumerate}
        \item If $\mu_1 > \mu_2$,
        \begin{enumerate}
            \item
            And a job at the collaborative station is completed first,
            \begin{enumerate}
                \item \label{state1.1:larger_mu1} There exists an optimal control policy that is of collaborative threshold type with finite threshold $N_1(\ell)+1$ depending on $\ell$ only.
                \item \label{state1.1:larger_mu1_larger_ell_action0} If in addition to $\mu_1 > \mu_2$, $\frac{h_1}{\mu_1} \leq \frac{\ell+1}{C_2}\frac{h_2}{\mu_2}$, an optimal policy always chooses independent service.
            \end{enumerate}
            \item
            And a job at the collaborative station is completed first,
            \begin{enumerate}
                \item \label{state1.2:larger_mu1} There exists an optimal control policy that is of collaborative threshold type with finite threshold $N_1(\ell-1)+1$ depending on $\ell$ only.
                \item \label{state1.2:larger_mu1_larger_ell_action0} If in addition to $\mu_1 > \mu_2$, $\frac{h_1}{\mu_1} \leq \frac{\ell}{C_2}\frac{h_2}{\mu_2}$, an optimal policy always chooses independent service.
            \end{enumerate}
        \end{enumerate}
    \end{enumerate}
\end{theorem}

\begin{proof}
    The proof is divided into several parts. Statements \ref{state1.1:larger_mu1} and \ref{state1.2:larger_mu1} follow from Proposition \ref{prop:larger_mu1_dec}.
    Statements \ref{state1.1:larger_mu2_smaller_ell} and \ref{state1.2:larger_mu2_smaller_ell} are implied by Proposition \ref{prop:larger_mu2_always_action1}.
    Proposition \ref{prop:larger_mu2_dec} supports Statements \ref{state1.1:larger_mu2_larger_ell} and \ref{state1.2:larger_mu2_larger_ell}.
    Furthermore, Corollary \ref{cor:always_action0} confirms the remaining statements: \ref{state1.1:larger_mu2_larger_ell_action0}, \ref{state1.2:larger_mu2_larger_ell_action0}, \ref{state1.1:larger_mu1_larger_ell_action0}, and \ref{state1.2:larger_mu1_larger_ell_action0}.
\end{proof}

Notice that while the results in Theorem \ref{thm:policy1_i} hold on a case-by-case basis, there are some themes. For example, if service can begin immediately at Station 2, the rate of service is higher and the expected cost per service is lower, then the decision-maker should begin service at Station 2. This is akin to the existence of an optimal control of collaborative threshold type with threshold infinity. In addition, if $n=1$ (a Station 1 job is completed first) and $\frac{h_2}{\mu_2} \leq \frac{h_1}{\mu_1} \leq \frac{\ell+1}{C_2}\frac{h_2}{\mu_2}$ we have the existence of an optimal control of collaborative threshold type with threshold one ($N = 0$). The results of Theorem \ref{thm:policy1_i} therefore characterize the threshold policies, specifying when those thresholds are one (corresponding to $N=0$), infinite or take other finite values. 
When the thresholds are finite (including the case when $N=0$), Theorem \ref{thm:policy1_i} indicates that a larger initial queue discourages collaborative service.
The results are also depicted in Table \ref{table:policy1}.
\begin{theorem}\label{thm:policy2_i}
    Suppose $\frac{h_1}{\mu_1} \leq \frac{h_2}{\mu_2}$ and consider any $x =(i,k,\ell ) \in X_D$, the following holds.
    \begin{enumerate}
        \item If $\mu_1 < \mu_2$,
        \begin{enumerate}
            \item And a job at the independent station is completed first,
            \begin{enumerate}
                \item \label{state2.1:larger_mu2_smaller_ell}
                If $\ell < C_2$, there exists an optimal control policy that is of collaborative threshold type with finite threshold $N_3(\ell)+1$ depending on $\ell$ only.
                \item \label{state2.1:larger_mu2_larger_ell}
                If $\ell \geq C_2$, an optimal policy always chooses the independent service.
            \end{enumerate}
            \item And a job at the collaborative station is completed first,
            \begin{enumerate}
                \item \label{state2.2:larger_mu2_smaller_ell}
                If $\ell \leq C_2$, there exists an optimal control policy that is of collaborative threshold type with finite threshold $N_3(\ell-1)+1$ depending on $\ell$ only.
                \item \label{state2.2:larger_mu2_larger_ell}
                If $\ell > C_2$, an optimal policy always chooses the independent service.
            \end{enumerate}
        \end{enumerate}
        \item \label{state2:larger_mu1}
        If $\mu_1 \geq \mu_2$,
        an optimal policy always chooses the independent service.
    \end{enumerate}
\end{theorem}

\begin{proof}
    The proof is structured in several parts. Proposition \ref{prop:larger_mu2_always_action0} supports Statements \ref{state2.1:larger_mu2_larger_ell} and \ref{state2.2:larger_mu2_larger_ell}.
    Proposition \ref{prop:larger_mu2_inc} implies Statements \ref{state2.1:larger_mu2_smaller_ell} and \ref{state2.2:larger_mu2_smaller_ell}.
    Furthermore, Corollary \ref{cor:always_action0} supports the remaining statement: \ref{state2:larger_mu1}.
\end{proof}

The results of Theorem \ref{thm:policy2_i} cover the case when the expected cost of collaborative service are higher than independent service. However, they can each be stated in terms of threshold policies. Since there are no queueing concerns with independent service, if the expected cost is lower and either the rate of service higher or one would have to queue for collaborative service, the decision-maker should always choose independent service ($N=\infty$). In other cases, there exists a finite threshold, suggesting that a larger initial queue encourages collaborative service. See also Table \ref{table:policy2}.
\begin{table} [htbp]
\caption{Policy structure for any state $(i,k,\ell)\in \X_{D}$ if $\frac{h_1}{\mu_1} \geq \frac{h_2}{\mu_2}$}\label{table:policy1}
\centering
\begin{tabular}{|c|c|c|c|c|}
\hline

\multicolumn{1}{|l|}{} & \multicolumn{1}{c|}{Service that finishes first} & \multicolumn{1}{c|}{\textbf{$\ell< C_2$}} & \multicolumn{1}{c|}{\textbf{$\ell = C_2$}} & \multicolumn{1}{c|}{\textbf{$\ell > C_2$}} \\ \hline
\multirow{2}{*}{\textbf{$\mu_1 \leq \mu_2$}}   & Station 1    & \multicolumn{1}{c|}{$a=1$}            & \multicolumn{2}{c|}{$a=\mathbb{I}\{i< N_2(\ell)+1\}$} \\
\cline{3-5}
    & Station 2 &   \multicolumn{2}{c|}{$a=1$} &  \text{$a=\mathbb{I}\{i< N_2(\ell-1)+1\}$}\\ \hline

\multirow{2}{*}{\textbf{$\mu_1 > \mu_2$}}   & Station 1  &  \multicolumn{3}{c|}{$a=\mathbb{I}\{i< N_1(\ell)+1\}$} \\ 
    & Station 2 & \multicolumn{3}{c|}{$a=\mathbb{I}\{i< N_1(\ell-1)+1\}$} \\ \hline    

\end{tabular}
\end{table}

\begin{table} [htbp]
\caption{Policy structure for any state $(i,k,\ell)\in \X_{D}$ if $\frac{h_1}{\mu_1} \leq \frac{h_2}{\mu_2}$}\label{table:policy2}
\centering
\begin{tabular}{|c|c|c|c|c|}
\hline

\multicolumn{1}{|l|}{} & \multicolumn{1}{c|}{Service that finishes first} & \multicolumn{1}{c|}{\textbf{$\ell< C_2$}} & \multicolumn{1}{c|}{\textbf{$\ell = C_2$}} & \multicolumn{1}{c|}{\textbf{$\ell > C_2$}} \\ \hline
\multirow{2}{*}{\textbf{$\mu_1 < \mu_2$}}   & Station 1    & $a=\mathbb{I}\{i \geq N_3(\ell)+1\}$            & \multicolumn{2}{c|}{\text{$a=0$}}  \\ 
\cline{3-5}
    & Station 2 & \multicolumn{2}{c|}{$a=\mathbb{I}\{i \geq N_3(\ell-1)+1\}$}      & \text{$a=0$} \\ \hline
\multirow{2}{*}{\textbf{$\mu_1 \geq \mu_2$}}   & Station 1  &  \multicolumn{3}{c|}{\multirow{2}{*}{$a=0$}}\\ 
    & Station 2 & \multicolumn{3}{c|}{\multirow{2}{*}{}} \\ \hline 
\end{tabular}
\end{table}

In addition to the monotone threshold structure in the number of jobs in the queue, we also have the monotonicity result with respect to the job counts at the two service stations.
This result suggests that an optimal control favors collaboration when Station 1 has more jobs, while Station 2 has fewer.
\begin{theorem}\label{thm:policy_k_ell}
    Consider any $x=(i,k,\ell) \in \X_D$ and $x' = (i,k+1,\ell-1) \in X_D$. The following holds
    \begin{enumerate}
        \item When an independent service job finishes first, an optimal decision at $x=(i,k+1,\ell-1)$ chooses the collaborative service if an optimal decision at $(i,k,\ell)$ does.
        \item A similar conclusion holds when a collaborative service job finishes first.
    \end{enumerate}
\end{theorem}
\begin{proof}
    This result is a direct application of Proposition \ref{prop:mono_ell} (see next subsection).
\end{proof}
We conclude this section with a special case to Theorems \ref{thm:policy1_i}--\ref{thm:policy_k_ell}, focusing on the scenario without queueing effects. See Table \ref{table:policy_enough_C2}.
\begin{theorem}
    Suppose $C_2 \geq C_1$. Consider any $x =(i,k,\ell ) \in X_D$, the following holds:
    \begin{enumerate}
        \item If $\frac{h_1}{\mu_2} \geq \frac{h_2}{\mu_2}$,
        \begin{enumerate}
            \item If $\mu_1 \leq \mu_2$, an optimal policy always chooses the collaborative service.
            \item If $\mu_1 > \mu_2$, a similar result holds as in Theorem \ref{thm:policy1_i} with the same assumption.
        \end{enumerate}
        \item If $\frac{h_1}{\mu_2} \leq \frac{h_2}{\mu_2}$,
        \begin{enumerate}
            \item If $\mu_1 < \mu_2$, a similar result holds as in Theorem \ref{thm:policy2_i} with the same assumption.
            \item If $\mu_1 \geq \mu_2$, an optimal policy always chooses the independence service.
        \end{enumerate}
        \item Theorem \ref{thm:policy_k_ell} still holds.
    \end{enumerate}
\end{theorem}
\begin{proof}
    Each statement follows from the corresponding statement in Corollary \ref{cor:enough_Type_II} under the same conditions.
\end{proof}
\begin{table} [htbp]
\caption{Policy structure for any state $(i,k,\ell)\in \X_{D}$ if $C_2 \geq C_1$}\label{table:policy_enough_C2}
\centering
\begin{tabular}{|c|c|c|c|}
\hline

\multicolumn{1}{|l|}{} & \multicolumn{1}{c|}{Service that finishes first} & \multicolumn{1}{c|}{\textbf{$\frac{h_1}{\mu_1} > \frac{h_2}{\mu_2}$}} & \multicolumn{1}{c|}{\textbf{$\frac{h_1}{\mu_1} \leq \frac{h_2}{\mu_2}$}} \\ \hline
\multirow{2}{*}{\textbf{$\mu_1 < \mu_2$}}   & Station 1    & \multirow{4}{*}{$a=1$}            & \text{$a=\mathbb{I}\{i> N_3(\ell)+1\}$} \\
    & Station 2 &    &  \text{$a=\mathbb{I}\{i> N_3(\ell-1)+1\}$}\\ \cline{1-2} \cline{4-4}

\multirow{2}{*}{\textbf{$\mu_1 = \mu_2$}}   & Station 1    &             & \multirow{4}{*}{$a=0$} \\
    & Station 2 &    &  \\ \cline{1-3}
    
\multirow{2}{*}{\textbf{$\mu_1 > \mu_2$}}   & Station 1  &  \text{$a=\mathbb{I}\{i\leq N_1(\ell)+1\}$} & \\ 
    & Station 2 &  \text{$a=\mathbb{I}\{i\leq N_1(\ell-1)+1\}$} &  \\ \hline    

\end{tabular}
\end{table}

Computation of each threshold can be challenging. In Section \ref{sec:heuristics} we provide simple heuristics to estimate thresholds with demonstrated performance.

\subsection{Supporting results} \label{sec:support-results}
This subsection examines the behavior of the difference $D(i,k,\ell)$ defined in \eqref{def:D} as $i$ (the number of jobs in the queue) varies. As shown in Section \ref{sec:main-results}, the results presented in Propositions \ref{prop:larger_mu1_dec}--\ref{prop:larger_mu2_inc} and Corollary \ref{cor:always_action0} support the proofs of Theorem \ref{thm:policy1_i} and Theorem \ref{thm:policy2_i}. See also Table \ref{table:sign1_D} for the case $\frac{h_1}{\mu_1} > \frac{h_2}{\mu_2}$ and Table \ref{table:sign2_D} for $\frac{h_1}{\mu_1} \leq \frac{h_2}{\mu_2}$. The case where $C_2 \geq C_1$ is summarized as a special scenario in Corollary \ref{cor:enough_Type_II} and Table \ref{table:sign_D_enough_C2}. 

Moreover, $D(i,k,\ell)$ is monotone diagonally in $k$ and $\ell$ (the job counts in Stations 1 and 2). 
That is $D(i,k+1,\ell-1) \geq D(i,k,\ell)$ for $(i,k+1,\ell-1) \in \X_{diff}$ and $ (i,k,\ell)\in \X_{diff}$. See Proposition \ref{prop:mono_ell} for details, which supports Theorem \ref{thm:policy_k_ell}.
Proofs of these propositions and corollaries are provided in Appendix \ref{sec:support-results}.
\begin{table} [htbp]
\caption{Sign of the difference $D(i,k,\ell)$ where $x\in \tilde{\X}_{diff}$ if $\frac{h_1}{\mu_1} > \frac{h_2}{\mu_2}$} \label{table:sign1_D}
\centering
\begin{tabular}{|c|c|c|}
\hline

\multicolumn{1}{|l|}{} & \multicolumn{1}{c|}{\textbf{$\ell < C_2$}} & \multicolumn{1}{c|}{\textbf{$\ell\geq C_2$}} \\ \hline

\textbf{$\mu_1 \leq \mu_2$}    & \text{$D > 0$ (Prop. \ref{prop:larger_mu2_always_action1})}               & \text{$ 
\begin{cases} 
D > 0 & \text{if } i < N_2(\ell)\\ 
D \leq 0    & \text{otherwise}
\end{cases}$  (Prop. \ref{prop:larger_mu2_dec})}  \\
\hline

\textbf{$\mu_1 > \mu_2$}    & \multicolumn{2}{c|}{$ 
\begin{cases} 
D > 0    & \text{if } i < N_1(\ell) \\
D \leq 0 & \text{otherwise}

\end{cases}$  (Prop. \ref{prop:larger_mu1_dec})}     \\
\hline
\end{tabular}
\end{table}

\begin{table} [htbp]
\caption{Sign of the difference $D(i,k,\ell)$ where $x\in \tilde{\X}_{diff}$ if $\frac{h_1}{\mu_1} \leq \frac{h_2}{\mu_2}$} \label{table:sign2_D}
\centering
\begin{tabular}{|c|c|c|}
\hline

\multicolumn{1}{|l|}{} & \multicolumn{1}{c|}{\textbf{$\ell < C_2$}} & \multicolumn{1}{c|}{\textbf{$\ell\geq C_2$}} \\ \hline

\textbf{$\mu_1 < \mu_2$}    & \text{$ 
\begin{cases} 
D \leq 0    & \text{if } i < N_3(\ell)\\
D > 0 & \text{otherwise}
\end{cases}$  (Prop. \ref{prop:larger_mu2_inc})}               & \text{$D\leq0$  (Prop. \ref{prop:larger_mu2_always_action0})}  \\
\hline

\textbf{$\mu_1 \geq \mu_2$}    & \multicolumn{2}{c|}{$D \leq0$  (Cor. \ref{cor:always_action0})}     \\
\hline
\end{tabular}
\end{table}

\begin{table} [htbp]
\caption{Sign of the difference $D(i,k,\ell)$ where $x\in \tilde{\X}_{diff}$ if $C_2 \geq C_1$ (Cor. \ref{cor:enough_Type_II})} \label{table:sign_D_enough_C2}
\centering
\begin{tabular}{|c|c|c|}
\hline

\multicolumn{1}{|l|}{} & \multicolumn{1}{c|}{\textbf{$\frac{h_1}{\mu_1} > \frac{h_2}{\mu_2}$}} & \multicolumn{1}{c|}{\textbf{$\frac{h_1}{\mu_1} \leq \frac{h_2}{\mu_2}$}} \\ \hline

\textbf{$\mu_1 < \mu_2$}    &  \multirow{3}{*}{$D > 0$}               &  \text{$ 
\begin{cases} 
D > 0 & \text{if } i< N_3(\ell)\\ 
D \leq 0    & \text{otherwise} 
\end{cases}$}  \\
\cline{1-1} \cline{3-3}

\textbf{$\mu_1 = \mu_2$}    &               & \multirow{3}{*}{$D\leq0$}  \\
\cline{1-2}

\textbf{$\mu_1 > \mu_2$}    & \text{$  
\begin{cases} 
D > 0 & \text{if } i< N_1(\ell)\\ 
D \leq 0    & \text{otherwise} 
\end{cases}$}               &   \\
\hline
\end{tabular}
\end{table}

\begin{proposition}\label{prop:larger_mu1_dec}
    Assume $\mu_1 \geq \mu_2$. (No need to assume $\frac{h_1}{\mu_1} \geq \frac{h_2}{\mu_2}$.) Consider any $x = (i,k,\ell) \in \tilde{\X}_{diff}$.
    \begin{enumerate}
        \item The following inequality holds for all $i \geq 0$,
            \begin{align}
             \lefteqn{D(i,k,\ell) - D(i+1,k,\ell)}\nonumber\\
                &\quad = v(i,k,\ell) - v(i,k-1,\ell+1)- \Big[v(i+1,k,\ell) - v(i+1,k-1,\ell+1) \Big]\geq 0. \label{eq:larger_mu1_dec}
            \end{align}
        In other words, the difference $D(i,k,\ell)$ is non-increasing in $i\geq0$.
        \item In addition, there exists a finite threshold $N_1(\ell)$ such that the difference $D(i,k,\ell)$ where $\ell \geq C_2$ is non-positive if and only if $i \geq N_1(\ell)$ and in fact becomes negative for all $i$ large enough. If we further assume $\mu_1 > \mu_2$, a similar result holds for $\ell < C_2$.
    \end{enumerate}
\end{proposition}

\begin{proposition} \label{prop:larger_mu2_always_action1}
    Assume $\mu_2 \geq \mu_1$ and $\frac{h_1}{\mu_1} \geq \frac{h_2}{\mu_2}$. For any $x=(i,k,\ell) \in \X_{diff}$ where $\ell < C_2$, the following inequality holds for all $i \geq 0$,
    \begin{align}
    D(i,k,\ell) = v(i,k,\ell) - v(i,k-1,\ell+1) \geq 0. \label{eq:larger_mu2_always_action1}
    \end{align} 
    Moreover, if $\frac{h_1}{\mu_1} > \frac{h_2}{\mu_2}$ holds strictly, then inequality \eqref{eq:larger_mu2_always_action1} is also strict, meaning that $D(i,k,\ell) > 0$.
\end{proposition}

\begin{proposition}\label{prop:larger_mu2_dec}
    Assume $\mu_2 \geq \mu_1$ and $\frac{h_1}{\mu_1} \geq \frac{h_2}{\mu_2}$. Consider any $x = (i,k,\ell) \in \tilde{\X}_{diff}$ where $\ell \geq C_2$. 
    \begin{enumerate} 
        \item \label{state:larger_mu2_dec}
        The following inequality holds for all $i \geq 0$,
    \begin{align}
    \lefteqn{D(i,k,\ell) - D(i+1,k,\ell)}\nonumber\\
        &\quad = v(i,k,\ell) - v(i,k-1,\ell+1)- \Big[v(i+1,k,\ell) - v(i+1,k-1,\ell+1) \Big]\geq 0.   \label{eq:larger_mu2_dec}
    \end{align}
    In other words, the difference $D(i,k,\ell)$ is non-increasing in $i$.
        \item \label{state:larger_mu2_action0_large_i}
        In addition, there exists a finite threshold $N_2(\ell)$ such that the difference $D(i,k,\ell)$ is non-positive if and only if $i \geq N_2(\ell)$ and in fact becomes negative for all $i$ large enough.
    \end{enumerate}
\end{proposition}

\begin{proposition} \label{prop:larger_mu2_always_action0}
    Assume $\mu_2 \geq \mu_1$ and $\frac{h_1}{\mu_1} \leq \frac{h_2}{\mu_2}$. For any $x = (i,k,\ell) \in \X_{diff}$ where $\ell \geq C_2$, 
    \begin{align}
    D(i,k,\ell) = v(i,k,\ell) - v(i,k-1,\ell+1) \leq 0. \label{eq:larger_mu2_always_action0}
\end{align} 
\end{proposition}
\begin{proposition}\label{prop:larger_mu2_inc}
    Assume $\mu_2 \geq \mu_1$ and $\frac{h_1}{\mu_1} \leq \frac{h_2}{\mu_2}$. Consider any $x = (i,k,\ell) \in \tilde{\X}_{diff}$ where $\ell < C_2$.
    \begin{enumerate}
        \item The following inequality holds for all $i \geq 0$,
        \begin{align}
            \lefteqn{D(i,k,\ell) - D(i+1,k,\ell)}\nonumber\\
            &\quad = v(i,k,\ell) - v(i,k-1,\ell+1)- \Big[v(i+1,k,\ell) - v(i+1,k-1,\ell+1) \Big] \leq 0.  \label{eq:larger_mu2_inc}
    \end{align}
    In other words, the difference $D(i,k,\ell)$ is non-decreasing in $i$.
    \item If further we assume $\mu_2>\mu_1$, there exists a finite threshold $N_3(\ell)$ such that the difference $D(i,k,\ell)$ is non-positive if and only if $i < N_3(\ell)$ and becomes positive for all $i$ large enough.
    \end{enumerate}
\end{proposition}

\begin{corollary}\label{cor:always_action0}
    Consider any state $x =(i,k,\ell) \in \X_{diff}$. Suppose either of the following conditions holds:
    \begin{enumerate}
        \item If $\mu_1 \geq \mu_2$,
        \item If $\mu_1 \leq \mu_2$ and $\ell \geq C_2$.
    \end{enumerate}
    Then $\frac{h_1}{\mu_1} \leq \frac{\max\{\ell+1,C_2\}}{C_2}\frac{h_2}{\mu_2}$ implies $D(i,k,\ell) \leq 0$ for all $i \geq 0$.
\end{corollary}

\begin{proposition}\label{prop:mono_ell}
    Consider any $x = (i,k,\ell) \in \X_{diff}$ and $x'=(i,k+1,\ell-1) \in \X_{diff}$. The following inequality holds for all $i\geq0$:
    \begin{align}
     \lefteqn{D(i,k+1,\ell-1) - D(i,k,\ell)}& \nonumber \\
        & \quad = v(i,k+1,\ell-1) - v(i,k,\ell) - \Big[v(i,k,\ell) - v(i,k-1,\ell+1)\Big] 
        \geq 0. 
        \label{eq:mono_ell}
    \end{align}
\end{proposition}

\begin{corollary}\label{cor:enough_Type_II}
    Assume $C_2 \geq C_1$. Consider any $x = (i,k,\ell)\in \X_{diff}$. The following holds
    \begin{enumerate}
        \item If $\frac{h_1}{\mu_2} > \frac{h_2}{\mu_2}$,
        \begin{enumerate}
            \item $\mu_1 \leq \mu_2$, $D(i,k,\ell) > 0$ (similarly as Proposition \ref{prop:larger_mu2_always_action1});
            \item $\mu_1 > \mu_2$, a similar result holds as in Proposition \ref{prop:larger_mu1_dec}.
        \end{enumerate}
        \item If $\frac{h_1}{\mu_2} \leq \frac{h_2}{\mu_2}$,
        \begin{enumerate}
            \item $\mu_1 < \mu_2$, a similar result holds as in Proposition \ref{prop:larger_mu2_inc}.
            \item $\mu_1 \geq \mu_2$, $D(i,k,\ell) \leq 0$ (similarly as Corollary \ref{cor:always_action0});
        \end{enumerate}
        \item When $x'=(i,k+1,\ell-1) \in \X_{diff}$, a similar result holds as in Proposition \ref{prop:mono_ell}: $D(i,k+1,\ell-1) - D(i,k,\ell)\geq 0$ for all $i$.
    \end{enumerate}
\end{corollary}

\section{Heuristics Design} \label{sec:heuristics}
This section aims at designing simple yet robust heuristics that effectively approximate an optimal policy. As established, determining an optimal control reduces to analyzing the sign of $D(i,k,\ell)$ for $x=(i,k,\ell) \in \tilde{\X}_{diff}$. Specifically:
\begin{itemize}
    \item When $\frac{h_1}{\mu_1} > \frac{h_2}{\mu_2}$, the sign of $D$, given in Table \ref{table:sign1_D}, implies the existence of $N(k)$ (dependent on $k$) such that $D(i,k,\ell) \leq 0$ if and only if $i \geq N(k)$, where $N(k)$ is defined as
    \begin{align}
        N(k) := \begin{cases}
            N_1(\ell) \quad &\text{if } \mu_1 >\mu_2,\\
            N_2(\ell) \quad &\text{if } \mu_1 \leq \mu_2 \text{ and } \ell = C_1-k \geq C_2,\\
            \infty \quad &\text{if } \mu_1 \leq \mu_2 \text{ and } \ell = C_1 - k < C_2.
        \end{cases} \label{eq:threshold_N_k}
    \end{align} 
    Here, we equivalently change the index from $\ell$ to $k$, using the relationship $k+\ell=C_1$ to maintain consistency with future definitions and arguments. 
    \item When $\frac{h_1}{\mu_1} \leq \frac{h_2}{\mu_2}$, the sign of $D$, summarized in Table \ref{table:sign2_D}, indicates that $D(i,k,\ell) > 0$ if and only if $i \geq \tilde{N}(\ell)$ for some $\tilde{N}(\ell)$, where $\tilde{N}(\ell)$ is defined as
    \begin{align}
        \tilde{N}(\ell) := \begin{cases}
            N_3(\ell) \quad &\text{if } \mu_1 < \mu_2 \text{ and } \ell < C_2,\\
            \infty \quad &\text{otherwise.}
        \end{cases} \label{eq:threshold_N_ell}
    \end{align} 
\end{itemize}

A key observation from Tables \ref{table:sign1_D} and \ref{table:sign2_D} is that, regardless of parameter configuration, $D(i,k,\ell)$ changes sign at most once as $i$ increases.
Moreover, the direction of this sign change (from positive to non-positive or vice versa) depends on whether $\frac{h_1}{\mu_1} > \frac{h_2}{\mu_2}$. 
\subsection{Motivation and definition of the linear approximation}
Leveraging the structural property of at most one sign change in $D$, for any $x = (i,k,\ell) \in \tilde{\X}_{diff}$, we define $H(i,k,\ell)$, an affine approximation (in $i$) of $D(i,k,\ell)$, as  
\begin{align}
    H(i,k,\ell) := i \cdot f(k) + g(k), \label{eq:f_and_g}
\end{align}
where both $f(k)$ and $g(k)$ depend only on $k$ (or equivalently $\ell$) but not on $i$. We estimate the sign change of $D$ using that of $H$. 

As shown in Tables \ref{table:sign1_D} and \ref{table:sign2_D}, there are regions where the sign of $D(i,k,\ell)$ is known to be independent of $i$, allowing $H(i,k,\ell)$ to be designed accordingly with the same sign. For example, if $D(i,k,\ell) > 0$ for all $i \geq 0$ in a given region, any $f(k) \geq 0$ and $g(k) > 0$ will suffice under the same parameter configuration. In particular, when $f(k)=0$, the function reduces to $H(i,k,\ell) = g(k)$, which is a constant in $i$ and thus does not change sign. Similar logic applies when $D \leq 0$ (for all $i$). We therefore focus on designing $H$ in parameter regions where $D(i,k,\ell)$ changes sign as $i$ increases, requiring a finite threshold --- i.e., where $f(k) \neq 0$ is necessary.

A specific form of $H$ is proposed in Definition \ref{def:H} below, where the slope $f(k)$ and the intercept $g(k)$ in \eqref{eq:f_and_g} are readily identifiable from the expression in \eqref{eq:def_H}.
For instance, when $\ell = C_1-k < C_2$, we have $f(k) = c$ and $g(k) = b$.
\begin{definition} \label{def:H}
    Consider any $x = (i,k,\ell)\in \tilde{\X}_{diff}$, so that $k+\ell=C_1$. Let
     \begin{align}
        H(i,k,\ell) &:= 
        \left\{
        \begin{array}{lcl}
            -1 & \text{if} & \ell \geq C_2 \text{ and } \frac{h_1}{\mu_1} \leq \frac{\ell+1}{C_2}\frac{h_2}{\mu_2},\\
           (i-y_k)c'+b' & \text{if} & \ell \geq C_2 \text{ and } \frac{h_1}{\mu_1} > \frac{\ell+1}{C_2}\frac{h_2}{\mu_2}, \\ 
          i c+b& \text{if} & \ell < C_2, \\
        \end{array}
        \right. \label{eq:def_H}
    \end{align}where
    \begin{align}
        c&:=\frac{h_0}{C_1}\left(\frac{1}{\mu_1} - \frac{1}{\mu_2}\right), &   b&:=\frac{h_1}{\mu_1} - \frac{h_2}{\mu_2}, \label{def:b_and_c}\\
        c'&:=-\frac{h_0}{C_2\mu_2}, & b'&:=\frac{h_1-h_2}{\mu_1} - \frac{C_1h_2}{C_2\mu_2},\label{def:b'_and_c'}
    \end{align}
    \begin{align}
        y_k &:= (k-1)+\min\{\ell,C_2\} m 
        = \left\{
        \begin{array}{lcl}
            (k-1) + \ell m, & \text{if} & \ell < C_2,\\
            (k-1) + C_2 m, & \text{if} & \ell \geq C_2,
        \end{array}
        \right.\label{def:sequence_y}
    \end{align}
    and $m$ is the ratio of the service rates defined as
\begin{align}
    m:=\frac{\mu_2}{\mu_1}, \label{def:ratio_of_rates}
\end{align}
\end{definition}

The definition of $H$ in \eqref{eq:def_H} is well-motivated by the underlying structure of $D$ and merits a case-by-case explanation to elucidate its construction.
\begin{itemize}
    \item When $\frac{h_1}{\mu_1} > \frac{h_2}{\mu_2}$ (see Table \ref{table:sign1_D}),
    \begin{enumerate} [label= \textbf{Case} \arabic*:, leftmargin=3\parindent]
        \item \label{case:queue_in_collab_dec}
        For $\ell \geq C_2$, $D(i,k,\ell)$ is not only monotone, non-increasing in $i$, but also eventually becomes negative for large $i$, implying a finite $N(k)$ in \eqref{eq:threshold_N_k}. (see Proposition \ref{prop:larger_mu1_dec} if $\mu_1>\mu_2$ or Proposition \ref{prop:larger_mu2_dec} if $\mu_1\leq\mu_2$). Given this, the linear function $H(i,k,\ell)$ should also be designed to be monotone, non-increasing in $i$ (i.e., $f(k)\leq 0$). 
        \begin{enumerate} [label= \textbf{Subcase} \alph*:, leftmargin=3\parindent]
            \item If $\frac{h_1}{\mu_1} \leq \frac{\ell+1}{C_2}\frac{h_2}{\mu_2}$ (first case of \eqref{eq:def_H}), in particular, we enforce $f(k) = 0$ and $g(k) = -1$, so that $H(i,k,\ell) = -1 \leq 0$, aligning with $D(i,k,\ell) \leq 0$ for all $i$ (Corollary \ref{cor:always_action0}).
            \item If $\frac{h_1}{\mu_1} > \frac{\ell+1}{C_2}\frac{h_2}{\mu_2}$ (second case of \eqref{eq:def_H}), note $f(k) = c' < 0$ by \eqref{def:b'_and_c'}.
        \end{enumerate}
        \item For $\ell < C_2$ (last case in \eqref{eq:def_H}), we let $f(k) = c$ and $g(k) = b > 0$.
        \begin{enumerate} [label= \textbf{Subcase} \alph*:, leftmargin=3\parindent]
            \item If $\mu_1 > \mu_2$, a similar result as \ref{case:queue_in_collab_dec} when $\ell \geq C_2$ (Proposition \ref{prop:larger_mu1_dec}) holds. Moreover, $N(k)>0$ in \eqref{eq:threshold_N_k} since $D(0,k,\ell) = b > 0$ by \eqref{eq:diff_bdry}. Correspondingly, we set $f(k) = c <0$ and $H(0,k,\ell) = g(k) = b =D(0,k,\ell)>0$, so that $H$ becomes non-positive at a positive finite value of $i$.
            \item If $\mu_1\leq\mu_2$, then $c \geq 0$ and $b > 0$, ensuring $H(i,k,\ell)=ci+b > 0$, aligning with $D(i,k,\ell) > 0$ (Proposition \ref{prop:larger_mu2_always_action1}).
        \end{enumerate}
    \end{enumerate}
    \item When $\frac{h_1}{\mu_1} \leq \frac{h_2}{\mu_2}$ (see  Table \ref{table:sign2_D}), the design of $H$ follows Definition \ref{def:H} as well.
    \begin{enumerate} [label= \textbf{Case} \arabic*:, leftmargin=3\parindent]
        \item For $\ell \geq C_2$, we always have $H(i,k,\ell) = -1 \leq 0$ since $\frac{h_1}{\mu_1} \leq \frac{h_2}{\mu_2} < \frac{\ell+1}{C_2}\frac{h_2}{\mu_2}$ (first case of \eqref{eq:def_H}). This corresponds to $D(i,k,\ell)\leq 0$ for all $i$ (Proposition \ref{prop:larger_mu2_always_action0}). 
        \item For $\ell < C_2$ (last case in \eqref{eq:def_H}), we let $f(k) = c$ and $g(k) = b \leq 0$.
        \begin{enumerate} [label= \textbf{Subcase} \alph*:, leftmargin=3\parindent]
            \item If $\mu_1 < \mu_2$, $D(i,k,\ell)$ is monotone non-decreasing and eventually becomes positive as $i$ increases (Proposition \ref{prop:larger_mu2_inc}). This implies a finite $N(\ell)$ in \eqref{eq:threshold_N_ell}. Furthermore, $N(\ell) > 0$ since $D(0,k,\ell) = b \leq 0$ by \eqref{eq:diff_bdry}.  Accordingly, the slope $f(k) = c > 0$ and the intercept $g(k) = H(0,k,\ell) = b = D(0,k,\ell) \leq 0$, so that $H$ becomes positive at a non-negative number of $i$.
            \item If $\mu_1 \geq \mu_2$, $H(i,k,\ell) = ci+b \leq 0$, since $c \leq 0$ and $b \leq 0$ as given in \eqref{def:b_and_c}, preserving consistency with $D(i,k,\ell) \leq 0$ (Corollary \ref{cor:always_action0}).
        \end{enumerate}
    \end{enumerate}
\end{itemize}

\subsection{Using \texorpdfstring{$H$}{H} to approximate decisions with \texorpdfstring{$D$}{D}}
Informally, the \textbf{actual} and \textbf{estimated integer thresholds}, defined next are the smallest non-negative integer at which $D(i,k,\ell)$ and $H(i,k,\ell)$, as functions of $i$, respectively, change sign from positive to non-positive or vice versa. 
These definitions --- formalized below in Definition \ref{def:thresholds_lowcost} for $\frac{h_1}{\mu_1} > \frac{h_2}{\mu_2}$ and in Definition \ref{def:thresholds_highcost} for $\frac{h_1}{\mu_1} \leq \frac{h_2}{\mu_2}$ --- enable us to derive bounds on these thresholds for subsequent analysis.

Building on our construction of $H$ in Definition \ref{def:H}, for each fixed $k$ or $\ell$, the estimated integer thresholds can be explicitly expressed in terms of system parameters, suggesting the sign of $H$ as summarized in Table \ref{table:sign1_H} for $\frac{h_1}{\mu_1} > \frac{h_2}{\mu_2}$ and Table \ref{table:sign2_H} for $\frac{h_1}{\mu_1} > \frac{h_2}{\mu_2}$, where $R_1$ and $R_2(k)$ are defined later in \eqref{eq:thresh_no_queue_in_collab} and \eqref{eq:thresh_queue_in_collab}, respectively. 
In light of Corollary \ref{cor:enough_Type_II}, $C_2 \geq C_1$ can be interpreted as a special case of these results. The sign of $H$, along with the corresponding integer thresholds in this case, are summarized in Table \ref{table:sign_H_enough_C2}.

 \begin{table} [htbp]
\caption{Sign of the approximation $H(i,k,\ell)$ where $x\in \tilde{\X}_{diff}$ if $\frac{h_1}{\mu_1} > \frac{h_2}{\mu_2}$} \label{table:sign1_H}
\centering
\begin{tabular}{|c|c|c|}
\hline

\multicolumn{1}{|l|}{} & \multicolumn{1}{c|}{\textbf{$\ell < C_2$}} & \multicolumn{1}{c|}{\textbf{$\ell\geq C_2$}} \\ \hline

\textbf{$\mu_1 \leq \mu_2$}    & \text{$H > 0$ }               & \multirow{3}{*}{
$\begin{cases} 
H> 0 & \text{if } \frac{h_1}{\mu_1} > \frac{\ell+1}{C_2}\frac{h_2}{\mu_2}, \text{ and } i < \max\left\{\left\lceil R_2(k)\right\rceil, 0\right\}\\ 
H\leq 0    & \text{otherwise} 
\end{cases}$ }  \\
\cline{1-2}

\textbf{$\mu_1 > \mu_2$}    &    {$ 
\begin{cases} 
H> 0 & \text{if } i < \left\lceil R_1\right\rceil\\ 
H\leq 0    & \text{otherwise} 
\end{cases}$} &    \\
\hline
\end{tabular}
\end{table}

\begin{table} [htbp]
\caption{Sign of the approximation $H(i,k,\ell)$ where $x\in \tilde{\X}_{diff}$ if $\frac{h_1}{\mu_1} \leq \frac{h_2}{\mu_2}$} \label{table:sign2_H}
\centering
\begin{tabular}{|c|c|c|}
\hline

\multicolumn{1}{|l|}{} & \multicolumn{1}{c|}{\textbf{$\ell < C_2$}} & \multicolumn{1}{c|}{\textbf{$\ell\geq C_2$}} \\ \hline

\textbf{$\mu_1 < \mu_2$}    & \text{$
\begin{cases} 
H \leq 0 & \text{if } i<\left\lfloor R_1\right\rfloor+1\\ 
H > 0    & \text{otherwise} 
\end{cases}$}               & \multirow{3}{*}{{$H \leq 0$}}  \\
\cline{1-2}

\textbf{$\mu_1 \geq \mu_2$}    & \text{$H \leq 0$} &    \\
\hline
\end{tabular}
\end{table}

\begin{table} [htbp]
\caption{Sign of the difference $H(i,k,\ell)$ where $x\in \tilde{\X}_{diff}$ if $C_2 \geq C_1$} \label{table:sign_H_enough_C2}
\centering
\begin{tabular}{|c|c|c|}
\hline

\multicolumn{1}{|l|}{} & \multicolumn{1}{c|}{\textbf{$\frac{h_1}{\mu_1} > \frac{h_2}{\mu_2}$}} & \multicolumn{1}{c|}{\textbf{$\frac{h_1}{\mu_1} \leq \frac{h_2}{\mu_2}$}} \\ \hline

\textbf{$\mu_1 < \mu_2$}    &  \multirow{3}{*}{$H > 0$}               &  \text{$ 
\begin{cases} 
H\leq 0 & \text{if } i < \left\lfloor R_1\right\rfloor+1\\ 
H> 0    & \text{otherwise} 
\end{cases}$}  \\
\cline{1-1} \cline{3-3}

\textbf{$\mu_1 = \mu_2$}    &               & \multirow{3}{*}{$H\leq0$}  \\
\cline{1-2}

\textbf{$\mu_1 > \mu_2$}    & \text{$ 
\begin{cases} 
H > 0 & \text{if } i < \left\lceil R_1\right\rceil\\ 
H \leq 0    & \text{otherwise} 
\end{cases}$}               &   \\
\hline
\end{tabular}
\end{table}

\begin{definition}\label{def:thresholds_lowcost}
    Suppose $\frac{h_1}{\mu_1} > \frac{h_2}{\mu_2}$. Consider $k = 1,2, \ldots, C_1$ with $\ell = C_1 - k$.
    The \textbf{actual} and \textbf{estimated integer thresholds} (indexed by $k$) are
    \begin{align}
        i_D(k) := \min\{i|D(i,k,\ell) \leq 0, \text{where }(i,k,\ell)\in \tilde{\X}_{diff}\}, \label{def:iD_lowcost}
    \end{align}
    and
    \begin{align}
        i_H(k) &:= \min\{i|H(i,k,\ell) \leq 0, \text{where }(i,k,\ell)\in \tilde{\X}_{diff}\}, \nonumber
    \end{align}
    respectively.
\end{definition}
Note that in Definition \ref{def:thresholds_lowcost}, a little arithmetic yields
\begin{align}
        i_H(k) &=\left\{
        \begin{array}{lcl}
           0 & \text{if} & \ell \geq C_2 \text{ and } \frac{h_1}{\mu_1} \leq \frac{\ell+1}{C_2}\frac{h_2}{\mu_2},\\
           \max\left\{\left\lceil R_2(k)\right\rceil, 0\right\} & \text{if} & \ell \geq C_2 \text{ and } \frac{h_1}{\mu_1} > \frac{\ell+1}{C_2}\frac{h_2}{\mu_2}, \\
           \left\lceil R_1\right\rceil& \text{if} & \ell <C_2 \text{ and } \mu_1 >\mu_2, \\
           \infty & \text{if} & \ell <C_2 \text{ and } \mu_1 \leq\mu_2, \\
        \end{array}
        \right.  \label{def:iH_lowcost}
\end{align}
where, 
    \begin{equation} \label{eq:thresh_no_queue_in_collab}
        R_1:=-\frac{b}{c} = -\frac{m h_1 - h_2}{\left(m - 1\right)h_0}\cdot C_1, 
    \end{equation}
and 
    \begin{equation} \label{eq:thresh_queue_in_collab}
        R_2(k) := -\frac{b'}{c'} + y_k  = \frac{C_2 m (h_1-h_2)-C_1h_2}{h_0}+ (k-1+ C_2 m ). 
    \end{equation}
    Here, the last equality in \eqref{eq:thresh_no_queue_in_collab} applies the expressions of $b$ and $c$ in \eqref{def:b_and_c} and $m$ in \eqref{def:ratio_of_rates}. Also, note that $R_1 > 0$, since $b>0$ and $c <0$ under the conditions that $\frac{h_1}{\mu_1} > \frac{h_2}{\mu_2}$ and $\mu_1 > \mu_2$.
    Similarly, \eqref{eq:thresh_queue_in_collab} follows by substituting definitions of $b', c'$ in \eqref{def:b'_and_c'} and $y_k$ in \eqref{def:sequence_y} (using $m$ as well) in the last step. 
\begin{definition} \label{def:thresholds_highcost}
    Suppose $\frac{h_1}{\mu_1} \leq \frac{h_2}{\mu_2}$ and consider any $\ell = 0,1, \ldots, C_1-1$. 
    Define the \textbf{actual} and \textbf{estimated integer thresholds} (index by $\ell$, where $k = C_1 - \ell$) as:
    \begin{align}
        \tilde{i}_D(\ell) := \min\{i|D(i,k,\ell) > 0, \text{where }(i,k,\ell)\in \tilde{\X}_{diff}\}, \label{def:iD_highcost}
    \end{align}
    and 
    \begin{align}
        \tilde{i}_H(\ell) &:= \min\{i|H(i,k,\ell) > 0, \text{where }(i,k,\ell)\in \tilde{\X}_{diff}\}, \nonumber \\
    \end{align}
    respectively. 
\end{definition}
Once again, for Definition \eqref{def:thresholds_highcost}, a little algebra yields
\begin{align}
          \tilde{i}_H(\ell) &=
        \begin{cases}
           \left\lfloor R_1\right\rfloor+1 & \text{if }  \ell <C_2 \text{ and } \mu_1 < \mu_2, \\
           \infty & \text{otherwise},
        \end{cases}\label{def:iH_highcost}
\end{align}
where $R_1$ is defined in \eqref{eq:thresh_no_queue_in_collab}.  Notice $R_1 \geq 0$ here, since $\frac{h_1}{\mu_1} \leq \frac{h_2}{\mu_2}$ and $\mu_1 < \mu_2$ imply $b \leq 0$ and $c > 0$ (see \eqref{def:b_and_c}).    
Note that $i_H(k)$ ($\tilde{i}_H(\ell)$) and the sign of $H$ is characterized based on whether $\ell < C_2$ or $\ell \geq C_2$, rather than on the relationship between $\mu_1$ and $\mu_2$, as used to define the thresholds in Tables \ref{table:sign1_D} and \ref{table:sign2_D}. This approach more effectively captures the significant impact of (a second) queueing on the waiting times of subsequent jobs.

With the definitions in place, we now present our main results of this section, summarized in Theorems \ref{thm:heur_lowcost} and \ref{thm:heur_highcost} below, along with Corollary \ref{cor:heur_enough_Type_II}. These results compare the estimated integer threshold $i_H(k)$ (resp. $\tilde{i}_H(\ell)$) with the actual integer threshold $i_D(k)$ (resp. $\tilde{i}_D(\ell)$) in regions of the parameter and state spaces where these thresholds are finite --- outside of which both are infinite --- thereby validating the accuracy of the approximation. In addition, several analytical properties of $i_D(k)$ and $\tilde{i}_D(\ell)$ are established. We begin by introducing two sufficient conditions and supporting notation that underpin parts of Theorem \ref{thm:heur_lowcost}.

\begin{definition} \label{def:probs_k}
    For $k = 1,2, \ldots, C_1$, define the following probabilities $p_k,q_k,r_k$ indexed by $k$:
    \begin{align}
        p_k 
        & := \left\{
        \begin{array}{lcl}
            \frac{C_1 \mu_1\mu_2}{d(k,\ell)d(k-1,\ell+1)} & \text{if} & \ell < C_2, \\
            \frac{C_2 \mu_1\mu_2}{d(k,\ell)d(k-1,\ell+1)} & \text{if} & \ell \geq C_2, \\
        \end{array} 
        \right. \nonumber \\
        & \ = \left\{
        \begin{array}{lcl}
             \frac{C_1 m}{\big(k+\ell m\big)\big((k-1) + (\ell+1)m\big)} & \text{if} & \ell < C_2, \\
             \frac{C_2 m}{\big(k+C_2 m\big)\big((k-1)+C_2 m\big)} & \text{if} & \ell\geq C_2,
        \end{array}
        \right. \label{def:pk}\\
        q_k &:= \frac{(k-1) \mu_1}{d(k-1,\ell+1)} 
        = \frac{k-1}{(k-1)+\min\{\ell+1,C_2\}m},
        \label{def:qk}\\
        r_k &:= \frac{\min\{\ell,C_2\}\mu_2}{d(k,\ell)} 
        = \frac{\min\{\ell,C_2\} m}{k + \min\{\ell,C_2\}m},
        \label{def:rk}
    \end{align}
    where $\ell = C_1-k$, and $m$ is the ratio of service rates defined in \eqref{def:ratio_of_rates}. Note we suppress the dependence on $m$ to simplify the notation, albeit they are functions of $m$.
\end{definition}

\begin{condition} \label{cond:sufficient_equal_thresh}
    Consider $k$ such that $\ell = C_1-k \geq C_2$. For $r : = \frac{C_2 m}{1+ C_2 m}$, define the following two conditions:
        \begin{align}
            y_{k}r^{i_H(k)-k} &< \frac{h_0}{h_0-h_2}\Big(R_2(k)-\big(i_H(k)-1\big)\Big), \label{cond:highcost_queue}\\
            \frac{r_k}{1-q_k}y_{k}r^{i_H(k)-k} &\leq \frac{h_0}{h_2-h_0}\Big(i_H(k)-R_2(k)\Big), \label{cond:highcost_collab}
        \end{align}
    where $R_2(k)$ is given in \eqref{eq:thresh_queue_in_collab},
    and $q_k$ and $r_k$ are specified in \eqref{def:qk} and \eqref{def:rk}, respectively.
\end{condition}

\begin{theorem} \label{thm:heur_lowcost}
    Suppose $\frac{h_1}{\mu_1} > \frac{h_2}{\mu_2}$ and consider any $k = 1,2, \ldots,C_1$. 
    \begin{enumerate}
        \item \label{state:thresh_mono}
        $i_D(k)$ is non-decreasing in $k$, i.e., $i_D(k-1) \leq i_D(k)$, for $k = 2, 3, \ldots, C_1$.
        \item  \label{state:heur_bd}
        Fix $k$ (and $\ell$) such that $\ell = C_1-k<C_2$. If $\mu_1>\mu_2$, then $i_H(k) \leq i_D(k) \leq i_H(k) + (C_1-1)$. In the case where $\mu_1 \leq \mu_2$, both $i_D(k)$ and $i_H(k)$ are infinite.
        \item Fix $k$ such that $\ell = C_1-k \geq C_2$.
        \begin{enumerate}[ref=\theenumi(\alph*)]
            \item \label{state:always_neg_queue_in_collab}
            We have $\frac{h_1}{\mu_1} \leq \frac{\ell+1}{C_2}\frac{h_2}{\mu_2}$ if and only if $i_D(k) = 0$. 
            In either case (if $\frac{h_1}{\mu_1} \leq \frac{\ell+1}{C_2}\frac{h_2}{\mu_2}$ or if $i_D(k) = 0$), $i_H(k) = 0$;
            \item \label{state:has_pos_queue_in_collab}
            Suppose $\frac{h_1}{\mu_1} > \frac{\ell+1}{C_2} \frac{h_2}{\mu_2}$, i.e., $i_D(k) \geq 1$.
            \begin{enumerate}[ref=\theenumii\roman*]
                \item \label{state:thresh_queue_in_collab}
                If $h_0 \geq h_2$ (resp. $h_0 \leq h_2$), we have $i_H(k) \geq i_D(k)$ (resp. $i_H(k) \leq i_D(k)$).
                Consequently, if $h_0=h_2$, then $i_H(k) = i_D(k)$.
                \item \label{state:equal_thresh_queue_in_collab}
                If $h_0 > h_2$ (resp. $h_0 < h_2$), and $k$ satisfies \eqref{cond:highcost_queue} (resp.\eqref{cond:highcost_collab}) in Condition \ref{cond:sufficient_equal_thresh}, then $i_H(k) = i_D(k)$. 
                \item \label{state:sequence_iD}
                If either $h_0 = h_2$ or $h_0 > h_2$ (resp. $h_0 < h_2$) with $k$ satisfying \eqref{cond:highcost_queue} (resp. \eqref{cond:highcost_collab}) in Condition \ref{cond:sufficient_equal_thresh}, then $i_D(k)=i_D(k-1)+1$, for all $k \geq 2$ such that $\ell = C_1-k \geq C_2$.
            \end{enumerate}            
        \end{enumerate}
    \end{enumerate}
\end{theorem}

\begin{theorem}\label{thm:heur_highcost}
    Suppose $\frac{h_1}{\mu_1} \leq \frac{h_2}{\mu_2}$. Consider any $\ell = 0,1, \ldots,C_1-1$. 
    \begin{enumerate}
        \item \label{state:thresh_mono_highcost}
        $\tilde{i}_D(\ell)$ is non-decreasing in $\ell$. That is, $\tilde{i}_D(\ell-1) \leq \tilde{i}_D(\ell)$, where $\ell = 1,2, \ldots, C_1-1$.
        \item \label{state:heur_bd2}
        If $\mu_1<\mu_2$ and $\ell <C_2$, we have $\tilde{i}_H(\ell) \leq \tilde{i}_D(\ell) \leq \tilde{i}_H(\ell) + (C_1-1)$; otherwise (either $\mu_1 \geq \mu_2$ and $\ell < C_2$ or $\ell \geq C_2$), both $\tilde{i}_H(\ell)$ and $\tilde{i}_D(\ell)$ are infinite.
    \end{enumerate}
    
\end{theorem}

\begin{corollary} \label{cor:heur_enough_Type_II}
    Assume $C_2 \geq C_1$. Consider any $x \in \X_{diff}$. The following holds
    \begin{enumerate}
        \item If $\frac{h_1}{\mu_2} > \frac{h_2}{\mu_2}$,
        \begin{enumerate}
            \item \label{state:thresh_mono_enough_Type_II}
            $i_D(k-1) \leq i_D(k)$, for $k = 2,3, \ldots, C_1$.
            \item \label{state:heur_bd_enough_Type_II}
            if $\mu_1 > \mu_2$, we have $i_H(k) \leq i_D(k) \leq i_H(k) + (C_1-1)$, implying $i_H(k) = i_D(k)$ if $C_1=1$; and both $i_D(k)$ and $i_H(k)$ are infinite if $\mu_1 \leq \mu_2$. 
        \end{enumerate}
        \item If $\frac{h_1}{\mu_2} \leq \frac{h_2}{\mu_2}$,
        \begin{enumerate}
            \item \label{state:thresh_mono_ell_enough_Type_II}
            $\tilde{i}_D(\ell-1) \leq \tilde{i}_D(\ell)$, for $\ell = 1, 2, \ldots, C_1-1$.
            \item \label{state:heur_bd2_enough_Type_II}
            if $\mu_1 < \mu_2$, we have $\tilde{i}_H(\ell) \leq \tilde{i}_D(\ell) \leq \tilde{i}_H(\ell) + (C_1-1)$, implying $\tilde{i}_H(\ell) = \tilde{i}_D(\ell)$ if $C_1=1$; and both $\tilde{i}_H(\ell)$ and $\tilde{i}_D(\ell)$ are infinite if $\mu_1 \geq \mu_2$. 
        \end{enumerate}
    \end{enumerate}
\end{corollary}

\subsection{Preliminary results used to prove bounds}\label{sec:prelim}
This subsection presents several lemmas comparing $D(i,k,\ell)$ with $H(i,k,\ell)$. Since they are both monotone in the same direction in regions where an estimation of the threshold is required, this helps with the comparison of $i_H(k)$ and $i_D(k)$ (or $\tilde{i}_H(\ell)$ and $\tilde{i}_D(\ell)$). 

Along with $p_k,q_k$ and $r_k$ in Definition \ref{def:probs_k}, we define quantities $b_k$ and $c_k$ below to further simplify the notation in the Bellman equations.

\begin{definition}\label{def:probs_quantities_k}
    For $k = 1, 2, \ldots,C_1-1$, define
    \begin{align}
        b_k&:= 
        \left\{
        \begin{array}{lcl}
            b & \text{if} & \ell = C_1 - k < C_2, \\
           b' & \text{if} & \ell = C_1 - k \geq C_2, 
        \end{array}
        \right. \label{def:bk}\\
         c_k &:= 
        \left\{
        \begin{array}{lcl}
            c & \text{if} & \ell = C_1 - k < C_2, \\
            c' & \text{if} & \ell = C_1 - k \geq C_2, \\
        \end{array}
        \right. \label{def:ck}
    \end{align}
    where $b,c$ are given in \eqref{def:b_and_c} and $b',c'$ are defined in \eqref{def:b'_and_c'}.
\end{definition}

The following lemma further investigates the properties of $D(i,k,\ell)$. 
Statement \ref{state:diff_bdry} specifies its values at $i=0$. 
Statement \ref{state:diff_ub_larger_mu2} establishes an affine upper bound when $\mu_2 \geq \mu_1$, which coincides with $H(i,k,\ell)$ if $\frac{h_1}{\mu_1} \leq \frac{h_2}{\mu_2}$ also holds. This result will be revisited when comparing $\tilde{i}_H(\ell)$ and $\tilde{i}_D(\ell)$.
Statement \ref{state:recursive_pos} reformulates the expression for $D(i,k,\ell)$ when $D(i-1,k,\ell) \geq 0$, using notations in Definitions \ref{def:probs_k} and \ref{def:probs_quantities_k}.
The proofs of Statements \ref{state:diff_bdry} and \ref{state:diff_ub_larger_mu2} are provided in Appendix \ref{sec:proof-prelims}, while the proof of Statement \ref{state:recursive_pos} is in Appendix \ref{sec:prelim-heuristics}.

\begin{lemma} \label{lemma:diff}
    For any $x =(i,k,\ell)\in \X_{diff}$, the following holds.
    \begin{enumerate}
        \item \label{state:diff_bdry}
        If $i=0$, 
        \begin{align}
            D(0,k,\ell) & = \left\{
          \begin{array}{lcl}
            \frac{h_1}{\mu_1} - \frac{h_2}{\mu_2} & \text{if} & \ell < C_2, \\
           \frac{h_1}{\mu_1} - \frac{(\ell+1) h_2}{C_2\mu_2} & \text{if} & \ell \geq C_2, \\
            \end{array}
            \right. \nonumber \\
            & = \frac{h_1}{\mu_1} - \frac{\max\{\ell+1,C_2\}}{C_2}\frac{h_2}{\mu_2}. \label{eq:diff_bdry}
        \end{align}
        \item \label{state:diff_ub_larger_mu2}
        Suppose $\mu_2\geq\mu_1$ (no need to assume $\frac{h_1}{\mu_1} \geq \frac{h_2}{\mu_2}$). The difference $D(i,k,\ell)$ is upper bounded by an affine function in $i$:
        \begin{align}
                D(i,k,\ell) 
                \leq ic+b, \label{eq:diff_ub_larger_mu2}
            \end{align}
        where $c$ and $b$ are defined in \eqref{def:b_and_c}.
        \item \label{state:recursive_pos}
        If $i\geq 1$ and $(i,k,\ell) \in \tilde{\X}_{diff}$, $D(i-1,k,\ell) \geq 0$ implies
        \begin{align}
            D(i,k,\ell) = p_k(ic_k+b_k) + q_k \max\{0,D(i-1,k-1,\ell+1)\} + r_k D(i-1,k,\ell), \label{eq:recursive_pos}
        \end{align}
        where $p_k, q_k$ and $r_k$ are given in Definition \ref{def:probs_k} and $b_k$ and $c_k$ are provided in Definition \ref{def:probs_quantities_k}.
    \end{enumerate}
\end{lemma}

In addition to $y_k$ defined in \eqref{def:sequence_y}, we introduce another (non-negative) sequence $z_\ell$, which facilitates the comparison between $\tilde{i}_H(\ell)$ and $\tilde{i}_D(\ell)$, particularly under the condition$\frac{h_1}{\mu_1} \leq \frac{h_2}{\mu_2}$.

\begin{definition} 
    For $\ell = 0, 1, \ldots, C_2-1$, define
        \begin{align}
            z_\ell := \frac{k-1}{m}+\ell, \label{def:sequence_z}
        \end{align} 
        where $k=C_1-\ell$ and $m>0$ is given in \eqref{def:ratio_of_rates}.
\end{definition}

The following lemma explores additional bounds on $D(i,k,\ell)$ under various parameter conditions, which directly relate to $H(i,k,\ell)$ within the corresponding parameter space. Additionally, it further investigates properties of $i_D(k)$ and $i_H(k)$. The proof of the entire lemma is given in Appendix \ref{sec:prelim-heuristics}.

\begin{lemma} \label{lemma:heur}
    Consider $x =(i,k,\ell)\in \tilde{\X}_{diff}$ with $k+\ell=C_1$. The following results concern the sequences $y_k$ and $z_\ell$ defined in \eqref{def:sequence_y} and \eqref{def:sequence_z}, respectively. Recall that the constants $c,b$ (from \eqref{def:b_and_c}) and $c',b'$ (from \eqref{def:b'_and_c'}), are independent of the system state.
    \begin{enumerate}[ref=\arabic*]
        \item Fix any $k = 1, 2, \ldots, C_1$ (indexed by $k$).
        \begin{enumerate}[ref=\theenumi(\alph*)]
            \item \label{state:y_ub}
            If $\mu_1\geq \mu_2$, then $y_k \leq C_1-1$.
            \item \label{state:heur_linear_ub}
            If $\mu_1 >\mu_2$ and $D(i,k,\ell) \geq 0$, then $D(i,k,\ell) \leq (i-y_k)c+b$.
            \item \label{state:queue_in_collab_bd}
            Suppose $\frac{h_1}{\mu_1} \geq \frac{h_2}{\mu_2}$ and $\ell \geq C_2$.
            \begin{enumerate}[ref=\theenumii\roman*]
                \item \label{state:queue_in_collab_costly_queue_ub}
                If $h_0 \geq h_2$ and $D(i,k,\ell) \geq 0$, then $D(i,k,\ell) \leq (i-y_k)c'+b'$.
                \item \label{state:queue_in_collab_costly_collab_lb}
                If $h_0 \leq h_2$, then $D(i,k,\ell) \geq (i-y_k)c'+b'$.
            \end{enumerate}
            Consequently, if $h_0 = h_2$ and $D(i,k,\ell) \geq 0$, then $D(i,k,\ell) = (i-y_k)c'+b'$.
            \item \label{state:cond_k-1}
            Suppose $\ell \geq C_2$, $\frac{h_1}{\mu_1} > \frac{\ell+1}{C_2}\frac{h_2}{\mu_2}$, and $i_H(k) \geq 1$. The following hold: 
            \begin{enumerate}[ref=\theenumii\roman*]
                \item \label{state:sequence_iH}
                $i_H(k) = i_H(k-1)+1$.
                \item \label{tate:cond_k_k-1}
                If $k$ satisfies \eqref{cond:highcost_queue} (or \eqref{cond:highcost_collab}) in Condition \ref{cond:sufficient_equal_thresh}, then $k-1$ does as well. 
            \end{enumerate}
            \item \label{state:asymp}
            Suppose $\ell \geq C_2$, $\frac{h_1}{\mu_1} > \frac{\ell+1}{C_2}\frac{h_2}{\mu_2}$, $h_0 > h_2$ (resp. $h_0 < h_2$), and $k$ satisfies \eqref{cond:highcost_queue} (resp. \eqref{cond:highcost_collab}) in Condition \ref{cond:sufficient_equal_thresh}. 
            The following hold:
            \begin{enumerate}[ref=\theenumii\roman*]
                \item \label{state:exp_bound}
                If 
                $0 \leq i < i_D(k)$, then $D(i,k,\ell) = (i-y_k)c'+b' + \varphi_k(i)$, where $\varphi_k(i) \geq \frac{h_2-h_0}{C_2 \mu_2}y_k r^{i-k+1}$ (resp. $\varphi_k(i) \leq \frac{h_2-h_0}{C_2 \mu_2}y_k r^{i-k+1}$).
                \item \label{state:thresh_reverse}
                $i_H(k) \leq i_D(k)$ (resp. $i_H(k) \geq i_D(k)$).
            \end{enumerate}
        \end{enumerate}
        \item Fix any $\ell = 0, 1, \ldots, C_2-1$, i.e., $\ell < C_2$ (indexed by $\ell$).
        \begin{enumerate}[ref=\theenumi(\alph*)]
            \item \label{state:z_ub}
            If $\mu_2 \geq \mu_1$, then $z_\ell \leq C_1-1$.
            \item \label{state:heur_linear_lb}
            If $\mu_1 <\mu_2$, $\frac{h_1}{\mu_1} \leq \frac{h_2}{\mu_2}$, and $D(i,k,\ell) \leq 0$, then $D(i,k,\ell) \geq (i-z_\ell)c+b$.
        \end{enumerate}
    \end{enumerate}
\end{lemma}

\subsection{Independent service costs higher than collaborative \texorpdfstring{$\big(\frac{h_1}{\mu_1} > \frac{h_2}{\mu_2}\big)$}{(h1/mu1 >= h2/mu2)}}
Suppose throughout the remainder of this subsection that $\frac{h_1}{\mu_1} > \frac{h_2}{\mu_2}$. 
We begin by proving Statement \ref{state:thresh_mono} in Theorem \ref{thm:heur_lowcost}, which characterizes how the actual threshold changes as the number of jobs in the independent (or equivalently, collaborative) service station varies.

\begin{proof}[Proof of Statement \ref{state:thresh_mono} in Theorem \ref{thm:heur_lowcost}.]
    This is a direct application of Proposition \ref{prop:mono_ell} which explains that holding the total number of customers in service fixed, $D$ is a non-decreasing function of $k$.
\end{proof}
\subsubsection{No queueing for collaborative service in \texorpdfstring{$(i,k-1,\ell+1)$}{(i,k-1,ell+1)} \texorpdfstring{($\ell < C_2$)}{(ell < C2)}}
Noting from Table \ref{table:sign1_D}, for $\ell < C_2$, it remains to consider the case $\mu_1 > \mu_2$, where $D(i,k,\ell)$ changes sign as $i$ increases. 
Statement \ref{state:heur_bd} in Theorem \ref{thm:heur_lowcost} establishes both upper and lower bounds for $i_D(k)$ relative to $i_H(k)$ in this case, while both thresholds are infinity when $\mu_1 \leq \mu_2$. The proof is provided below.

\begin{proof}[Proof of Statement \ref{state:heur_bd} in Theorem \ref{thm:heur_lowcost}.]
    Consider first $\mu_1 > \mu_2$.
    
    We prove the direction $i_H(k) \leq i_D(k)$, where $k$ is such that $\ell = C_1-k < C_2$, by showing $H\big(i_D(k),k,\ell\big) \leq 0$. Notice the boundary values in \eqref{eq:diff_bdry} imply that for $\ell < C_2$,
    \begin{align*}
        D(0,k,\ell) = \frac{h_1}{\mu_1} - \frac{h_2}{\mu_2} = b = H(0,k,\ell) > 0,
    \end{align*}
    where the second equality relies on the definition of $b$ in \eqref{def:b_and_c}. This ensures $i_H(k) \geq 1$ and $i_D(k) \geq 1$ (see \eqref{def:iH_lowcost} and \eqref{def:iD_lowcost}). By the definition of $i_D(k)$, we have $D\big(i_D(k),k,\ell\big)\leq 0$, and $D\big(i_D(k)-1,k,\ell\big) > 0$. Therefore, the equation \eqref{eq:recursive_pos} holds in particular at $i_D(k)$:
    \begin{align}
        0 & \geq D\big(i_D(k),k,\ell\big) \nonumber \\
         & = p_k \big(i_D(k) \cdot c_k+b_k\big) + q_k \max\{0,D\big(i_D(k)-1,k-1,\ell+1\big) + r_k D\big(i_D(k)-1,k,\ell\big) \nonumber \\
         &\geq p_k \big(i_D(k) \cdot c_k+b_k\big) \label{eq:neg_iD}
    \end{align}where the last inequality holds by recalling from Definition \ref{def:probs_k}, $p_k, q_k$ and $r_k$ are all non-negative, as defined in \eqref{def:pk}, \eqref{def:qk} and \eqref{def:rk}, respectively.
    Furthermore, applying $b_k = b$ and $c_k = c$ when $\ell < C_2$ (from their definitions in \eqref{def:bk} and \eqref{def:ck}, respectively) in \eqref{eq:neg_iD} yields
    \begin{align*}
        0 \geq p_k \big(i_D(k) \cdot c+b\big) = p_k H\big(i_D(k),k,\ell\big).
    \end{align*}
    The result that $H\big(i_D(k),k,\ell\big) \leq 0$ thus follows.

    We prove the remaining inequality $i_D(k) \leq i_H(k) + (C_1-1)$ by showing that $D(i,k,\ell) > 0$ implies $H\big(i-(C_1-1)\big) = \big(i-(C_1-1)\big)c+b > 0$, or equivalently, $i < \frac{b}{-c} + (C_1-1)$ since $\mu_1 > \mu_2$ implies $c<0$ (see \eqref{def:b_and_c}). To that end, referring to Statement \ref{state:heur_linear_ub} in Lemma \ref{lemma:heur}, we observe that $D(i,k,\ell) > 0$ implies $0 < D(i,k,\ell) \leq (i-y_k)c+b$. Consequently,
    \begin{align*}
        i < \frac{b}{-c} + y_k < \frac{b}{-c} + (C_1-1),
    \end{align*}
    where the rightmost inequality applies Statement \ref{state:y_ub} in Lemma \ref{lemma:heur}.

    Suppose now $\mu_1 \leq \mu_2$ for $\ell <C_2$. Then by Proposition \ref{prop:larger_mu2_always_action1}, we have $D(i,j,k,\ell) > 0$, implying $i_D(k) = \infty$ by its definition in \eqref{def:iD_lowcost}. Moreover, $i_H(k) = \infty$ as given in \eqref{def:iH_lowcost}. The result thus follows.
\end{proof}

\begin{figure*} 
    \centering
    \begin{subfigure}[htbp]{0.49\textwidth}
        \centering
        \includegraphics [scale=.6] {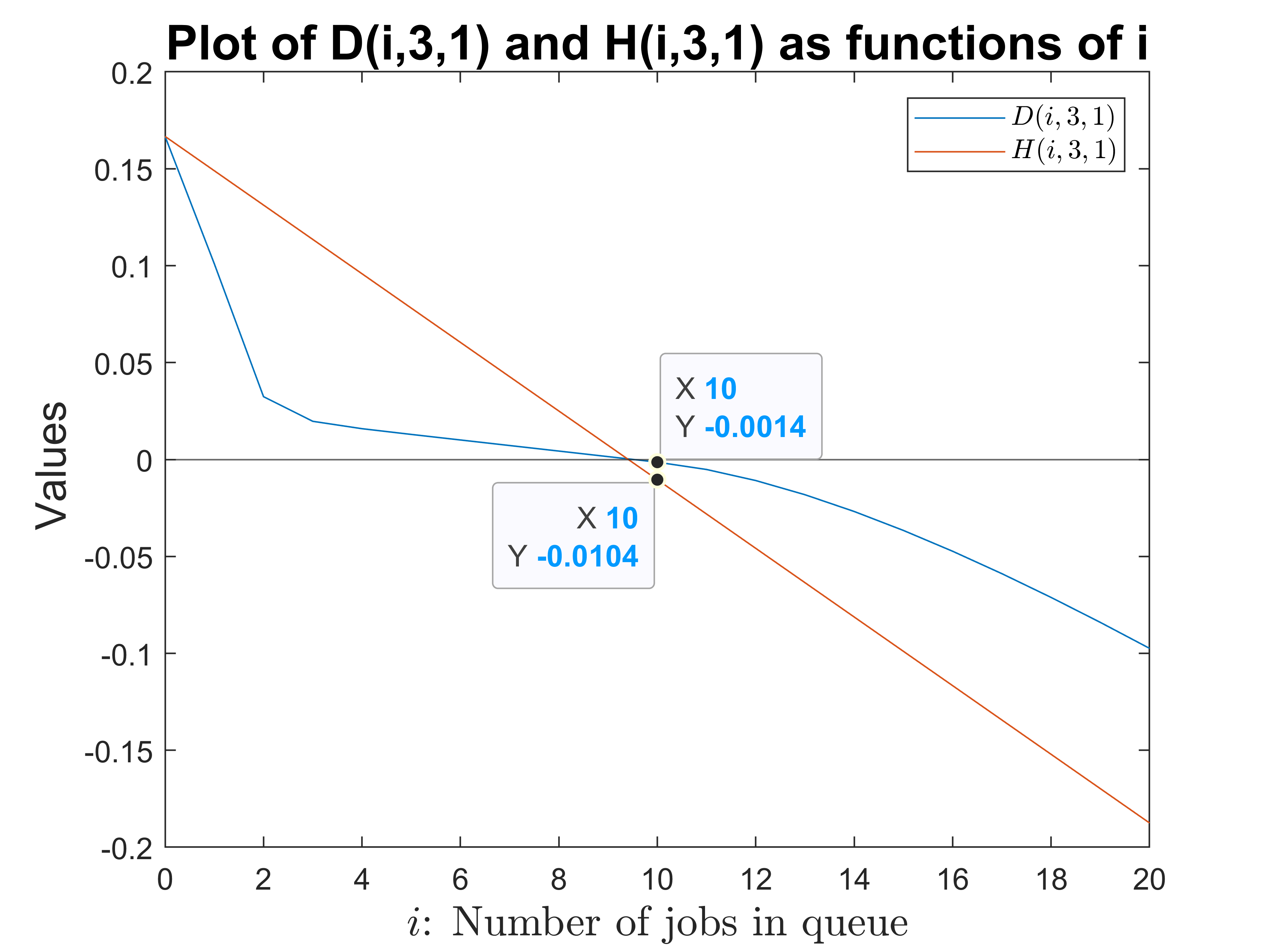}
        \caption{Example \ref{eg:best1}.}
    \end{subfigure}
    \hfill
    \begin{subfigure}[htbp]{0.49\textwidth}
        \centering
        \includegraphics [scale=.6] {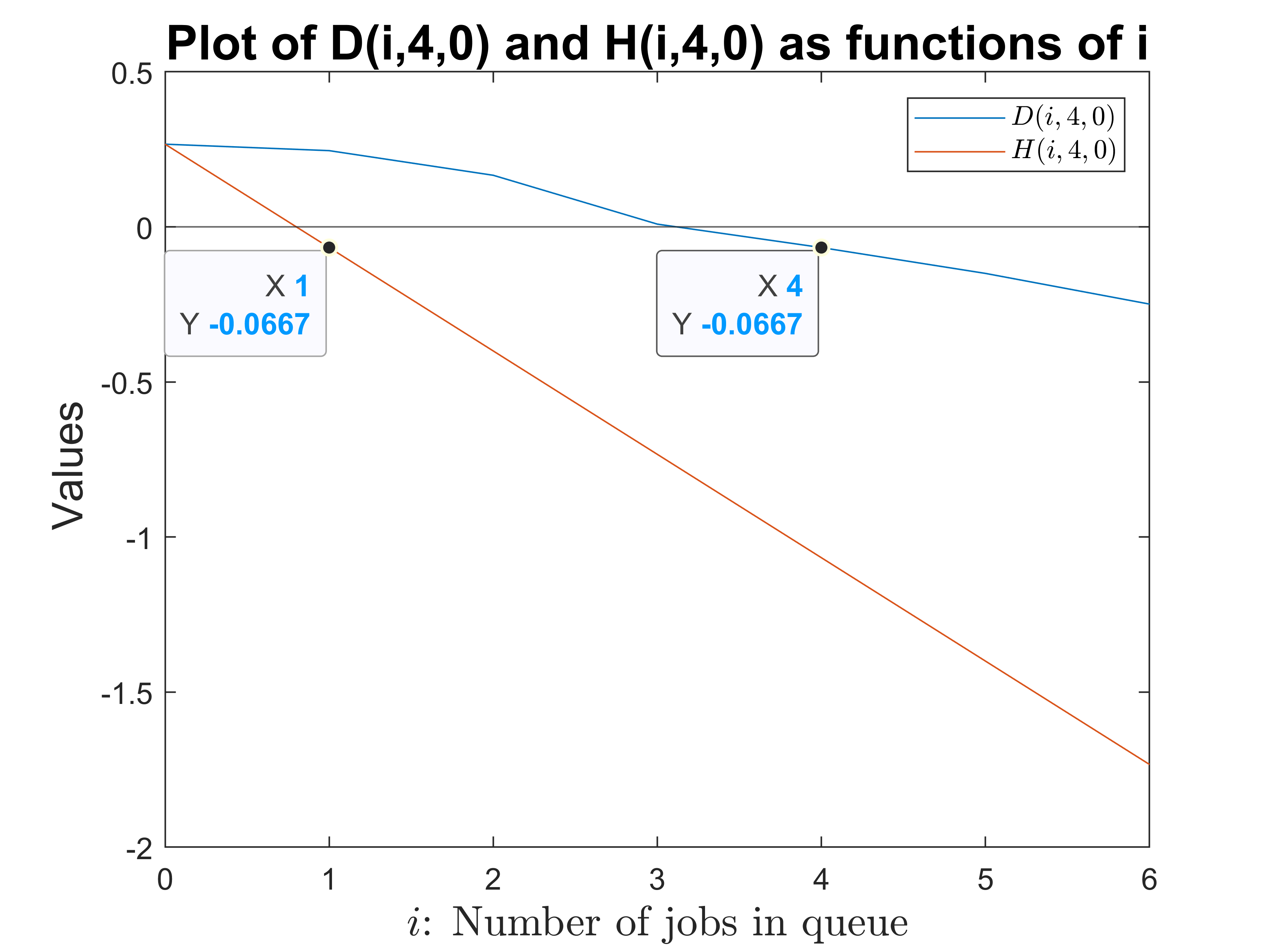}
        \caption{Example \ref{eg:worst1}.}
    \end{subfigure}
    \caption{Plot of $D(i,k,\ell)$ and $H(i,k,\ell)$ as functions of $i$ under parameter configurations in Examples \ref{eg:best1} and \ref{eg:worst1}.}\label{fig:eg_lowcost_no_wait}
\end{figure*}

The bounds established in Statement \ref{state:heur_bd} in Theorem \ref{thm:heur_lowcost} are indeed attainable, as demonstrated in the following two examples. For a visual representation, see Figure \ref{fig:eg_lowcost_no_wait}.

\begin{example}\label{eg:best1}
    Let $C_1 = 4$, $C_2 = 2$, $\mu_1 = 3$, $\mu_2 = 0.96$, $h_0 = 0.1$, $h_1 = 1$, $h_2 = 0.16$, satisfying $\frac{h_1}{\mu_1}>  \frac{h_2}{\mu_2}$. Under these parameters, we obtain $i_D(3) = i_H(3)=10$.
\end{example}
\begin{example}\label{eg:worst1}
    Let $C_1 = 4$, $C_2 = 2$, $\mu_1 = 3$, $\mu_2 = 0.6$, $h_0 = 1$, $h_1 = 1$, $h_2 = 0.04$, satisfying $\frac{h_1}{\mu_1}>  \frac{h_2}{\mu_2}$. Under these parameters, we obtain $i_D(4)=4$ and $i_H(4) = 1$, so that $i_D(4) = i_H(4)+(C_1-1)$.
\end{example}
\begin{remark}\label{remark:worst_lowcost}
    Experimenting with different parameters demonstrates sufficient accuracy of the estimation $i_H$ for the optimal threshold $i_D(k)$. Notably, poorer performance (e.g., the worst-case scenario described in Example \ref{eg:worst1}) tends to occur primarily when $i_D(k)$ (or $i_H(k)$) is small, such as $i_D(k) \leq C_1$ (or $i_H(k) \leq 1)$. 
\end{remark}
\subsubsection{Queueing effects for collaborative service in \texorpdfstring{$(i,k-1,\ell+1)$}{(i,k-1,ell+1)} \texorpdfstring{($\ell \geq C_2$)}{(ell >= C2)}} \label{sec:heur_queue_in_collab}
There are two cases that distinguish the discussion for $\ell \geq C_2$, based on the inequality $\frac{h_1}{\mu_1} > \text{ ($\leq$) }\frac{\ell+1}{C_2}\frac{h_2}{\mu_2}$. In both cases, however, $D(i,k,\ell)$ is monotone non-increasing in $i$ (Proposition \ref{prop:larger_mu1_dec} if $\mu_1\geq \mu_2$ and Proposition \ref{prop:larger_mu2_dec} if $\mu_2 > \mu_1$).
\begin{enumerate}[label= \textbf{Case} \arabic*:,ref={ \arabic*}, leftmargin=3\parindent]
    \item \label{case:always_neg}
    If $\frac{h_1}{\mu_1} \leq \frac{\ell+1}{C_2}\frac{h_2}{\mu_2}$, then $D(i,k,\ell) \leq 0$ for all $i \geq 0$ (Corollary \ref{cor:always_action0}), implying $i_D(k) = 0$ (recall its definition in \eqref{def:iD_lowcost}). From \eqref{def:iH_lowcost}, we also have $i_H(k) = 0$.
    \item \label{case:pos_to_neg}
    If $\frac{h_1}{\mu_1} > \frac{\ell+1}{C_2}\frac{h_2}{\mu_2}$, then $D(0,k,\ell) >0$ (see \eqref{eq:diff_bdry}), corresponding to $i_D(k) \geq 1$. 
\end{enumerate} 
The proof of Statement \ref{state:always_neg_queue_in_collab} in Theorem \ref{thm:heur_lowcost} is now in reach and is provided below.
\begin{proof} [Proof of Statement \ref{state:always_neg_queue_in_collab} in Theorem \ref{thm:heur_lowcost}.]
    Suppose first $i_D(k)=0$. From \eqref{def:iD_lowcost}, this implies $D(0,k,\ell) \leq 0$, yielding that $\frac{h_1}{\mu_1} \leq \frac{\ell+1}{C_2}\frac{h_2}{\mu_2}$ (the contrapositive of the implication in Case \ref{case:pos_to_neg}). 
    Conversely, if $\frac{h_1}{\mu_1} \leq \frac{\ell+1}{C_2}\frac{h_2}{\mu_2}$, then by the discussion in Case \ref{case:always_neg}, we have $i_H(k) = i_D(k) = 0$.
\end{proof}
It remains to compare $i_H(k)$ with $i_D(k)$ when $\frac{h_1}{\mu_1} > \frac{\ell+1}{C_2}\frac{h_2}{\mu_2}$. In this case, given that $i_D(k) \geq 1$, one could alternatively redefine $H$ so that $i_H(k) = \max\{\lceil R_2(k)\rceil, 1\}$, instead of $i_H(k) = \max\{\lceil R_2(k)\rceil, 0\}$ as given in \eqref{def:iH_lowcost}, where $R_2$ is specified in \eqref{eq:thresh_queue_in_collab}. 
However, this minor refinement does not impact the essential ideas or results. We therefore retain the original definition for simplicity and consistency.

When $\frac{h_1}{\mu_1} > \frac{\ell+1}{C_2}\frac{h_2}{\mu_2}$, Statement \ref{state:thresh_queue_in_collab} in Theorem \ref{thm:heur_lowcost} establishes that $i_H(k)$ serves as an upper bound for $i_D(k)$ when staying at Station 0 is more costly than at Station 2 ($h_0 \geq h_2$), and a lower bound when the opposite ($h_0 \leq h_2$) holds. The corresponding proof is given below.

\begin{figure*} 
    \centering
    \begin{subfigure}[htbp]{0.49\textwidth}
        \centering
        \includegraphics [scale=.6] {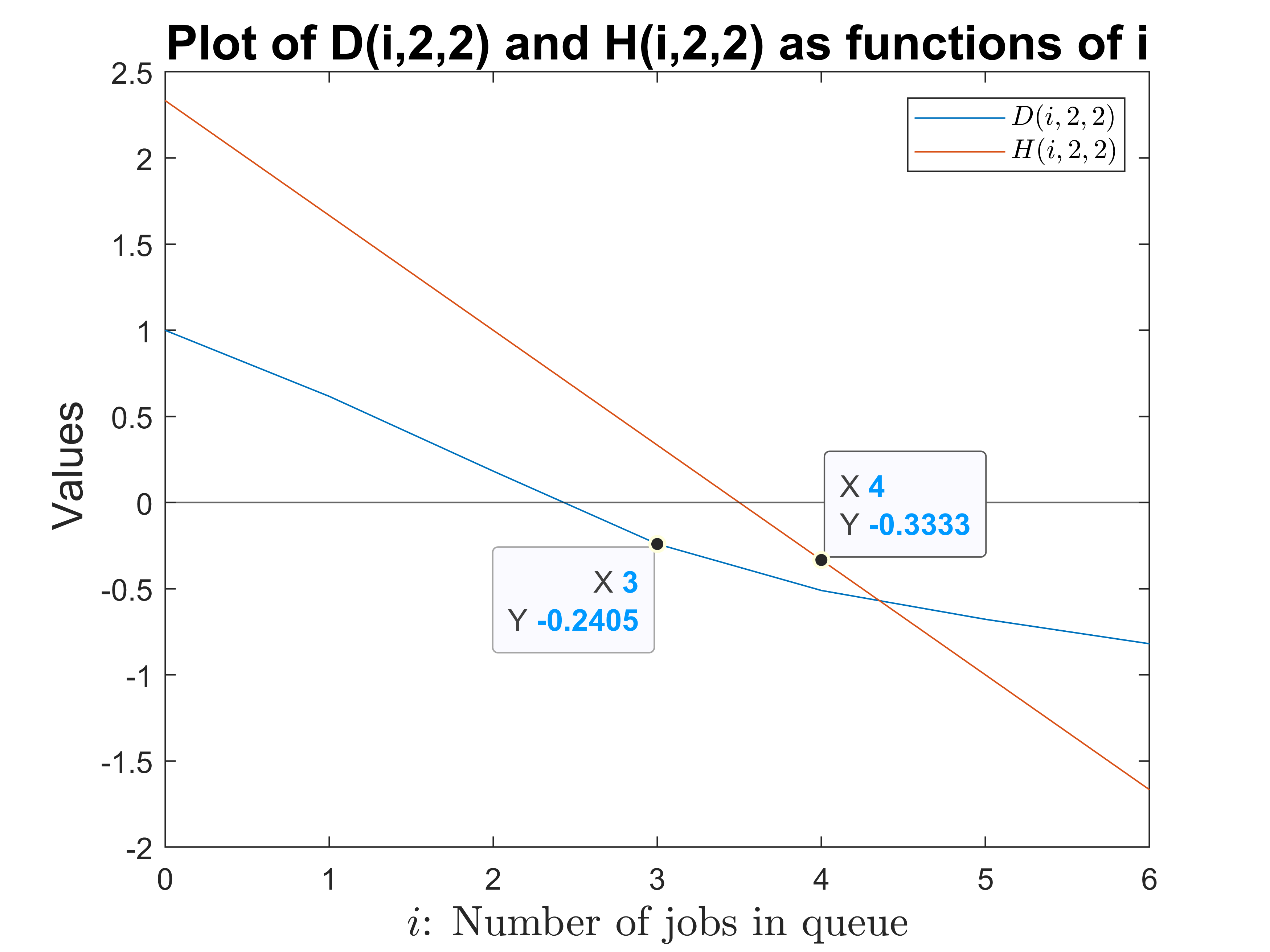}
        \caption{$i_D(2)=3$ and $i_H(2) = 4$ when $h_1 = 2$.}
    \end{subfigure}
    \hfill
    \begin{subfigure}[htbp]{0.49\textwidth}
        \centering
        \includegraphics [scale=.6] {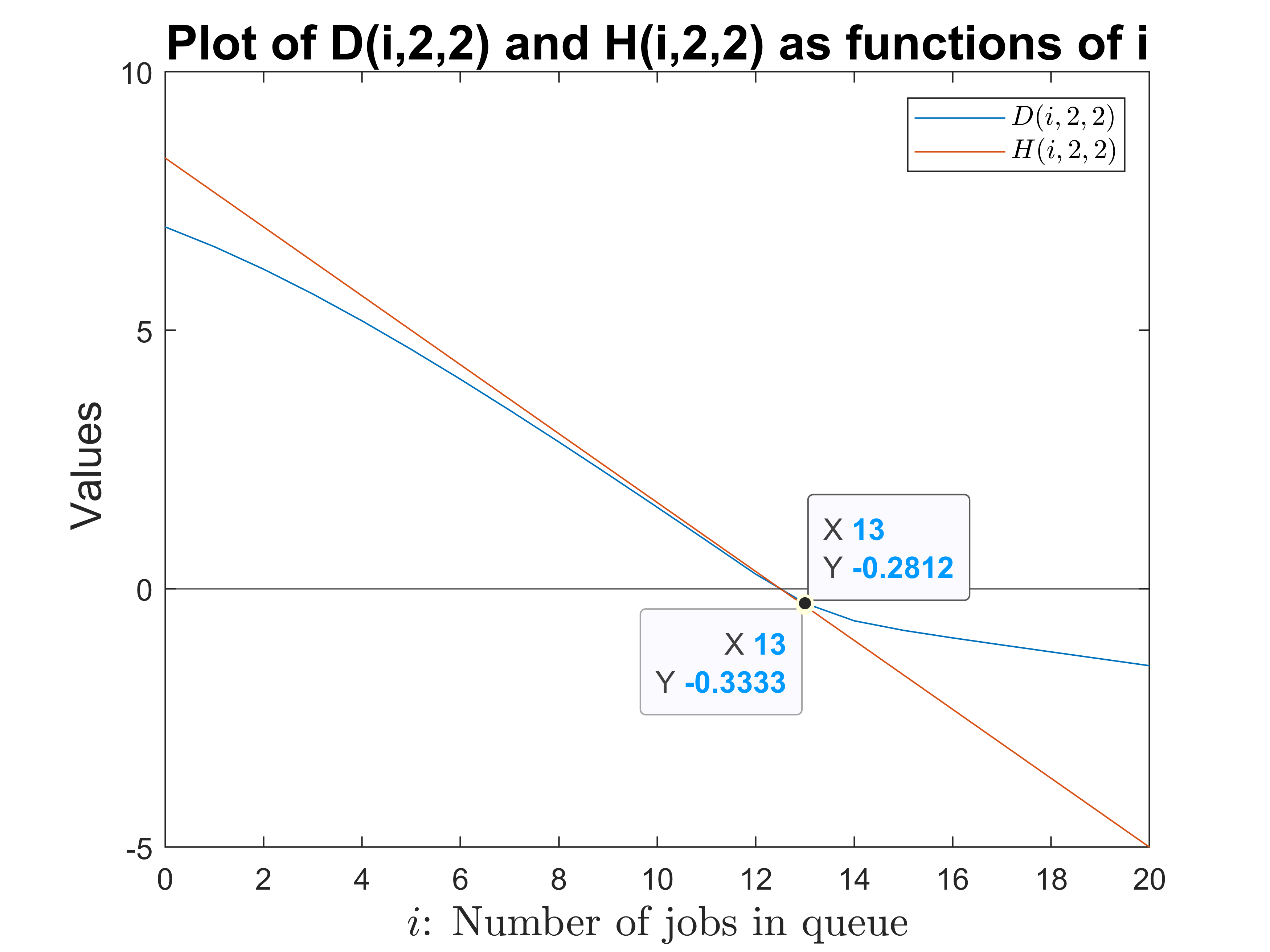}
        \caption{$i_D(2) = i_H(2)= 13$ when $h_1 = 8$.}
    \end{subfigure}
    \caption{Plot of $D(i,k,\ell)$ and $H(i,k,\ell)$ as functions of $i$ under parameter configurations in Example \ref{eg:costly_queue} such that $h_0>h_2$.}\label{fig:costly_queue}
\end{figure*}

\begin{figure*} 
    \centering
    \begin{subfigure}[htbp]{0.49\textwidth}
        \centering
        \includegraphics [scale=.6] {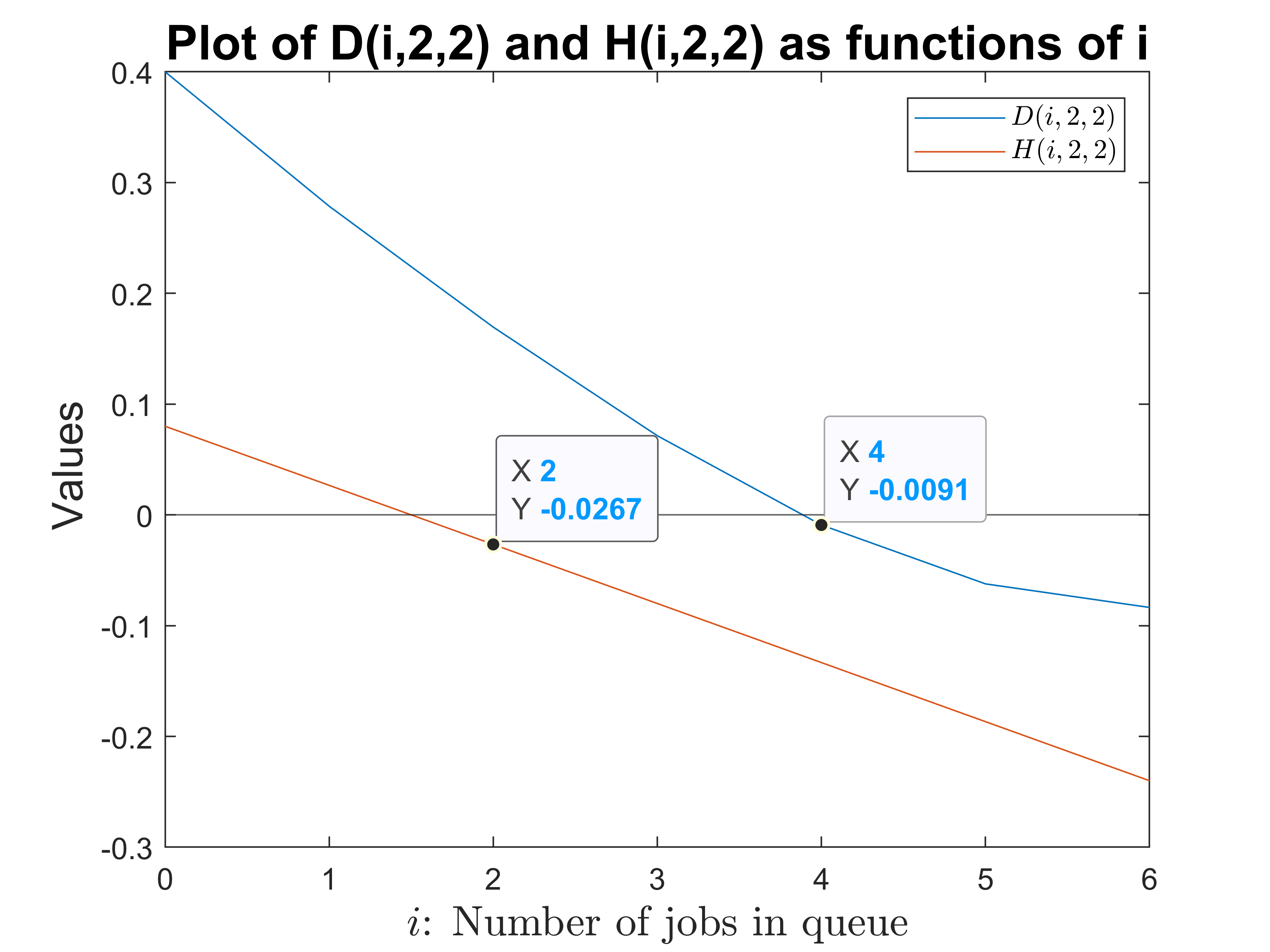}
        \caption{$i_D(2)=4$ and $i_H(2) = 2$ when $h_1 = 0.8$.}
    \end{subfigure}
    \hfill
    \begin{subfigure}[htbp]{0.49\textwidth}
        \centering
        \includegraphics [scale=.6] {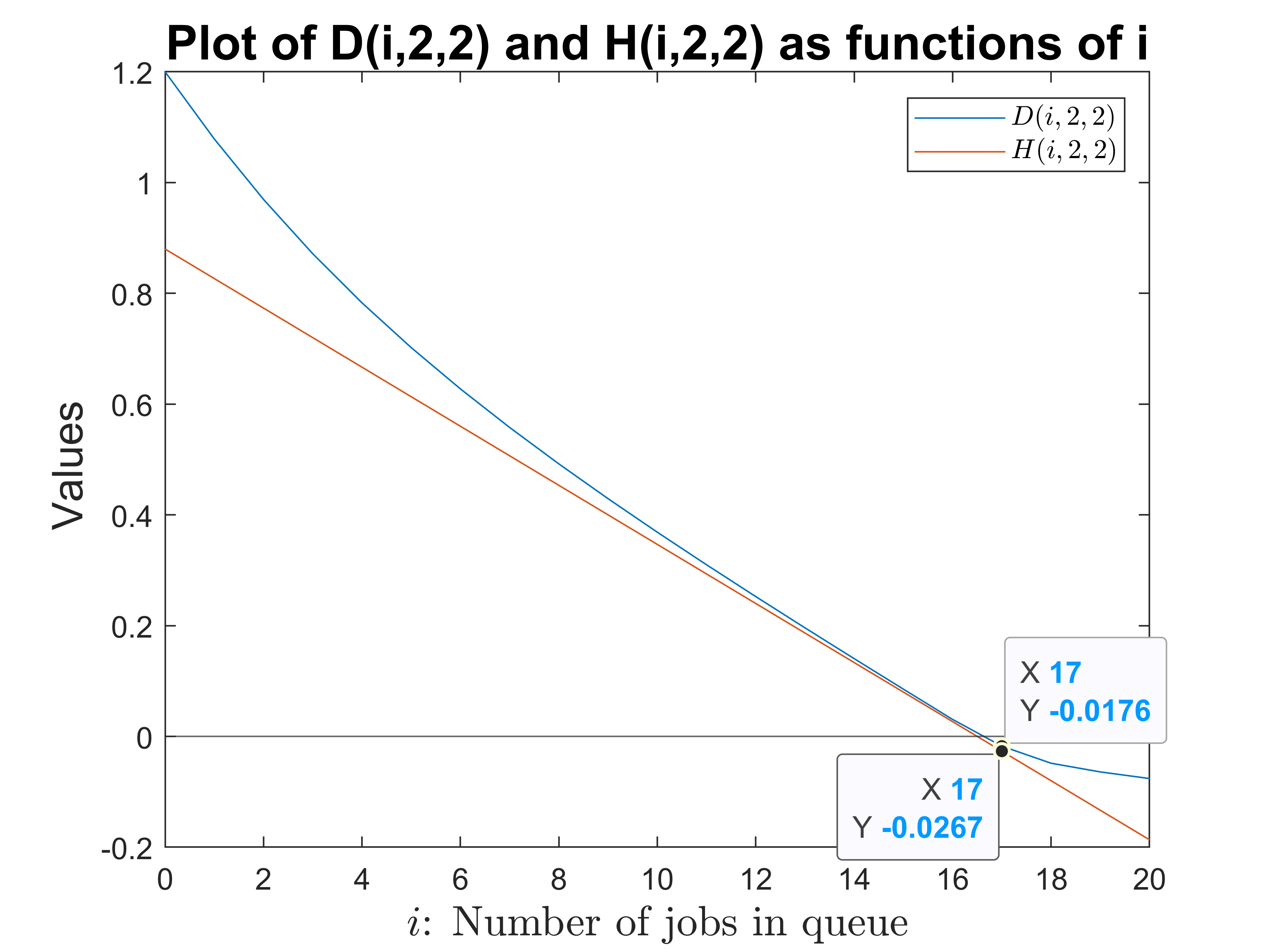}
        \caption{$i_D(2) = i_H(2)= 17$ when $h_1 = 1.6$.}
    \end{subfigure}
    \caption{Plot of $D(i,k,\ell)$ and $H(i,k,\ell)$ as functions of $i$ under parameter configurations in Example \ref{eg:costly_collab} such that $h_0<h_2$.}\label{fig:costly_collab}
\end{figure*}

\begin{figure}[htbp]
\centering
    \includegraphics[scale=0.6]{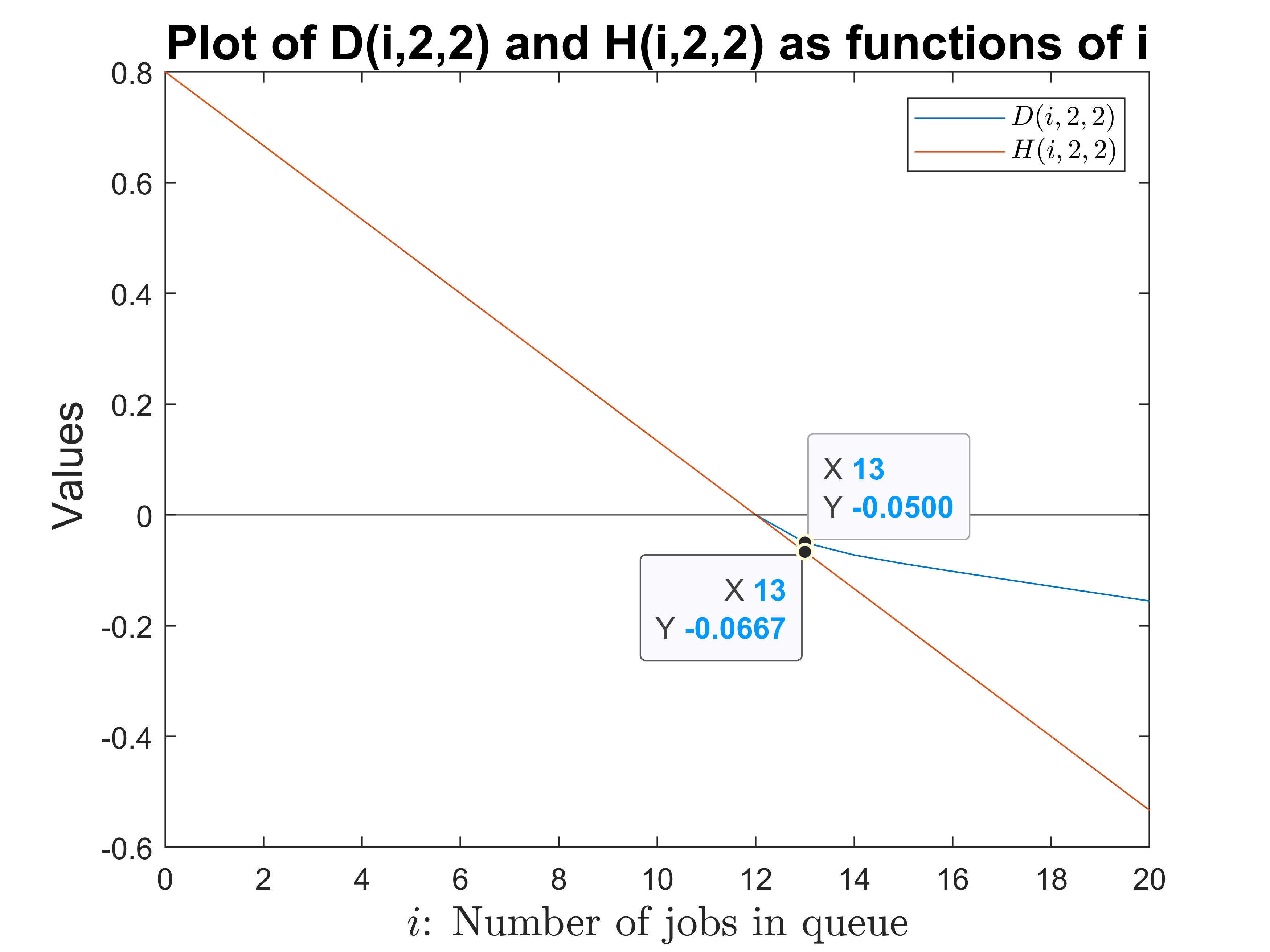}
\caption{Plot of $D(i,k,\ell)$ and $H(i,k,\ell)$ as functions of $i$ under parameter configurations in Example \ref{eg:equally_costly} such that $h_0=h_2$.} \label{Fig:equally_costly}
\end{figure}

\begin{proof}[Proof of Statement \ref{state:thresh_queue_in_collab} in Theorem \ref{thm:heur_lowcost}.]
    For any $k$ such that $\ell = C_1-k \geq C_2$, consider first the case where $h_0 \geq h_2$. We aim to show $i_H(k) \geq i_D(k)$ by proving that $D(i,k,\ell)>0$ implies $H(i,k,\ell) >0$. Suppose $D(i,k,\ell) > 0$ for some $i$; this in turn implies $\frac{h_1}{\mu_1} > \frac{\ell+1}{C_2}\frac{h_2}{\mu_2}$, since otherwise we would have $D(i,k,\ell) \leq 0$ for all $i$ (Corollary \ref{cor:always_action0} or Statement \ref{state:always_neg_queue_in_collab} in Theorem \ref{thm:heur_lowcost}). Referring to the definition of $H$ in \eqref{def:H}, we have $H(i,k,\ell) = (i-y_k)c'+b' \geq D(i,k,\ell) > 0$ by Statement \ref{state:queue_in_collab_costly_queue_ub} in Lemma \ref{lemma:heur}.
    
    If $h_0 \leq h_2$, a symmetric argument establishes $i_H(k) \leq i_D(k)$ by showing that $H(i,k,\ell) > 0$ implies $D(i,k,\ell) > 0$, using Statement \ref{state:queue_in_collab_costly_collab_lb} in Lemma \ref{lemma:heur}.
\end{proof}
Condition \ref{cond:sufficient_equal_thresh} presents two sufficient conditions for $i_D(k) = i_H(k)$ where $\ell=C_1-k\geq C_2$: Inequality \eqref{cond:highcost_queue} applies when $h_0 > h_2$, while inequality \eqref{cond:highcost_collab} applies when $h_0 < h_2$.
This is precisely the content of Statement~\ref{state:equal_thresh_queue_in_collab} in Theorem~\ref{thm:heur_lowcost}, with the proof provided below.
\begin{proof} [Proof of Statement \ref{state:equal_thresh_queue_in_collab} in Theorem \ref{thm:heur_lowcost}.]
    First, consider the case $h_0 > h_2$. Statement \ref{state:thresh_reverse} in Lemma \ref{lemma:heur} establishes $i_D(k) \geq i_H(k)$ since $k$ satisfies \eqref{cond:highcost_queue}, while Statement \ref{state:thresh_queue_in_collab} in Theorem \ref{thm:heur_lowcost} confirms the reverse inequality $i_D(k) \leq i_H(k)$, yielding $i_D(k) = i_H(k)$.
    A symmetric argument applies when $h_0 < h_2$ and if $k$ satisfies \eqref{cond:highcost_collab}.

    Finally, in both cases, the necessity of the condition for achieving $i_D(k) = i_H(k)$ when $k=1$ follows directly from Statement \ref{state:thresh_reverse} in Lemma \ref{lemma:heur}.
\end{proof}
\begin{remark} \label{remark:equal_thresh}
    With explicit expressions, Condition \ref{cond:sufficient_equal_thresh} provides a practical sufficient condition for ensuring the exact estimation of $i_D(k)$, depending on whether $h_0 > h_2$. 
    Notice that $i_H(k) \geq 1$ is a must in order for $i_H(k) = i_D(k)$, considering $i_D(k) \geq 1$. Consequently, Condition \ref{cond:sufficient_equal_thresh} reveals that the right-hand side captures two key factors: 
    \begin{enumerate}
        \item \label{fac:ratio}
        The ratio $\frac{h_2}{h_0}$;
        \item \label{fac:dist_to_integers}
        The distance between the zero of $H(i,k,\ell)$ --- namely, $R_2(k)$ (see \eqref{eq:thresh_queue_in_collab}) --- and its nearest non-negative integers $i_H(1)-1$ and $i_H(1)$.
    \end{enumerate}
   The left-hand side is proportional to a term that decays exponentially with $i_H(k)-k$. Intuitively, as $i_H(k)$ increases, the condition is more likely to hold, leading to a more accurate approximation of $i_D(k)$, provided that other factors, especially the two factors on the right-hand side, remain unchanged.
\end{remark}
Example \ref{eg:costly_queue} below illustrates Statements \ref{state:thresh_queue_in_collab} and \ref{state:equal_thresh_queue_in_collab} in Theorem \ref{thm:heur_lowcost} both when \eqref{cond:highcost_queue} in Condition \ref{cond:sufficient_equal_thresh} is satisfied or when it is violated for the case where $h_0 > h_2$. Similarly, for Example \ref{eg:costly_collab} for the case when $h_0 < h_2$. Example \ref{eg:equally_costly} considers the special case $h_0 = h_2$.
Figures \ref{fig:costly_queue}--\ref{Fig:equally_costly} visually depict these examples, and compare $H(i,k,\ell)$ with $D(i,k,\ell)$, supporting Statement \ref{state:asymp} in Lemma \ref{lemma:heur}.

\begin{example} [$h_0>h_2$] \label{eg:costly_queue}
    Let $C_1 = 4$, $C_2 = 2$, $\mu_1 = 1$, $\mu_2 = 1.5$, $h_0 = 2$, $h_1 = 2$, $h_2 = 1$, satisfying $\frac{h_1}{\mu_1}> \frac{h_2}{\mu_2}$. Under these parameters, we obtain $-\frac{b'}{c'}+y_2 = 3.5$, so that (from \eqref{def:iH_lowcost}) $i_H(2) =\lceil3.5\rceil = 4$, and $i_D(2)=3$. Note, inequality \eqref{cond:highcost_queue} in Condition \ref{cond:sufficient_equal_thresh} is not satisfied (by $k=2$). 
    On the other hand, if we set $h_1 = 8$ while keeping all other parameters unchanged, we obtain $-\frac{b'}{c'}+y_2 = 12.5 > 3.5$. Notably, both $y_2$ and $r$ on the left-hand side, as well as the entire right-hand side of \eqref{cond:highcost_queue}, remain unchanged; however, the inequality now holds. Consequently, we find $i_D(2) = i_H(2)= 13 = \lceil12.5\rceil$. In both cases, the inequality $i_D(k) \leq i_H(k)$ holds, supporting Statement \ref{state:thresh_queue_in_collab} in Theorem \ref{thm:heur_lowcost}.
\end{example}
\begin{example} [$h_0<h_2$] \label{eg:costly_collab}
    Let $C_1 = 4$, $C_2 = 2$, $\mu_1 = 1$, $\mu_2 = 1.5$, $h_0 = 0.16$, $h_1 = 0.8$, $h_2 = 0.4$, satisfying $\frac{h_1}{\mu_1}> \frac{h_2}{\mu_2}$. Under these parameters, we have $-\frac{b'}{c'}+y_2 = 1.5$ so that $i_H(2) =  \lceil1.5\rceil = 2 $, and $i_D(2)=4$. Moreover, inequality \eqref{cond:highcost_collab} in Condition \ref{cond:sufficient_equal_thresh} is not satisfied (by $k=2$). 
    However, setting $h_1 = 1.6$ while keeping all other parameters unchanged yields $-\frac{b'}{c'}+y_2 = 16.5 > 1.5$. Again, all $y_2$, $r$, and $\frac{r_2}{1-q_2}$ on the left-hand side, as well as the entire right-hand side of \eqref{cond:highcost_collab}, remain the same; yet the inequality now holds. Consequently, we get $i_D(2) = i_H(2) = \lceil16.5\rceil= 17 $. In both cases, the inequality $i_D(k) \geq i_H(k)$ holds, supporting \ref{state:thresh_queue_in_collab} in Theorem \ref{thm:heur_lowcost}.
\end{example}
\begin{example} [$h_0 = h_2$] \label{eg:equally_costly}
    Let $C_1 = 4$, $C_2 = 2$, $\mu_1 = 1$, $\mu_2 = 1.5$, $h_0 = 0.2$, $h_1 = 1$, $h_2 = 0.2$, satisfying $\frac{h_1}{\mu_1}> \frac{h_2}{\mu_2}$. Under these parameters, we have $i_D(2) = i_H (2) = 13$, supporting the conclusion from Statement \ref{state:thresh_queue_in_collab} in Theorem \ref{thm:heur_lowcost}. 
\end{example}

We conclude this subsection by proving the last statement in Theorem \ref{thm:heur_lowcost}, which quantifies the increment $i_D(k) - i_D(k-1)$ under Condition \ref{cond:sufficient_equal_thresh}.
\begin{proof} [Proof of Statement \ref{state:sequence_iD} in Theorem \ref{thm:heur_lowcost}.]
    Suppose first $h_0 = h_2$. From Statement \ref{state:thresh_queue_in_collab} in Theorem \ref{thm:heur_lowcost}, we have 
    \begin{align}
        i_H(k-1) &= i_D(k-1),\label{eq:equal_thresh_k-1} \\ 
        i_H(k) &= i_D(k) \geq 1, \label{eq:equal_thresh_k}
    \end{align}where the inequality in \eqref{eq:equal_thresh_k} is implied by Statement \ref{state:always_neg_queue_in_collab} in Theorem \ref{thm:heur_lowcost} since $\frac{h_1}{\mu_1} > \frac{\ell+1}{C_2}\frac{h_2}{\mu_2}$.
    Equation \eqref{eq:equal_thresh_k} leads directly to Statement \ref{state:sequence_iH} in Lemma \ref{lemma:heur}:
    \begin{align}
        i_H(k) = i_H(k-1)+1. \label{eq:sequence_iH}
    \end{align}
    Combining \eqref{eq:equal_thresh_k-1}, \eqref{eq:equal_thresh_k}, and \eqref{eq:sequence_iH} completes the proof in this case.

    Consider now $h_0>h_2$ (resp. $h_0 < h_2$), where inequality \eqref{cond:highcost_queue} (resp. \eqref{cond:highcost_collab}) in Condition \ref{cond:sufficient_equal_thresh} holds for $k$.
    Equation \eqref{eq:equal_thresh_k} also holds by Statement \ref{state:equal_thresh_queue_in_collab} in Theorem \ref{thm:heur_lowcost} under the corresponding condition for $k$, together with Statement \ref{state:always_neg_queue_in_collab} in Lemma \ref{lemma:heur} again. Consequently, both results in Statement \ref{state:cond_k-1} of Lemma \ref{lemma:heur} follow: Equation \eqref{eq:sequence_iH} holds, and $k-1$ also satisfies \eqref{cond:highcost_queue} (resp. \eqref{cond:highcost_collab}) in Condition \ref{cond:sufficient_equal_thresh}.
    It remains to verify \eqref{eq:equal_thresh_k-1}. If $i_D(k-1) = 0$, then $i_H(k-1)=0$ by Statement \eqref{state:always_neg_queue_in_collab} in Theorem \ref{thm:heur_lowcost}; if $i_D(k-1) \geq 1$, \eqref{eq:equal_thresh_k-1} follows from Statement \ref{state:equal_thresh_queue_in_collab} in the same theorem under the corresponding condition for $k-1$.
\end{proof}

\subsection{Collaborative service costs higher than independent \texorpdfstring{$\big(\frac{h_1}{\mu_1} \leq \frac{h_2}{\mu_2}\big)$}{(h1/mu1 <= h2/mu2)}}
If $\frac{h_1}{\mu_1} \leq \frac{h_2}{\mu_2}$, the only remaining scenario for estimating $D(i,k,\ell)$ is when $\mu_1 < \mu_2$ and $\ell < C_2$. 
In parallel to Statement \ref{state:heur_bd} in Theorem \ref{thm:heur_lowcost}, Statement \ref{state:heur_bd2} in Theorem \ref{thm:heur_highcost} establishes a two-sided bound comparing the estimated threshold $\tilde{i}_H(\ell)$ with the actual threshold $\tilde{i}_D(\ell)$ in this setting --- both thresholds are infinite otherwise. The proof is given below.

\begin{figure*} 
    \centering
    \begin{subfigure}[b]{0.49\textwidth}
        \centering
        \includegraphics [scale=.6] {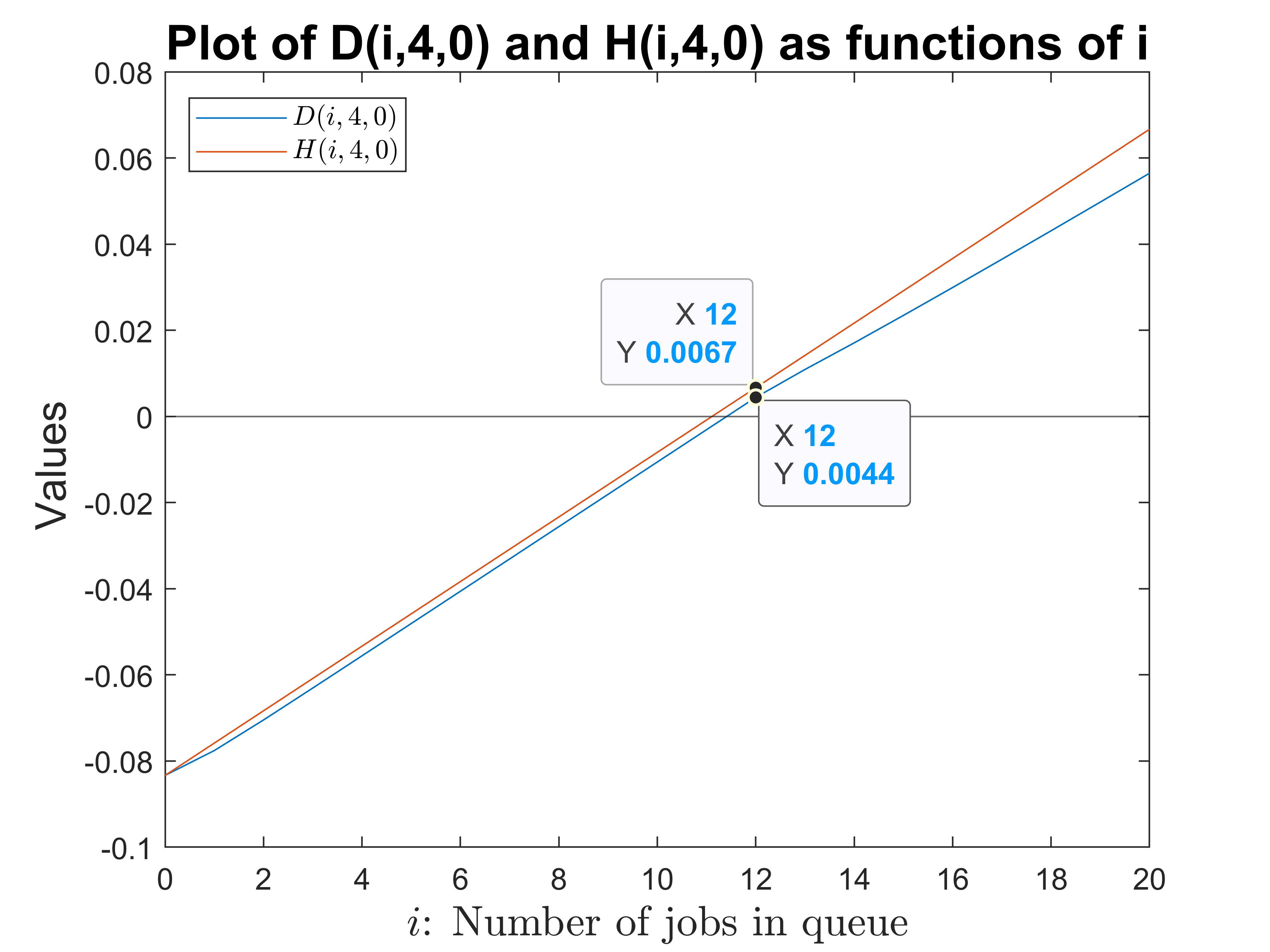}
        \caption{Example \ref{eg:best2}.}
    \end{subfigure}
    \hfill
    \begin{subfigure}[b]{0.49\textwidth}
        \centering
        \includegraphics [scale=.6] {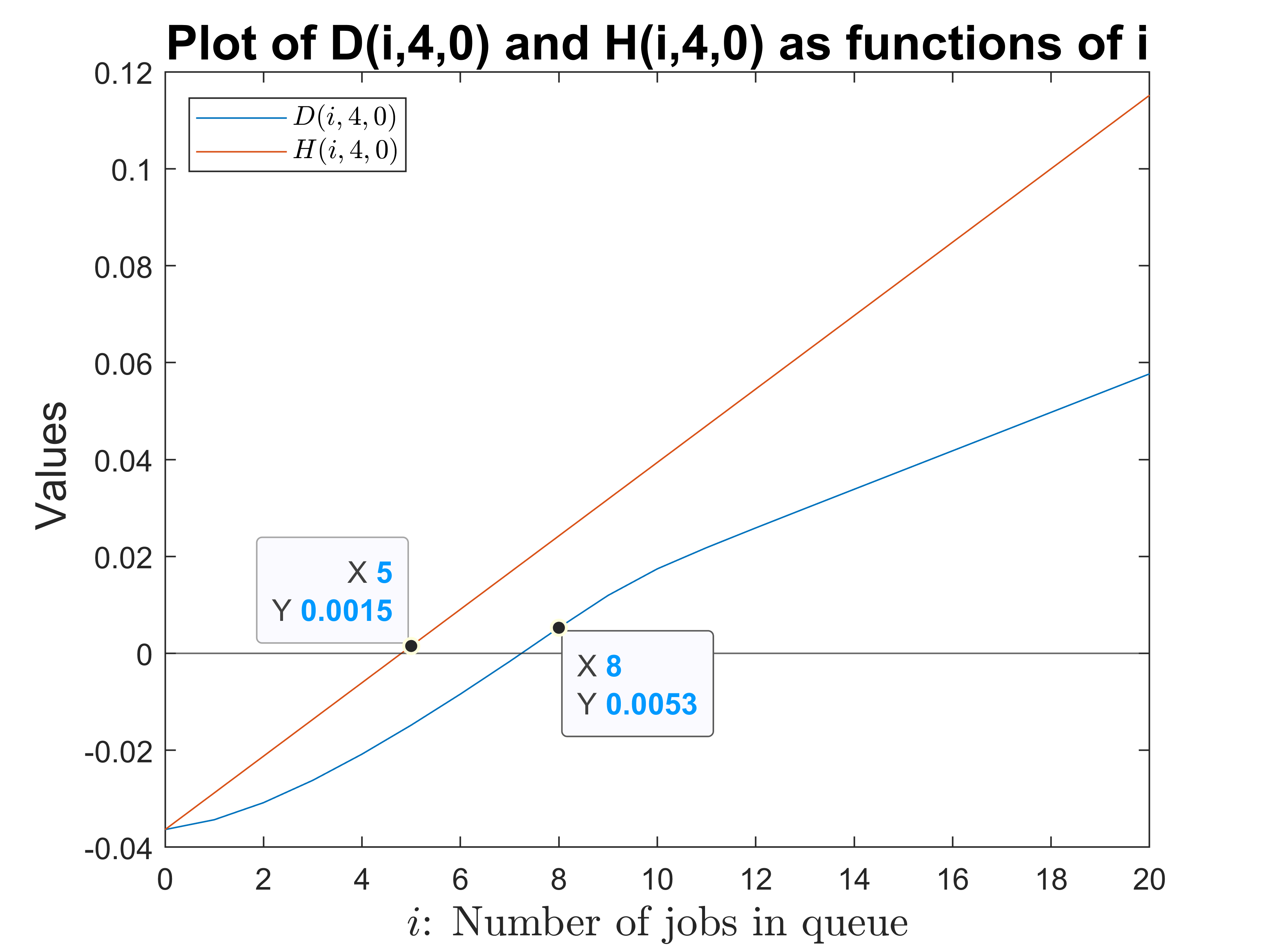}
        \caption{Example \ref{eg:worst2}.}
    \end{subfigure}
    \caption{Plot of $D(i,k,\ell)$ and $H(i,k,\ell)$ as functions of $i$ under parameter configurations in Examples \ref{eg:best2} and \ref{eg:worst2}.}\label{fig:eg_highcost_no_wait}
\end{figure*}

\begin{proof}[Proof of Statement \ref{state:heur_bd2} in Theorem \ref{thm:heur_highcost}.]
    Suppose $\mu_1 < \mu_2$ and consider $\ell < C_2$.
    To see the direction that $\tilde{i}_H(\ell) \leq \tilde{i}_D(\ell)$, notice that $H(i,k,\ell)<0$ implies $D(i,k,\ell) \leq H(i,k,\ell)< 0$ by Statement \ref{state:diff_ub_larger_mu2} in Lemma \ref{silemma:diff}. 
    
    To establish the remaining part of the result, namely $\tilde{i}_D(\ell) \leq \tilde{i}_H(\ell) + (C_1-1)$, we show that $D(i,k,\ell) < 0$ implies $H\big(i-(C_1-1)\big) = \big(i-(C_1-1)\big)c+b < 0$, or equivalently, $i < \frac{-b}{c} + (C_1-1)$ since $c>0$ (applying $\mu_2 > \mu_1$ in \eqref{def:b_and_c}). 
    As in the proof of Statement \ref{state:heur_bd} in Theorem \ref{thm:heur_lowcost}, applying Statement \ref{state:heur_linear_lb} from Lemma \ref{lemma:heur} yields $0 > D(i,k,\ell) \geq (i-z_\ell)c+b$, provided that $D(i,k,\ell) < 0$. Consequently, we obtain $i < \frac{-b}{c} + z_\ell$. The result follows from the fact that $z_\ell \leq C_1-1$, as stated in Statement \ref{state:z_ub} in the same lemma.

    In other regions of the parameter and state spaces --- specifically, when either $\mu_1<\mu_2$ with $\ell \geq C_2$ or $\mu_1 \geq \mu_2$ --- Corollary \ref{cor:always_action0} implies $D(i,k,\ell) \leq 0$ for all $i$, leading to $\tilde{i}_D(\ell) = \infty$. This aligns with the definition $\tilde{i}_H(\ell) = \infty$ in \eqref{def:iH_highcost}.
\end{proof}
We provide examples attaining the equalities in Statement \ref{state:heur_bd2} in Theorem \ref{thm:heur_highcost} as well, where the visualization is depicted in Figure \ref{fig:eg_highcost_no_wait}.
\begin{example}\label{eg:best2}
    Let $C_1 = 4$, $C_2 = 2$, $\mu_1 = 3$, $\mu_2 = 30$, $h_0 = 0.1$, $h_1 = 1$, $h_2 = 12.5$ satisfying $\frac{h_1}{\mu_1}<  \frac{h_2}{\mu_2}$. We have $\tilde{i}_D(0) = \tilde{i}_H(0)=12$.
\end{example}
\begin{example}\label{eg:worst2}
    Let $C_1 = 4$, $C_2 = 2$, $\mu_1 = 3$, $\mu_2 = 3.3$, $h_0 = 1$, $h_1 = 1$, $h_2 = 1.22$ satisfying $\frac{h_1}{\mu_1}<  \frac{h_2}{\mu_2}$. We have $\tilde{i}_D(0)=8$ and $\tilde{i}_H(0) = 5$, so that $\tilde{i}_D(0) = \tilde{i}_H(0)+(C_1-1)$.
\end{example}
\begin{remark}
    Experiments indicate that the worst-case scenario primarily arises when $\mu_1$ and $\mu_2$ are close; see Example \ref{eg:worst2}.
\end{remark}

Finally, Statement \ref{state:thresh_mono_highcost} in Theorem \ref{thm:heur_highcost}, which, like Statement \ref{state:thresh_mono} in Theorem \ref{thm:heur_lowcost}, follows directly from Proposition \ref{prop:mono_ell} as well.

\subsection{No queueing effects for all states \texorpdfstring{$(C_2\geq C_1)$}{(C2 >= C1)}}
We conclude this section by examining the corresponding results for the case without queueing effects across all states (i.e., $C_2 \geq C_1$), as special cases of those derived for $C_2 < C_1$ with $\ell < C_2$.
In particular, if $C_2 \geq C_1 = 1$, the approximation becomes exact: the estimated threshold coincides with the actual threshold in all regions of the parameter space where the latter is finite, and both thresholds are infinite otherwise.

\begin{proof} [Proof of Corollary \ref{cor:heur_enough_Type_II}.]
    Since $x \in \X_{diff}$ implies that $k = C_1 - \ell \geq 1$, it follows that $\ell =C_1-k\leq C_1 - 1 \leq C_2 -1$ since $C_2 \geq C_1$, ensuring $\ell < C_2$. Therefore, the Bellman equations for $C_2 \geq C_1$ coincide with those when $C_2 < C_1$ with $\ell <C_2$. 
    Consequently, the proofs for $C_2 \geq C_1$ follow directly from those for $C_2 < C_1$ with $\ell < C_2$ under the same parameter assumptions. Specifically,
    Statements \ref{state:thresh_mono_enough_Type_II} and \ref{state:heur_bd_enough_Type_II} are implied by Statements \ref{state:thresh_mono} and \ref{state:heur_bd} in Theorem \ref{thm:heur_lowcost}, respectively, and Statements \ref{state:thresh_mono_ell_enough_Type_II} and \ref{state:heur_bd2_enough_Type_II} follow from Statements \ref{state:thresh_mono_highcost} and \ref{state:heur_bd2} in Theorem \ref{thm:heur_highcost}, respectively.
\end{proof}

\section{Numerical Analysis} \label{sec:numerical}
We evaluate the performance of our proposed heuristics through a numerical study, comparing the value functions --- i.e., total expected costs --- of states $(i_0, k_0, \ell_0)$, under various threshold heuristics (including ours) and an optimal policy. The tested parameters are listed in Table \ref{table:params}, covering all combinations of configurations. Note that the optimal policy remains unchanged when all holding costs are multiplied by a positive constant (while service rates remain fixed); therefore, without loss of generality, we set $h_1 = 1$. For the same reason, we fix $\mu_1 = 10$ and allow $\mu_2$ to vary.

\begin{table}[htbp]
\caption{Parameter configurations for numerical test.}
\label{table:params}
\centering
\begin{tabular}{cc}
\hline
\multicolumn{1}{c}{\textbf{Parameters}} & \multicolumn{1}{c}{\textbf{Values}} \\ \hline
$h_0$                               & 0.01, 0.02, 0.05, 0.1, 0.2, 0.5, 1  \\
$h_1$                               & 1                                   \\
$h_2$                               & 0.1, 0.2, 0.5, 1, 1.5, 2            \\
$\mu_1$                              & 10                              \\
$\mu_2$                              & 4, 5, 6, 8, 10, 12, 15, 20, 25    \\ \hline    
\end{tabular}
\end{table}

The performance is measured by the relative error defined as follows.
\begin{align*}
    err = \frac{v^{\pi}(i_0,k_0,\ell_0) - v(i_0,k_0,\ell_0)}{v(i_0,k_0,\ell_0)}, 
\end{align*}where $v^{\pi}$ is the value function under any policy $\pi$ and $v$ denotes the optimal value function. To distinguish our heuristic from others that are commonly adopted in practice, we refer to it as $\pi'$. We compute total expected costs and relative errors for all states $(i_0,k_0,\ell_0) \in \X_D$ where $i_0 = 20$ and $30$. The results, summarized in Tables \ref{table:cmp_state20_larger_mu1}--\ref{table:cmp_state30_highcost}, show that $\pi'$ achieves near-optimal accuracy, with an average cost within 0.5\% of the optimal, and remains robust across parameter variations. On the other hand, alternative threshold heuristics (defined later) exhibit high sensitivity to changes in system parameters and can incur costs up to 100\% higher than the optimal.

Recall the definition of the difference $D(i,k,\ell)$ for $(i,k,\ell) \in \tilde{\X}_{diff}$ in \eqref{def:D}, the signs of $D(i-1,k,\ell)$ and $D(i-1,k+1,\ell-1)$ determine where to serve the job if a service is completed first at Station 1 and Station 2, respectively. 
Therefore, without loss of generality, we proceed with the numerical analysis under the assumption that a Station 1 job is completed first, i.e. $n=1$ (recall the definition of $n$ following the action set $\A(i,k,\ell)$ in \eqref{def:action_set}).
A symmetric analysis applies for $n=2$.

\subsection{Independent service costs higher than collaborative \texorpdfstring{$\big(\frac{h_1}{\mu_1} > \frac{h_2}{\mu_2}\big)$}{(h1/mu1 > h2/mu2)}}
When $\frac{h_1}{\mu_1} > \frac{h_2}{\mu_2}$, recall the actual and estimated integer thresholds are defined in \eqref{def:iD_lowcost} and \eqref{def:iH_lowcost}, respectively, for $(i,k,\ell) \in \tilde{\X}_{diff}$. Assuming $n=1$, both an optimal policy and $\pi'$ are therefore of \textbf{collaborative threshold type with threshold} $T_D(k)$ and $T_H(k)$, respectively, where $T_D(k) = i_D(k) + 1$ and $T_H(k) = i_H(k) + 1$.
For comparison, we evaluate other benchmark policies $\pi_j$, also of \textbf{collaborative threshold type with threshold $T_j(k)$} for $j = 1,2, \ldots,4$. The results are summarized in Tables \ref{table:cmp_state20_larger_mu1}--\ref{table:cmp_state30_larger_mu2}.
The action taken at each state $x = (i,k,\ell) \in \X_D$ (so $i \geq 1$) under Policy $\pi_j$ is defined as follows. Recall $k+\ell=C_1$.

\begin{enumerate} [leftmargin=3cm, label=Policy $\pi_\arabic*$.]
    \item Always go to Station 1: $a=0$, $\forall n = 1,2$.
    \item Go to Station 1 if and only if the number of jobs waiting in queue exceeds the fixed threshold $10$: $a = \mathbb{I}\{i \leq 10\}$, $\forall n = 1,2$.
    \item Always go to Station 2: $a=1$, $\forall n = 1,2$.
    \item Go to Station 2 if and only if no waiting is required to commence Station 2 service: $a = \mathbb{I}\{\ell<C_2\}$ if $n=1$; and $a = \mathbb{I}\{\ell \leq C_2\}$ if $n=2$.
\end{enumerate}

\begin{table}[htbp]
    \caption{Comparison of relative errors $(\%)$ for states $(20,k_0,\ell_0)$ when $\frac{h_1}{\mu_1} > \frac{h_2}{\mu_2}$ and $\mu_1\geq\mu_2$.}
    \centering
    \renewcommand{\arraystretch}{1.0}
    \label{table:cmp_state20_larger_mu1}
    \begin{tabular}{|c|r|r|r|r|r|r|r|}
        \hline
        \multicolumn{2}{|c|}{} & \multicolumn{1}{c|}{$C_1 = 2$} & \multicolumn{2}{c|}{$C_1 = 3$} & \multicolumn{3}{c|}{$C_1 = 4$} \\ \cline{3-8}
        \multicolumn{2}{|c|}{} & $C_2 = 1$ & $C_2 = 1$ & $C_2 = 2$ & $C_2 = 1$ & $C_2 = 2$ & $C_2 = 3$ \\
        \hline

        \multirow{3}{*}{\textbf{Policy $\pi'$}} & Max error & 0.19 & 1.47 & 0.69 & 3.33 & 0.99 & 1.13 \\
                                               & Avg error & 0.01 & 0.09 & 0.06 & 0.28 & 0.10 & 0.11 \\
                                               & Std error & 0.03 & 0.27 & 0.14 & 0.68 & 0.22 & 0.23 \\
        \hline

        \multirow{3}{*}{\textbf{Policy $\pi_1$}} & Max error & 235.16 & 147.05 & 379.17 & 96.42 & 277.65 & 451.97 \\
                                                & Avg error & 25.88  & 16.45  & 48.76  & 11.42  & 36.21  & 62.95 \\
                                                & Std error & 38.58  & 23.38  & 63.31  & 15.49  & 46.22  & 73.49 \\
        \hline

        \multirow{3}{*}{\textbf{Policy $\pi_2$}} & Max error & 109.30 & 65.51 & 181.64 & 75.62 & 127.42 & 218.39 \\
                                                & Avg error & 22.48  & 24.07  & 23.58  & 26.98  & 18.25  & 28.12 \\
                                                & Std error & 15.32  & 14.44  & 28.54  & 18.28  & 19.09  & 35.19 \\
        \hline

        \multirow{3}{*}{\textbf{Policy $\pi_3$}} & Max error & 263.48 & 363.05 & 140.66 & 430.67 & 180.97 & 94.70 \\
                                                & Avg error & 66.48  & 99.52  & 27.20  & 122.72 & 40.56  & 14.59 \\
                                                & Std error & 58.63  & 77.97  & 31.43  & 89.53  & 40.47  & 20.33 \\
        \hline

        \multirow{3}{*}{\textbf{Policy $\pi_4$}} & Max error & 101.56 & 88.43 & 120.08 & 68.03 & 137.02 & 114.47 \\
                                                & Avg error & 9.66   & 6.24   & 14.06  & 4.58   & 13.11  & 14.31 \\
                                                & Std error & 15.50  & 12.56  & 19.23  & 9.45   & 21.06  & 17.81 \\
        \hline
    \end{tabular}
\end{table}

\begin{table}[htbp]
    \caption{Comparison of relative errors $(\%)$ for states $(20,k_0,\ell_0)$ when $\frac{h_1}{\mu_1} > \frac{h_2}{\mu_2}$ and $\mu_1<\mu_2$.}
    \centering
    \label{table:cmp_state20_larger_mu2}
    \begin{tabular}{|c|r|r|r|r|r|r|r|}
        \hline
        \multicolumn{2}{|c|}{} & \multicolumn{1}{c|}{$C_1 = 2$} & \multicolumn{2}{c|}{$C_1 = 3$} & \multicolumn{3}{c|}{$C_1 = 4$} \\ \cline{3-8}
        \multicolumn{2}{|c|}{} & $C_2 = 1$ & $C_2 = 1$ & $C_2 = 2$ & $C_2 = 1$ & $C_2 = 2$ & $C_2 = 3$ \\
        \hline

        \multirow{3}{*}{\textbf{Policy $\pi'$}} & Max error & 1.35 & 5.73 & 1.61 & 5.06 & 3.65 & 2.45 \\
                                               & Avg error & 0.08 & 0.30 & 0.12 & 0.34 & 0.28 & 0.14 \\
                                               & Std error & 0.23 & 0.92 & 0.30 & 0.92 & 0.66 & 0.39 \\
        \hline

        \multirow{3}{*}{\textbf{Policy $\pi_1$}} & Max error & 723.99 & 500.79 & 1078.51 & 371.49 & 822.65 & 1255.11 \\
                                                & Avg error & 88.39  & 58.56  & 127.38  & 42.28  & 96.82  & 142.77 \\
                                                & Std error & 100.69 & 70.01  & 140.37  & 52.61  & 108.51 & 153.55 \\
        \hline

        \multirow{3}{*}{\textbf{Policy $\pi_2$}} & Max error & 356.55 & 238.55 & 533.24 & 170.89 & 398.28 & 621.66 \\
                                                & Avg error & 52.25  & 39.35  & 66.65  & 33.81  & 48.94  & 73.26 \\
                                                & Std error & 45.24  & 28.60  & 67.95  & 22.60  & 50.44  & 75.62 \\
        \hline

        \multirow{3}{*}{\textbf{Policy $\pi_3$}} & Max error & 73.52  & 134.47 & 36.17  & 188.21 & 66.05  & 23.79 \\
                                                & Avg error & 27.86  & 51.36  & 9.97   & 70.51  & 19.38  & 4.66 \\
                                                & Std error & 22.58  & 41.13  & 10.26  & 56.39  & 19.40  & 5.86 \\
        \hline

        \multirow{3}{*}{\textbf{Policy $\pi_4$}} & Max error & 201.34 & 226.78 & 200.95 & 214.64 & 271.92 & 184.83 \\
                                                & Avg error & 13.76  & 14.99  & 14.66  & 14.05  & 19.95  & 12.75 \\
                                                & Std error & 30.19  & 32.35  & 27.80  & 29.89  & 36.81  & 23.18 \\
        \hline
    \end{tabular}
\end{table}

\begin{table}[htbp]
    \caption{Comparison of relative errors $(\%)$ for states $(30,k_0,\ell_0)$ when $\frac{h_1}{\mu_1} > \frac{h_2}{\mu_2}$ and $\mu_1 \geq \mu_2$.}
    \centering
    \label{table:cmp_state30_larger_mu1}
    \begin{tabular}{|c|r|r|r|r|r|r|r|}
        \hline
        \multicolumn{2}{|c|}{} & \multicolumn{1}{c|}{$C_1 = 2$} & \multicolumn{2}{c|}{$C_1 = 3$} & \multicolumn{3}{c|}{$C_1 = 4$} \\ \cline{3-8}
        \multicolumn{2}{|c|}{} & $C_2 = 1$ & $C_2 = 1$ & $C_2 = 2$ & $C_2 = 1$ & $C_2 = 2$ & $C_2 = 3$ \\
        \hline

        \multirow{3}{*}{\textbf{Policy $\pi'$}} & Max error & 0.13 & 1.00 & 0.47 & 2.28 & 0.68 & 0.78 \\
                                               & Avg error & 0.00 & 0.06 & 0.03 & 0.18 & 0.06 & 0.07 \\
                                               & Std error & 0.02 & 0.18 & 0.09 & 0.46 & 0.15 & 0.15 \\
        \hline

        \multirow{3}{*}{\textbf{Policy $\pi_1$}} & Max error & 199.32 & 127.67 & 345.30 & 87.88 & 258.57 & 428.48 \\
                                                & Avg error & 20.53  & 13.48  & 40.97  & 9.68   & 30.73  & 55.24 \\
                                                & Std error & 32.30  & 19.80  & 58.31  & 13.14  & 42.89  & 71.65 \\
        \hline

        \multirow{3}{*}{\textbf{Policy $\pi_2$}} & Max error & 125.76 & 76.75 & 224.13 & 52.33 & 163.50 & 280.02 \\
                                                & Avg error & 17.96  & 17.14  & 25.60  & 17.62  & 19.14  & 33.52 \\
                                                & Std error & 17.50  & 10.64  & 36.69  & 10.89  & 25.53  & 46.65 \\
        \hline

        \multirow{3}{*}{\textbf{Policy $\pi_3$}} & Max error & 302.39 & 434.99 & 175.94 & 539.87 & 227.08 & 125.95 \\
                                                & Avg error & 83.16  & 127.23 & 37.06  & 160.02 & 55.68  & 21.70 \\
                                                & Std error & 68.35  & 94.10  & 39.63  & 111.64 & 51.82  & 27.67 \\
        \hline

        \multirow{3}{*}{\textbf{Policy $\pi_4$}} & Max error & 80.14  & 70.83  & 102.79 & 54.05  & 119.13 & 101.32 \\
                                                & Avg error & 8.80   & 5.12   & 13.18  & 3.49   & 11.16  & 13.97 \\
                                                & Std error & 12.41  & 9.63   & 17.15  & 7.09   & 18.18  & 16.96 \\
        \hline
    \end{tabular}
\end{table}

\begin{table}[htbp]
    \caption{Comparison of relative errors $(\%)$ for states $(30,k_0,\ell_0)$ when $\frac{h_1}{\mu_1} > \frac{h_2}{\mu_2}$ and $\mu_1 < \mu_2$.}
    \centering
    \label{table:cmp_state30_larger_mu2}
    \begin{tabular}{|c|r|r|r|r|r|r|r|}
        \hline
        \multicolumn{2}{|c|}{} & \multicolumn{1}{c|}{$C_1 = 2$} & \multicolumn{2}{c|}{$C_1 = 3$} & \multicolumn{3}{c|}{$C_1 = 4$} \\ \cline{3-8}
        \multicolumn{2}{|c|}{} & $C_2 = 1$ & $C_2 = 1$ & $C_2 = 2$ & $C_2 = 1$ & $C_2 = 2$ & $C_2 = 3$ \\
        \hline

        \multirow{3}{*}{\textbf{Policy $\pi'$}} & Max error & 0.90 & 3.90 & 1.09 & 3.52 & 2.53 & 1.69 \\
                                               & Avg error & 0.05 & 0.20 & 0.07 & 0.22 & 0.18 & 0.09 \\
                                               & Std error & 0.15 & 0.61 & 0.20 & 0.62 & 0.45 & 0.26 \\
        \hline

        \multirow{3}{*}{\textbf{Policy $\pi_1$}} & Max error & 638.11 & 456.80 & 998.90 & 347.19 & 780.56 & 1202.72 \\
                                                & Avg error & 81.88  & 54.78  & 122.72  & 39.92  & 93.82  & 142.47 \\
                                                & Std error & 90.14  & 63.90  & 135.04  & 48.50  & 106.24 & 153.97 \\
        \hline

        \multirow{3}{*}{\textbf{Policy $\pi_2$}} & Max error & 423.35 & 296.69 & 664.52 & 220.93 & 512.76 & 800.42 \\
                                                & Avg error & 61.01  & 43.41  & 85.67  & 34.85  & 63.90  & 98.20 \\
                                                & Std error & 56.73  & 37.38  & 88.54  & 27.13  & 68.01  & 101.88 \\
        \hline

        \multirow{3}{*}{\textbf{Policy $\pi_3$}} & Max error & 76.43  & 143.88 & 38.18  & 205.78 & 71.38  & 25.15 \\
                                                & Avg error & 32.12  & 60.29  & 12.47  & 84.42  & 24.37  & 6.41 \\
                                                & Std error & 23.34  & 43.52  & 11.24  & 61.01  & 21.52  & 6.77 \\
        \hline

        \multirow{3}{*}{\textbf{Policy $\pi_4$}} & Max error & 164.77 & 191.79 & 169.16 & 184.90 & 233.71 & 157.40 \\
                                                & Avg error & 10.54  & 11.56  & 12.23  & 10.80  & 16.63  & 11.28 \\
                                                & Std error & 25.28  & 27.66  & 24.88  & 25.73  & 33.36  & 21.26 \\
        \hline
    \end{tabular}
\end{table}

Referring to the definition of a policy of collaborative threshold-type in Definition \ref{def:threshold_policy}, it is evident that $T_1(k) = 1$, $T_2(k) = 11$, $T_3(k) = \infty$, $\forall k=1,2, \ldots,C_1$; these values are fixed and independent of the states. Meanwhile, $T_4(k)$ is defined as follows:
\begin{align}
    T_4(k) = \left\{
    \begin{array}{cc}
        \infty & \text{if } \ell < C_2,\\
        1 & \text{otherwise.}
    \end{array}\right. \label{eq:T4}
\end{align}
This relies on the fact that $a = \mathbb{I}\{\ell<C_2\} = \mathbb{I}\{1 \leq i < T_4(k)\}$ if $n=1$, and $a = \mathbb{I}\{\ell \leq C_2\} = \mathbb{I}\{1 \leq i < T_4(k+1)\}$ if $n=2$, as per the definition of $T_4(k)$.

One key observation from the experiments is that the threshold $i_D(k)$, where $k = 1,2, \ldots,C_1$, varies significantly depending on the relationship between $\ell$ and $C_2$. Specifically, if $n=1$ and $\ell \geq C_2$, all Type-II servers are occupied when deciding where to perform the service. In this scenario, sending a job to Station 2 not only delays service but also prevents the Type-I server from handling new jobs in the queue. As a result, $i_D(k)$ is non-increasing in $\ell$ (or equivalently, non-decreasing in $k$), indicating that more jobs in Station 2 discourage the collaborative service --- an observation that aligns with Statement \ref{state:thresh_mono} in Theorem \ref{thm:heur_lowcost}. Notably, there is a substantial drop in $i_D(k)$ at $\ell = C_2$ compared to $\ell = C_2-1$, i.e., $i_D\big(C_1-(C_2-1)\big) \gg i_D(C_1-C_2)$, unless $i_D\big(C_1-(C_2-1)\big)$ is already small (say equal to $0$) due to the downstream blocking effects.

Failing to consider the dependence of the actual threshold on $\ell$ (or equivalently $k$), the performance of the four benchmark policies of \textbf{collaborative threshold type} depends heavily on how their fixed thresholds $T_j(k)$ compare to the optimal thresholds $T_D(k)$. This explains why, under different scenarios --- including varying parameter configurations, server combinations, and initial states --- any of the first three fixed-threshold policies may perform either best or worst among the four. 
Policy $\pi_4$, on the other hand, accounts for the non-increasing monotonicity of $i_D(k)$ (and hence $T_D(k)$) with respect to $\ell$, particularly the sharp decline at $\ell = C_2$ due to the downstream blocking (see \eqref{eq:T4}). However, regardless of the parameter values, $\pi_4$ imposes an extreme rule by avoiding any idling of available Type-II servers or any downstream blocking.

Following the theoretical analysis of how $i_H(k)$ differs from $i_D(k)$, we now conduct a numerical investigation to assess how these differences affect the values (costs) under $\pi'$ against the optimal. 
The largest relative errors are typically observed when $i_D(k) > i_H(k) = 0$ for some $k$ such that $\ell = C_1 - k \geq C_2$. The parameter configurations that lead to these high errors exhibit a structural pattern, particularly when $\frac{h_1}{\mu_1} > \frac{\ell+1}{C_2}\frac{h_2}{\mu_2}$ and $h_0 < h_2$. The first inequality corresponds to Case \ref{case:pos_to_neg} in Section \ref{sec:heur_queue_in_collab},
yielding that $D(0,k,\ell)>0$ and $i_D(k) \geq 1$. Otherwise ($\frac{h_1}{\mu_1} \leq \frac{\ell+1}{C_2}\frac{h_2}{\mu_2}$), we have $i_H(k) = i_D(k) = 0$ by Statement \ref{state:always_neg_queue_in_collab} in Theorem \ref{thm:heur_lowcost}. 
The second condition, $h_0 < h_2$, also plays a critical role. Recall that when $h_0 \neq h_2$, the condition in Statement \ref{state:equal_thresh_queue_in_collab} of Theorem \ref{thm:heur_lowcost} (see \eqref{cond:highcost_queue} or \eqref{cond:highcost_collab} in Condition \ref{cond:sufficient_equal_thresh}) is more likely to fail for small $i_H(k)$ (see Remark~\ref{remark:equal_thresh} and Examples~\ref{eg:costly_queue}--\ref{eg:costly_collab}). 
When $h_0 > h_2$, Statement \ref{state:thresh_queue_in_collab} in Theorem \ref{thm:heur_lowcost} provides a positive lower bound for $i_H(k)$ --- in particular, $i_H(k) \geq i_D(k) \geq 1$ --- ensuring it cannot be too small. 
When $h_0<h_2$, however, it is possible that $i_H(k)=0$, in which case $i_H(k) \neq i_D(k)$ without needing to verify \eqref{cond:highcost_collab} in Condition \ref{cond:sufficient_equal_thresh}. 

These worst-case errors caused by $i_D(k) > i_H(k) = 0$ --- which arise when $\ell \geq C_2$, $\frac{h_1}{\mu_1} > \frac{\ell+1}{C_2}\frac{h_2}{\mu_2}$ and $h_0 < h_2$ --- are observed under both $\mu_1 < \mu_2$ and $\mu_1 \geq \mu_2$. 
Among the parameter combinations tested with $\mu_1 \geq \mu_2$ and $C_1 = 4$, $C_2 = 3$, the maximum error occurs when $\mu_1 = \mu_2 = 10$, $h_0 = 0.01$, $h_1 = 1$, and $h_2 = 0.5$, at state $(i_0, 1, 3)$ for both $i_0 = 20$ and $i_0 = 30$. In this case, for $k = 2, 3, 4$ (i.e., $\ell < C_2$), if $n=1$, both $i_0 \leq T_D(k)$ and $i_0 \leq T_H(k)$ hold, so $\pi'$ and the optimal policy choose the same action ($a = 1$). The only discrepancy arises at $k = 1$ ($\ell = 3 \geq C_2$), where $i_D(1) = 4$ but $i_H(1) = 0$.
A similar pattern emerges when $\mu_1 < \mu_2$, though the discrepancy becomes even more pronounced. Among all configurations tested with $\mu_1 < \mu_2$ and $C_1 = 4$, $C_2 = 3$, the maximum error again occurs at $(i_0, 1, 3)$ for both $i_0 = 20$ and $i_0 = 30$, when $\mu_1 = 10$, $\mu_2 = 12$, $h_0 = 0.01$, $h_1 = 1$, and $h_2 = 0.5$. As before, the only discrepancy arises at $k = 1$, where $i_D(1) = 9$ and $i_H(1) = 0$; for $k = 2, 3, 4$ and $n=1$, both $i_D(k)$ and $i_H(k)$ (and hence $T_D(k)$ and $T_H(k)$) are infinite (Statement \ref{state:heur_bd} in Theorem \ref{thm:heur_lowcost}), so $\pi'$ again coincides with the optimal, both suggesting $a = 1$.

Note that both the maximum and average errors are smaller when $\mu_1 \geq \mu_2$ compared to when $\mu_1 < \mu_2$. As discussed earlier, this discrepancy primarily arises when $\ell \geq C_2$ with $\frac{h_1}{\mu_1} > \frac{\ell+1}{C_2}\frac{h_2}{\mu_2}$ (i.e., $i_D(k) \geq 1$; see Case \ref{case:pos_to_neg} in Section \ref{sec:heur_queue_in_collab}) and $h_0 < h_2$. 
Under these conditions, $D(i,k,\ell)$ is monotone non-increasing in $i$, and satisfies $D(i,k,\ell) > 0$ if and only if $i < i_D(k)$. 
Furthermore, our approximation $H$ serves as a global lower bound for $D$ as a function of $i$ (Statement \ref{state:queue_in_collab_costly_collab_lb} in Lemma \ref{lemma:heur}), causing $i_H(k)$ to underestimate $i_D(k)$ (Statement \ref{state:thresh_queue_in_collab} in Theorem \ref{thm:heur_lowcost}).
The accuracy of this approximation, however, depends on the relative magnitudes of $\mu_1$ and $\mu_2$. Suppose $n = 1$. 
When $\mu_2 > \mu_1$, the faster collaborative service partially compensates for the downstream blocking, softening its impact on the upstream. This yields a gradual shift in the optimal decision toward independent service and a more gradual decline of $D(i,j,k,\ell)$ in $i$. Consequently, the portion of the graph where $D > 0$ is more curved, reducing the accuracy of the linear approximation $H$.
In contrast, when $\mu_1 \geq \mu_2$, collaborative service is not only slower but also introduces more downstream congestion. This leads to a sharper transition from collaborative to independent service, making $D$ less curved in $i$ and thus more tightly approximated by the linear bound $H$.

A final observation is that the relative error in costs for states $(i_0,k_0,\ell_0)$ is smaller for $i_0 = 30$ compared to $i_0=20$, indicating improved performance under $\pi'$ in scenarios with longer initial queues.
This is because $\pi'$ is designed as a \textbf{collaborative threshold type policy} that captures key structural properties of the optimal policy via adaptive thresholds, ensuring robustness. Even in parameter regimes where the approximation may be less accurate --- e.g., when $\frac{h_1}{\mu_1} > \frac{\ell+1}{C_2}\frac{h_2}{\mu_2}$ and $h_0 < h_2$ --- a high initial queue diminishes the long-term impact of occasional suboptimal actions taken during periods of light upstream congestion, rendering them negligible over time.

\subsection{Collaborative service costs higher than independent \texorpdfstring{$\big(\frac{h_1}{\mu_1} \leq \frac{h_2}{\mu_2}\big)$}{(h1/mu1 < h2/mu2)}}

When $\frac{h_1}{\mu_1} \leq \frac{h_2}{\mu_2}$, the actual and estimated integer thresholds are now given in \eqref{def:iD_highcost} and \eqref{def:iH_highcost}, respectively, for $(i,k,\ell) \in \tilde{\X}_{diff}$. Consequently, if $n=1$, both an optimal policy and $\pi'$ follow the \textbf{independent threshold type with threshold} $\tilde{T}_D(\ell) = \tilde{i}_D(\ell) + 1$ and $\tilde{T}_H(\ell) = \tilde{i}_H(\ell) + 1$, respectively.
As before, we evaluate several benchmark policies $\tilde{\pi}_j$, defined below, which also follow the \textbf{independent threshold type with thresholds $\tilde{T}_j(\ell)$} for $j = 1,2, \ldots, 4$, for comparison. We only compare the results for $\mu_1 < \mu_2$, as summarized in Tables \ref{table:cmp_state20_highcost} and \ref{table:cmp_state30_highcost}, since otherwise both $\tilde{i}_D(\ell)$ and $\tilde{i}_H(\ell)$ are infinite. 
\begin{enumerate} [leftmargin=3cm, label=Policy $\tilde{\pi}_\arabic*$.]
    \item Always go to Station 1: $a=0$, $\forall n = 1,2$.
    \item Go to Station 1 if and only if the number of jobs waiting in queue exceeds the fixed threshold $10$: $a = \mathbb{I}\{i > 10\}$, $\forall n = 1,2$.
    \item Always go to Station 2: $a=1$, $\forall n = 1,2$.
    \item Go to Station 2 if and only if no waiting is required to commence Station 2 service: $a = \mathbb{I}\{\ell<C_2\}$ if $n=1$; and $a = \mathbb{I}\{\ell \leq C_2\}$ if $n=2$.
\end{enumerate}

\begin{table}[htbp]
    \caption{Comparison of relative errors $(\%)$ for states $(20,k_0,\ell_0)$ when $\frac{h_1}{\mu_1} < \frac{h_2}{\mu_2}$ and $\mu_1 < \mu_2$.}
    \centering
    \label{table:cmp_state20_highcost}
    \begin{tabular}{|c|r|r|r|r|r|r|r|}
        \hline
        \multicolumn{2}{|c|}{} & \multicolumn{1}{c|}{$C_1 = 2$} & \multicolumn{2}{c|}{$C_1 = 3$} & \multicolumn{3}{c|}{$C_1 = 4$} \\ \cline{3-8}
        \multicolumn{2}{|c|}{} & $C_2 = 1$ & $C_2 = 1$ & $C_2 = 2$ & $C_2 = 1$ & $C_2 = 2$ & $C_2 = 3$ \\
        \hline

        \multirow{3}{*}{\textbf{Policy $\pi'$}} & Max error & 0.05 & 0.07 & 0.19 & 0.09 & 0.21 & 0.38 \\
                                               & Avg error & 0.01 & 0.01 & 0.03 & 0.01 & 0.04 & 0.06 \\
                                               & Std error & 0.02 & 0.02 & 0.06 & 0.03 & 0.07 & 0.12 \\
        \hline

        \multirow{3}{*}{\textbf{Policy $\pi_1$}} & Max error & 13.66 & 7.41 & 12.40 & 4.46 & 7.56 & 9.37 \\
                                                & Avg error & 1.49  & 0.53 & 1.10  & 0.20 & 0.48 & 0.69 \\
                                                & Std error & 3.16  & 1.38 & 2.73  & 0.66 & 1.39 & 1.97 \\
        \hline

        \multirow{3}{*}{\textbf{Policy $\pi_2$}} & Max error & 91.73 & 144.59 & 57.04 & 187.62 & 79.32 & 44.26 \\
                                                & Avg error & 59.00 & 102.03 & 29.41 & 135.00 & 48.53 & 19.81 \\
                                                & Std error & 14.45 & 18.95  & 11.50 & 22.31  & 14.51 & 9.99 \\
        \hline

        \multirow{3}{*}{\textbf{Policy $\pi_3$}} & Max error & 190.82 & 308.29 & 113.75 & 409.46 & 164.29 & 88.41 \\
                                                & Avg error & 114.14 & 200.69 & 63.66  & 271.14 & 104.87 & 46.11 \\
                                                & Std error & 34.69  & 47.70  & 23.90  & 58.65  & 29.19  & 19.64 \\
        \hline

        \multirow{3}{*}{\textbf{Policy $\pi_4$}} & Max error & 32.49 & 22.69 & 41.57 & 17.38 & 32.46 & 45.05 \\
                                                & Avg error & 10.08 & 7.30  & 14.72 & 5.36  & 11.73 & 17.32 \\
                                                & Std error & 9.03  & 5.74  & 11.42 & 4.00  & 8.25  & 12.08 \\
        \hline
    \end{tabular}
\end{table}

\begin{table}[htbp]
    \caption{Comparison of relative errors $(\%)$ for states $(30,k_0,\ell_0)$ when $\frac{h_1}{\mu_1} < \frac{h_2}{\mu_2}$ and $\mu_1 < \mu_2$.}
    \centering
    \label{table:cmp_state30_highcost}
    \begin{tabular}{|c|r|r|r|r|r|r|r|}
        \hline
        \multicolumn{2}{|c|}{} & \multicolumn{1}{c|}{$C_1 = 2$} & \multicolumn{2}{c|}{$C_1 = 3$} & \multicolumn{3}{c|}{$C_1 = 4$} \\ \cline{3-8}
        \multicolumn{2}{|c|}{} & $C_2 = 1$ & $C_2 = 1$ & $C_2 = 2$ & $C_2 = 1$ & $C_2 = 2$ & $C_2 = 3$ \\
        \hline

        \multirow{3}{*}{\textbf{Policy $\pi'$}} & Max error & 0.03 & 0.04 & 0.10 & 0.05 & 0.11 & 0.20 \\
                                               & Avg error & 0.01 & 0.01 & 0.02 & 0.01 & 0.02 & 0.04 \\
                                               & Std error & 0.01 & 0.01 & 0.03 & 0.02 & 0.04 & 0.07 \\
        \hline

        \multirow{3}{*}{\textbf{Policy $\pi_1$}} & Max error & 16.65 & 9.72 & 17.17 & 6.34 & 11.33 & 14.89 \\
                                                & Avg error & 2.31  & 0.96 & 1.95  & 0.44 & 1.00  & 1.46 \\
                                                & Std error & 4.17  & 2.07 & 4.11  & 1.13 & 2.38  & 3.43 \\
        \hline

        \multirow{3}{*}{\textbf{Policy $\pi_2$}} & Max error & 131.11 & 215.88 & 81.19 & 290.75 & 119.12 & 64.07 \\
                                                & Avg error & 80.95  & 147.08 & 43.43 & 203.43 & 75.01  & 30.76 \\
                                                & Std error & 20.79  & 29.02  & 15.97 & 35.93  & 19.86  & 14.04 \\
        \hline

        \multirow{3}{*}{\textbf{Policy $\pi_3$}} & Max error & 198.28 & 330.11 & 121.44 & 449.58 & 180.73 & 95.99 \\
                                                & Avg error & 114.80 & 209.99 & 65.31  & 293.28 & 112.02 & 48.43 \\
                                                & Std error & 36.04  & 52.50  & 25.60  & 66.34  & 32.75  & 21.62 \\
        \hline

        \multirow{3}{*}{\textbf{Policy $\pi_4$}} & Max error & 32.09 & 22.65 & 41.87 & 17.46 & 32.82 & 46.05 \\
                                                & Avg error & 8.91  & 6.78  & 13.34 & 5.24  & 11.00 & 16.05 \\
                                                & Std error & 9.04  & 6.01  & 11.86 & 4.28  & 8.82  & 12.90 \\
        \hline
    \end{tabular}
\end{table}
From the definition of the independent threshold-type policy in Definition \ref{def:threshold_policy}, the corresponding thresholds are given by $\tilde{T}_1(\ell) = \infty$, $\tilde{T}_2(\ell) = 11$, $\tilde{T}_3(\ell) = 1$, which remain constant across states. Meanwhile, $\tilde{T}_4(\ell)$ is state-dependent and defined as $\tilde{T}_4(\ell) = 1$ if $\ell < C_2$ and $\tilde{T}_4(\ell) = \infty$ otherwise.

Recall that when $\ell \geq C_2$ under $\mu_1 < \mu_2$, we have $D(i,k,\ell) \leq 0$ for all $i$ by Proposition \ref{prop:larger_mu2_always_action0}. Therefore, if $n=1$, an optimal decision chooses the independent service ($a=0$), as also indicated in Statement \ref{state1.1:larger_mu1_larger_ell_action0} of Theorem \ref{thm:policy2_i}. In this case, Policy $\tilde{\pi}_1$ and Policy $\tilde{\pi}_4$ align with the optimal action, explaining their superior performance among the four benchmarks. In contrast, Policy $\tilde{\pi}_3$ always chooses $a=1$ even when $\ell \geq C_2$, which incurs the highest error across all parameter settings. 
When $\ell < C_2$ (still under $\mu_1 < \mu_2$), an optimal action is $a = \mathbb{I}\{i \geq \tilde{i}_D(\ell)\}$ for some finite threshold $\tilde{i}_D(\ell)$ (see Statement \ref{state:heur_bd2} in Theorem \ref{thm:heur_highcost}). The test states $i_0=20,30$ often lie near or below this threshold, leading to Policy $\tilde{\pi}_1$ outperforming Policy $\tilde{\pi}_4$.
However, its performance deteriorates under high-initial queueing conditions, where all three performance measures worsen --- supporting the theoretical insight that longer initial queues favor collaborative service.

The proposed heuristic policy $\pi'$ achieves costs within 0.5\% of the optimal, even when the worst-case bound identified in Statement \ref{state:heur_bd2} of Theorem \ref{thm:heur_highcost} is attained. For instance, the largest relative error under $C_1=4$ and $C_2=3$ with $\mu_1 < \mu_2$ occurs at state $(i_0, 1, 3)$ when $\mu_1=10$, $\mu_2=15$, $h_0=1$, $h_1=1$, and $h_2=2$, where $\tilde{i}_D(1)=7$ and $\tilde{i}_H(1)=4$.
As observed previously when $\frac{h_1}{\mu_1} > \frac{h_2}{\mu_2}$, a similar pattern holds when $\frac{h_1}{\mu_1} \leq \frac{h_2}{\mu_2}$: states with $i_0 = 30$ exhibit smaller relative errors than those with $i_0 = 20$. This further demonstrates the improved performance of $\pi'$ when the initial queue is high, where the impact of early suboptimal actions is diminished over time.

\section{Conclusion} \label{sec:conclusions}
This paper analyzed a controlled queueing system featuring with both flexible (Type-I) and dedicated (Type-II) servers. Type-I servers also serve as decision-makers, determining for each job whether to engage in independent service or collaborative service. A clearing system MDP model was utilized to formalize this sequential decision-making scenario, with the goal of completing existing jobs or returning to the normative state.

The structural properties of the optimal policy were fully characterized across the parameter space, revealing a threshold-based policy structure, allowing for an infinite threshold in some cases. Based on these insights, we proposed an adaptive heuristic threshold policy tailored to system parameters. Both theoretical and numerical analyses demonstrated the accuracy and robustness of the proposed thresholds in comparison to the optimal ones, particularly when the initial waiting queue is long. 
This study motivates future research in the following direction:

\textbf{Expanding to multi-stage queueing systems.} Our model assumes that decisions are made in real-time upon initial interaction with a new job. However, in practice, decision-making may require additional time or need to occur after an initial service stage, such as screening job conditions in manufacturing. This leads to a two-stage controlled queueing system, where decisions are made after an upstream service to determine whether additional processing should be performed independently or collaboratively. Exploring such multi-stage dynamics could further enhance the applicability of the proposed framework.

\bibliographystyle{unsrtnat}
\bibliography{policy_ref}

\newpage
\appendix 
\section{Appendix}\label{sec:appendix}

\subsection{Preliminaries for supporting results in Subsection \ref{sec:support-results}}\label{sec:proof-prelims}
To establish the propositions and corollaries in Subsection \ref{sec:support-results}, we first introduce several preliminary lemmas along with their proofs.
Notably, Lemma \ref{lemma:diff}, outlined in Subsection \ref{sec:prelim}, is also essential for the supporting results and is therefore proved here.

We begin by presenting Lemma \ref{lemma:diff} and Lemma \ref{lemma:mono_i} as the remaining preliminaries.
Lemma \ref{silemma:diff} further examines the properties of $D(i,k,\ell)$: Statement \ref{state:diff} reformulates $D(i,k,\ell)$ for $i \geq 1$ using Bellman equations, while Statement \ref{state:diff_lb} provides an affine lower bound in terms of $b'$ and $c'$.
This is followed by Lemma \ref{lemma:mono_i}, which establishes that the value function is non-decreasing $i$.

\begin{lemma}\label{silemma:diff}
    Consider any $x =(i,k,\ell)\in \X_{diff}$.
    \begin{enumerate}
    \item \label{state:diff}
    If $i \geq 1$ (so $x = (i,k,\ell) \in \tilde{\X}_{diff}$),
    \begin{align}
        \lefteqn{D(i,k,\ell) = \frac{-i\big(\mu_1-\mathbb{I}\{\ell<C_2\}\mu_2\big)h_0}{d(k,\ell)d(k-1,\ell+1)} + \frac{k h_1 + \ell h_2 }{d(k,\ell)} - \frac{(k-1)h_1+(\ell+1)h_2}{d(k-1,\ell+1)}}& \nonumber \\
        &\qquad +\frac{(k-1) \mu_1}{d(k-1,\ell+1)} \Big[\min\{v(i-1, k, \ell), v(i-1, k-1, \ell+1)\}\nonumber \\
        &\qquad \quad - \min\{v(i-1, k-1, \ell+1), v(i-1, k-2, \ell+2)\} \Big] \nonumber\\ 
        &\qquad + \frac{\min\{\ell,C_2\}\mu_2}{d(k,\ell)} \Big[\min\{v(i-1, k+1, \ell-1), v(i-1, k, \ell)\}\nonumber \\
        &\qquad \quad - \min\{v(i-1, k, \ell), v(i-1, k-1, \ell+1)\}\Big]. 
        \label{eq:diff}
    \end{align}
    In addition, the sum of the first three one-step cost terms on the first line of \eqref{eq:diff} equals
    \begin{align}
        \lefteqn{\frac{-i\big(\mu_1-\mathbb{I}\{\ell<C_2\}\mu_2\big)h_0}{d(k,\ell)d(k-1,\ell+1)} + \frac{k h_1 + \ell h_2 }{d(k,\ell)} - \frac{(k-1)h_1+(\ell+1)h_2}{d(k-1,\ell+1)}}& \nonumber \\
       &\quad = \left\{
       \begin{array}{lcl}
            \frac{C_1\mu_1\mu_2}{d(k,\ell)d(k-1,\ell+1)} \left(\frac{ih_0}{C_1}\left(\frac{1}{\mu_1} - \frac{1}{\mu_2}\right) + \frac{h_1}{\mu_1} - \frac{h_2}{\mu_2} \right) & \text{if} & \ell < C_2, \\
            \frac{C_2\mu_1\mu_2}{d(k,\ell)d(k-1,\ell+1)}\left(\frac{-ih_0}{C_2 \mu_2} + \frac{h_1-h_2}{\mu_1} - \frac{C_1h_2}{C_2 \mu_2} \right) & \text{if} & \ell \geq C_2. \\
        \end{array}
        \right. 
        \label{eq:diff_const_expand}
    \end{align}
    With the notations from Definitions \ref{def:probs_k} and \ref{def:probs_quantities_k}, the equation given in \eqref{eq:diff}, can be rearranged as follows:
     \begin{align}
        \lefteqn{D(i,k,\ell)= p_k(ic_k+b_k) }& \nonumber \\
        &\qquad + q_k \Big[\min\{v(i-1, k, \ell), v(i-1, k-1, \ell+1)\} \nonumber \\
        &\qquad \quad - \min\{v(i-1, k-1, \ell+1), v(i-1, k-2, \ell+2)\} \Big] \nonumber\\
        &\qquad + r_k \Big[\min\{v(i-1, k+1, \ell-1), v(i-1, k, \ell)\} \nonumber \\
        &\qquad \quad - \min\{v(i-1, k, \ell), v(i-1, k-1, \ell+1)\}\Big]. \label{eq:alge_diff}
    \end{align}
    \item \label{state:diff_lb}
    The difference $D(i,k,\ell)$ is also globally lower bounded by an affine function in $i$ as follows.
    \begin{align}
        D(i,k,\ell) \geq ic'+b',
        \label{eq:diff_lb}
    \end{align}
    \end{enumerate}
\end{lemma}

\begin{lemma}\label{lemma:mono_i}
    Consider any state $x = (i,k,\ell) \in \X$. The following holds for all $i\geq 0$,
    \begin{align}
        v(i+1,k,\ell) - v(i,k,\ell) \geq 0. \label{eq:mono_i}
    \end{align}
\end{lemma}

\noindent First, we prove Statement \ref{state:diff_bdry} in Lemma \ref{lemma:diff}, which provides the values $D(i,k,\ell)$ at $i=0$.

\begin{proof}[Proof of Statement \ref{state:diff_bdry} in Lemma \ref{lemma:diff}.]
    We can prove the result at the boundary by induction on the number of jobs in the system $M:=k+\ell$ for all $x\in \X$ where $i = 0$, and the details are omitted here for brevity.
\end{proof}

\noindent Next, we proceed with the proofs of both statements in Lemma \ref{silemma:diff}.

\begin{proof} [Proof of Statement \ref{state:diff} in Lemma \ref{silemma:diff}]
    For any $x = (i,k,\ell) \in \tilde{\X}_{diff}$ where $i \geq 1$, substituting Bellman equations for $v(i,k,\ell)$ and $v(i,k-1,\ell+1)$ yields
\begin{align*}
    \lefteqn{D(i,k,\ell) = v(i,k,\ell) - v(i,k-1,\ell+1)}& \nonumber \\
    &\quad = \frac{i h_0+ k h_1 + \ell h_2 }{d(k,\ell)} - \frac{ih_0+(k-1)h_1+(\ell+1)h_2}{d(k-1,\ell+1)}\nonumber \\
    &\qquad + \frac{k \mu_1}{d(k,\ell)} \min\{v(i-1, k, \ell), v(i-1, k-1, \ell+1)\} \nonumber \\
    & \qquad + \frac{\min\{\ell,C_2\}\mu_2}{d(k,\ell)} \min\{v(i-1, k+1, \ell-1), v(i-1, k, \ell)\}\nonumber \\
    &\qquad - \frac{(k-1) \mu_1}{d(k-1,\ell+1)} \min\{v(i-1, k-1, \ell+1), v(i-1, k-2, \ell+2)\} \nonumber \\
    & \qquad - \frac{\min\{\ell+1,C_2\}\mu_2}{d(k-1,\ell+1)} \min\{v(i-1, k, \ell), v(i-1, k-1, \ell+1)\}\nonumber\\
    &\quad = \frac{-i\big(\mu_1-\mathbb{I}\{\ell<C_2\}\mu_2\big)h_0}{d(k,\ell)d(k-1,\ell+1)} + \frac{k h_1 + \ell h_2 }{d(k,\ell)} - \frac{(k-1)h_1+(\ell+1)h_2}{d(k-1,\ell+1)} \nonumber \\
    &\qquad +\frac{(k-1) \mu_1}{d(k-1,\ell+1)} \Big[\min\{v(i-1, k, \ell), v(i-1, k-1, \ell+1)\}\nonumber \\
    &\qquad \quad - \min\{v(i-1, k-1, \ell+1), v(i-1, k-2, \ell+2)\} \Big] \nonumber\\ 
    &\qquad + \frac{\min\{\ell,C_2\}\mu_2}{d(k,\ell)} \Big[\min\{v(i-1, k+1, \ell-1), v(i-1, k, \ell)\}\nonumber \\
    &\qquad \quad - \min\{v(i-1, k, \ell), v(i-1, k-1, \ell+1)\}\Big]\\
    &\qquad + \left(\left(\frac{k \mu_1}{d(k,\ell)} - \frac{(k-1) \mu_1}{d(k-1,\ell+1)}\right) - \left(\frac{\min\{\ell+1,C_2\}\mu_2}{d(k-1,\ell+1)} - \frac{\min\{\ell,C_2\}\mu_2}{d(k,\ell)}\right)\right)\cdot \\
    &\qquad \quad\min\{v(i-1, k, \ell), v(i-1, k-1, \ell+1)\}.
\end{align*}
In the last equality, we group the terms that determine where to work after a Station 1 service completion at states $(i,k,\ell)$ and $(i,k-1,\ell+1)$ with a smaller coefficient $\frac{(k-1) \mu_1}{d(k-1,\ell+1)}$. A similar grouping applies to the terms representing where to work after a Station 2 service completion. The result follows from \eqref{eq:diff} by applying the following arithmetic computations:
\begin{align*}
    \left(\frac{k \mu_1}{d(k,\ell)} - \frac{(k-1) \mu_1}{d(k-1,\ell+1)}\right) - \left(\frac{\min\{\ell+1,C_2\}\mu_2}{d(k-1,\ell+1)} - \frac{\min\{\ell,C_2\}\mu_2}{d(k,\ell)}\right) = 0.
\end{align*}
Recalling $k+\ell=C_1$, some arithmetic yields \eqref{eq:diff_const_expand}. 
\end{proof}

\begin{proof}[Proof of Statement \ref{state:diff_lb} in Lemma \ref{silemma:diff}]
    Recall from the definition of $c'$ in \eqref{def:b'_and_c'} that $c' < 0$ holds.
    We proceed by induction on $i$. For the base case $i=0$, From Equation \eqref{eq:diff_bdry} (Statement \ref{state:diff_bdry} in Lemma \ref{lemma:diff}), we obtain:
    \begin{align*}
        D(0,k,\ell) = \frac{h_1}{\mu_1} - \frac{\max\{\ell+1,C_2\}}{C_2}\frac{h_2}{\mu_2} \geq \frac{h_1}{\mu_1} - \frac{C_1h_2}{C_2\mu_2} \geq \frac{h_1-h_2}{\mu_1} - \frac{C_1h_2}{C_2\mu_2} = b',
    \end{align*}where the last step follows from the definition of $b'$ in \eqref{def:b'_and_c'}.
    Therefore, the base case for $i=0$ is verified. Now, assume \eqref{eq:diff_lb} holds for $i-1$ and consider the difference in \eqref{eq:alge_diff} at $i \geq 1$. 
    Notice that the difference with coefficient $q_k$ in the second term in \eqref{eq:alge_diff} is in fact lower bounded by $ic'+b'$:
    \begin{align}
        \lefteqn{\min\{v(i-1, k, \ell), v(i-1, k-1, \ell+1)\}}& \nonumber \\
        &\qquad  - \min\{v(i-1, k-1, \ell+1), v(i-1, k-2, \ell+2)\} \nonumber \\
        & \quad = \min\Big\{v(i-1, k, \ell) - \min\{v(i-1, k-1, \ell+1), v(i-1, k-2, \ell+2)\}, \nonumber \\
        &\qquad \quad v(i-1, k-1, \ell+1) - \min\{v(i-1, k-1, \ell+1), v(i-1, k-2, \ell+2)\}\Big\} \nonumber \\
        &\quad \geq \min\{v(i-1, k, \ell) - v(i-1, k-1, \ell+1), v(i-1, k-1, \ell+1) - v(i-1, k-2, \ell+2)\} \nonumber \\
        &\quad \geq \min\{(i-1)c'+b', (i-1)c'+b'\} \nonumber \\
        &\quad \geq ic'+b',\label{eq:diff_lb_mu1_terms}
    \end{align}where the second-to-last step follows by using the inductive hypothesis (twice), and the last step holds since $c'<0$.
    Similarly, by a comparable argument, the difference in the term with coefficient $r_k$ in \eqref{eq:alge_diff} is also lower bounded by $ic'+b'$:
    \begin{align}
        \lefteqn{\min\{v(i-1, k+1, \ell-1), v(i-1, k, \ell)\}}& \nonumber \\
        &\qquad - \min\{v(i-1, k, \ell), v(i-1, k-1, \ell+1)\} \nonumber \\
        & \quad \geq ic'+b'.\label{eq:diff_lb_mu2_terms}
    \end{align}
    Applying both \eqref{eq:diff_lb_mu1_terms} and \eqref{eq:diff_lb_mu2_terms} in \eqref{eq:alge_diff} yields
    \begin{align*}
        \lefteqn{D(i,k,\ell) \geq p_k(ic_k+b_k) + (q_k+r_k) (ic'+b')}&\\
        &\quad \geq (p_k+q_k+r_k) (ic'+b')\\
        &\quad  = ic'+b',
    \end{align*}where the second inequality holds by \eqref{eq:cmp_bk_b_ck_c}, which states that $c_k \geq c'$ and $b_k \geq b'$, and the last step follows from the identity in \eqref{eq:identity}.
\end{proof}

\noindent Now, we return to proving the remaining statement in Lemma \ref{lemma:diff}.
\begin{proof}[Proof of Statement \ref{state:diff_ub_larger_mu2} in Lemma \ref{lemma:diff}]
    Referring to the expression of $c$ in Definition \ref{def:b_and_c}, it follows that $c = \frac{h_0}{C_1}\left(\frac{1}{\mu_1} - \frac{1}{\mu_2}\right) \geq 0$ since $\mu_2 \geq \mu_1$.
    The proof proceeds with induction on $i$. When $i = 0$, the boundary values in \eqref{eq:diff_bdry} imply that
    \begin{align*}
        D(0,k,\ell) = \frac{h_1}{\mu_1} - \frac{\max\{\ell+1,C_2\}}{C_2}\frac{h_2}{\mu_2} \leq \frac{h_1}{\mu_1} - \frac{h_2}{\mu_2} = b,
    \end{align*}where the last step relies on the definition of $b$ in \eqref{def:b_and_c}. This verifies the base case at $i=0$.
    Now, we prove the upper bound \eqref{eq:diff_ub_larger_mu2} at a generic $i$, where $i \geq 1$, by assuming it holds at $i-1$. 
    Consider the alternative expression of $D(i,k,\ell)$ given in \eqref{eq:alge_diff}.
    Notice the difference with coefficient $q_k$ in the second term is upper bounded by $ic+b$:
    \begin{align}
        \lefteqn{\min\{v(i-1, k, \ell), v(i-1, k-1, \ell+1)\}}& \nonumber \\
        &\qquad  - \min\{v(i-1, k-1, \ell+1), v(i-1, k-2, \ell+2)\} \nonumber \\
        & \quad = \max\Big\{\min\{v(i-1, k, \ell), v(i-1, k-1, \ell+1)\} - v(i-1, k-1, \ell+1), \nonumber \\
        &\qquad \quad \min\{v(i-1, k, \ell), v(i-1, k-1, \ell+1)\} - v(i-1, k-2, \ell+2)\Big\} \nonumber \\
        &\quad \leq \max\{v(i-1, k, \ell) - v(i-1, k-1, \ell+1), v(i-1, k-1, \ell+1) - v(i-1, k-2, \ell+2)\} \nonumber \\
        &\quad \leq \max\{(i-1)c+b, (i-1)c+b\} \nonumber \\
        &\quad \leq ic+b,\label{eq:diff_ub_mu1_terms}
    \end{align}where the second-to-last step utilizes the inductive hypothesis twice, and the last step follows from $c \geq 0$.
    Following a similar argument, the difference with coefficient $r_k$ in the third term is also upper bounded by $ic+b$:
   \begin{align}
        \lefteqn{\min\{v(i-1, k+1, \ell-1), v(i-1, k, \ell)\}}& \nonumber \\
        &\qquad \quad - \min\{v(i-1, k, \ell), v(i-1, k-1, \ell+1)\} \nonumber \\
        &\quad \leq ic+b. \label{eq:diff_ub_mu2_terms}
    \end{align}
    Applying both \eqref{eq:diff_ub_mu1_terms} and \eqref{eq:diff_ub_mu2_terms} in \eqref{eq:alge_diff} yields
    \begin{align*}
        \lefteqn{D(i,k,\ell) \leq p_k(ic_k+b_k) + (q_k+r_k) (ic+b)}&\\
        &\quad \leq (p_k+q_k+r_k) (ic+b)\\
        &\quad  = ic+b,
    \end{align*}where the second inequality holds by \eqref{eq:cmp_bk_b_ck_c}, which states that $c_k \leq c$ and $b_k \leq b$, and the last step follows from the identity in \eqref{eq:identity}.
\end{proof}

\noindent We conclude this subsection with the proof of Lemma \ref{lemma:mono_i}. 

\begin{proof} [Proof of Lemma \ref{lemma:mono_i}]
    The proof proceeds by induction on $i$, where the base case when $i=0$ utilizes the result of the boundary values in \eqref{eq:diff_bdry}. As the details are analogous to many others in this paper, they are omitted here.
\end{proof}

\subsection{Proof of supporting results in Subsection \ref{sec:support-results}} \label{sec:}
This subsection presents the proofs for the supporting results, as outlined in Subsection \ref{sec:support-results}.
\begin{proof} [Proof of Proposition \ref{prop:larger_mu1_dec}.]
We prove inequality \eqref{eq:larger_mu1_dec} by induction on $i$.
When $i=0$, applying the boundary values in \eqref{eq:diff_bdry} and the equation in \eqref{eq:diff} at $i=1$ yields
\begin{align}
    \lefteqn{v(0,k,\ell) - v(0,k-1,\ell+1) - v(1,k,\ell) + v(1,k-1,\ell+1)}& \nonumber\\
    & \quad = \frac{h_1}{\mu_1} - \frac{\max\{\ell+1,C_2\}}{C_2}\frac{h_2}{\mu_2} \nonumber\\
    & \qquad - \left(\frac{-\big(\mu_1-\mathbb{I}\{\ell<C_2\}\mu_2\big)h_0}{d(k,\ell)d(k-1,\ell+1)} + \frac{k h_1 + \ell h_2 }{d(k,\ell)} - \frac{(k-1)h_1+(\ell+1)h_2}{d(k-1,\ell+1)}\right) \nonumber\\
    &\qquad -\frac{(k-1) \mu_1}{d(k-1,\ell+1)} \Big[\min\{v(0, k, \ell), v(0, k-1, \ell+1)\} \nonumber\\
    &\qquad \qquad - \min\{v(0, k-1, \ell+1), v(0, k-2, \ell+2)\} \Big] \nonumber\\
    &\qquad - \frac{\min\{\ell,C_2\}\mu_2}{d(k,\ell)} \Big[\min\{v(0, k+1, \ell-1), v(0, k, \ell)\} - \min\{v(0, k, \ell), v(0, k-1, \ell+1)\}\Big] \nonumber\\
    &\quad \geq \frac{h_1}{\mu_1} - \frac{\max\{\ell+1,C_2\}}{C_2}\frac{h_2}{\mu_2} - \frac{k h_1 + \ell h_2 }{d(k,\ell)} + \frac{(k-1)h_1+(\ell+1)h_2}{d(k-1,\ell+1)} \nonumber\\
    &\qquad -\frac{(k-1) \mu_1}{d(k-1,\ell+1)} \Big[\min\{v(0, k, \ell), v(0, k-1, \ell+1)\} \nonumber\\
    &\qquad \qquad - \min\{v(0, k-1, \ell+1), v(0, k-2, \ell+2)\} \Big]  \label{eq:larger_mu1_dec_inte0}\\
    &\qquad - \frac{\min\{\ell,C_2\}\mu_2}{d(k,\ell)} \Big[\min\{v(0, k+1, \ell-1), v(0, k, \ell)\} - \min\{v(0, k, \ell), v(0, k-1, \ell+1)\}\Big], \nonumber
\end{align}
where the inequality holds since $\mu_1-\mathbb{I}\{\ell<C_2\}\mu_2 \geq \mu_1-\mu_2 \geq 0$. Consider the difference of the two minima 
in the term with coefficient $-\frac{(k-1) \mu_1}{d(k-1,\ell+1)}$ in \eqref{eq:larger_mu1_dec_inte0}. Notice that
\begin{align}
    \lefteqn{\min\{v(0, k, \ell), v(0, k-1, \ell+1)\} - \min\{v(0, k-1, \ell+1), v(0, k-2, \ell+2)\}}& \nonumber \\
    &\quad = \max\Big\{\min\{v(0, k, \ell), v(0, k-1, \ell+1)\} - v(0, k-1, \ell+1),  \nonumber \\
    &\qquad \qquad \min\{v(0, k, \ell), v(0, k-1, \ell+1)\} - v(0, k-2, \ell+2)\Big\} \nonumber \\
    &\quad \leq \max\Big\{v(0, k, \ell) - v(0, k-1, \ell+1), v(0, k-1, \ell+1) - v(0, k-2, \ell+2)\Big\} \nonumber \\
    &\quad = \max\Big\{\frac{h_1}{\mu_1} - \frac{\max\{\ell+1,C_2\}}{C_2}\frac{h_2}{\mu_2}, \frac{h_1}{\mu_1} - \frac{\max\{\ell+2,C_2\}}{C_2}\frac{h_2}{\mu_2}\Big\} \nonumber \\
    &\quad =  \frac{h_1}{\mu_1} - \frac{\max\{\ell+1,C_2\}}{C_2}\frac{h_2}{\mu_2} \label{eq:larger_mu1_dec_inte4}.
\end{align}
Similarly, the difference $\min\{v(0, k+1, \ell-1), v(0, k, \ell)\} - \min\{v(0, k, \ell), v(0, k-1, \ell+1)\}\leq\frac{h_1}{\mu_1} - \frac{\max\{\ell,C_2\}}{C_2}\frac{h_2}{\mu_2}$ in the term with coefficient $- \frac{\min\{\ell,C_2\}\mu_2}{d(k,\ell)}$ in \eqref{eq:larger_mu1_dec_inte0}.
Applying the two bounds back to \eqref{eq:larger_mu1_dec_inte0} yields
\begin{align*}
    \lefteqn{v(0,k,\ell) - v(0,k-1,\ell+1) - v(1,k,\ell) + v(1,k-1,\ell+1)}&\\
    &\quad \geq \frac{h_1}{\mu_1} - \frac{\max\{\ell+1,C_2\}}{C_2}\frac{h_2}{\mu_2} - \left(\frac{k h_1 + \ell h_2 }{d(k,\ell)} - \frac{(k-1)h_1+(\ell+1)h_2}{d(k-1,\ell+1)}\right)\\
    &\qquad -\frac{(k-1) \mu_1}{d(k-1,\ell+1)} \left(\frac{h_1}{\mu_1} - \frac{\max\{\ell+1,C_2\}}{C_2}\frac{h_2}{\mu_2}\right) - \frac{\min\{\ell,C_2\}\mu_2}{d(k,\ell)} \left(\frac{h_1}{\mu_1} - \frac{\max\{\ell,C_2\}}{C_2}\frac{h_2}{\mu_2}\right)\\
    &\quad = \left(1-\frac{k\mu_1}{d(k,\ell)} - \frac{\min\{\ell,C_2\}\mu_2}{d(k,\ell)}\right)\frac{h_1}{\mu_1} - \left(\frac{\ell\mu_2}{d(k,\ell)} - \frac{\min\{\ell,C_2\}\mu_2}{d(k,\ell)}\right)\frac{h_2}{\mu_2}\\
    &\qquad - \left(\frac{\max\{\ell+1,C_2\}}{C_2}\left(1-\frac{(k-1) \mu_1}{d(k-1,\ell+1)}\right)-\frac{(\ell+1)\mu_2}{d(k-1,\ell+1)}\right)\frac{h_2}{\mu_2} \\
    &\quad = 0,
\end{align*}
where the last equality follows by noticing that $\min\{\ell,C_2\}\cdot\max\{\ell,C_2\} = \ell C_2$ and similarly $\min\{\ell+1,C_2\}\cdot\max\{\ell+1,C_2\} = (\ell+1) C_2$. Hence the base case when $i=0$ holds in \eqref{eq:larger_mu1_dec}. Assume the result holds at $i-1$, consider it at $i \geq 1$. Applying \eqref{eq:diff} at $i$ and $i+1$ yields
    \begin{align}
    \lefteqn{v(i,k,\ell) - v(i,k-1,\ell+1) - v(i+1,k,\ell) + v(i+1,k-1,\ell+1)}& \nonumber \\
    & \quad = \frac{\big(\mu_1-\mathbb{I}\{\ell<C_2\}\mu_2\big)h_0}{d(k,\ell)d(k-1,\ell+1)} \nonumber \\ 
    &\qquad +\frac{(k-1) \mu_1}{d(k-1,\ell+1)} \Big[\min\{v(i-1, k, \ell), v(i-1, k-1, \ell+1)\} \nonumber \\
    &\qquad \quad - \min\{v(i-1, k-1, \ell+1), v(i-1, k-2, \ell+2)\} \nonumber \\
    &\qquad \quad - \min\{v(i, k, \ell), v(i, k-1, \ell+1)\} \nonumber \\
    &\qquad \quad + \min\{v(i, k-1, \ell+1), v(i, k-2, \ell+2)\}\Big] \nonumber \\
    &\qquad + \frac{\min\{\ell,C_2\}\mu_2}{d(k,\ell)} \Big[\min\{v(i-1, k+1, \ell-1), v(i-1, k, \ell)\} \nonumber \\
    &\qquad \quad - \min\{v(i-1, k, \ell), v(i-1, k-1, \ell+1)\} \nonumber \\
    &\qquad \quad - \min\{v(i, k+1, \ell-1), v(i, k, \ell)\} \nonumber \\
    &\qquad \quad + \min\{v(i, k, \ell), v(i, k-1, \ell+1)\}\Big] 
    \label{eq:larger_mu1_dec_inte1}
\end{align}
Consider the term with coefficient $\frac{(k-1) \mu_1}{d(k-1,\ell+1)}$ in \eqref{eq:larger_mu1_dec_inte1} 
first. There are two cases to examine based on whether $v(i-1, k, \ell) > v(i-1, k-1, \ell+1)$.
\begin{enumerate}[label= \textbf{Case} \arabic*:, leftmargin=3\parindent]
    \item If $v(i-1, k, \ell) > v(i-1, k-1, \ell+1)$, we choose an upper bound $v(i, k-1, \ell+1)$ in the third minimum and consider the two sub-cases
    \begin{enumerate}[label= \textbf{Subcase} (\alph*):, leftmargin=3\parindent]
        \item If $v(i, k-1, \ell+1) \leq v(i, k-2, \ell+2)$, replacing the second minimum with an upper bound $v(i-1, k-1, \ell+1)$ yields
        \begin{align*}
            \lefteqn{\min\{v(i-1, k, \ell), v(i-1, k-1, \ell+1)\}}\\
            &\qquad - \min\{v(i-1, k-1, \ell+1), v(i-1, k-2, \ell+2)\}\\
            &\qquad  - \min\{v(i, k, \ell), v(i, k-1, \ell+1)\}\\
            &\qquad + \min\{v(i, k-1, \ell+1), v(i, k-2, \ell+2)\}\\
            &\quad \geq v(i-1, k-1, \ell+1) - v(i-1, k-1, \ell+1)\\
            &\qquad - v(i, k-1, \ell+1) + v(i, k-1, \ell+1)\\
            &\quad = 0.
        \end{align*}
        \item If $v(i, k-1, \ell+1) > v(i, k-2, \ell+2)$, replacing the second minimum with an upper bound $v(i-1, k-2, \ell+2)$ yields
        \begin{align}
            \lefteqn{\min\{v(i-1, k, \ell), v(i-1, k-1, \ell+1)\}} \nonumber \\
            &\qquad - \min\{v(i-1, k-1, \ell+1), v(i-1, k-2, \ell+2)\} \nonumber \\
            &\qquad - \min\{v(i, k, \ell), v(i, k-1, \ell+1)\} \nonumber \\
            &\qquad + \min\{v(i, k-1, \ell+1), v(i, k-2, \ell+2)\} \nonumber \\
            &\quad \geq v(i-1, k-1, \ell+1) - v(i-1, k-2, \ell+2) \nonumber \\
            &\qquad - v(i, k-1, \ell+1) + v(i, k-2, \ell+2) \nonumber \\
            &\quad \geq 0, \label{eq:larger_mu1_dec_inte2}
        \end{align}by the inductive hypothesis.
    \end{enumerate}
    \item If $v(i-1, k, \ell) \leq v(i-1, k-1, \ell+1)$, we choose an upper bound $v(i, k, \ell)$ in the third minimum and consider the two subcases
    \begin{enumerate}
        \item If $v(i, k-1, \ell+1) \leq v(i, k-2, \ell+2)$, replacing the second minimum with an upper bound $v(i-1, k-1, \ell+1)$ yields
        \begin{align}
            \lefteqn{\min\{v(i-1, k, \ell), v(i-1, k-1, \ell+1)\}} \nonumber \\
            &\qquad - \min\{v(i-1, k-1, \ell+1), v(i-1, k-2, \ell+2)\} \nonumber \\
            &\qquad - \min\{v(i, k, \ell), v(i, k-1, \ell+1)\} \nonumber \\
            &\qquad  + \min\{v(i, k-1, \ell+1), v(i, k-2, \ell+2)\} \nonumber \\
            &\quad \geq v(i-1, k, \ell) - v(i-1, k-1, \ell+1) - v(i, k, \ell) + v(i, k-1, \ell+1) \nonumber \\
            &\quad \geq 0, \label{eq:larger_mu1_dec_inte3}
        \end{align}
        where the last inequality holds by the inductive hypothesis.
        \item If $v(i, k-1, \ell+1) > v(i, k-2, \ell+2)$, replacing the second minimum with an upper bound $v(i-1, k-2, \ell+2)$ yields
        \begin{align*}
            \lefteqn{\min\{v(i-1, k, \ell), v(i-1, k-1, \ell+1)\}}\\
            &\qquad - \min\{v(i-1, k-1, \ell+1), v(i-1, k-2, \ell+2)\}\\
            &\qquad - \min\{v(i, k, \ell), v(i, k-1, \ell+1)\}\\
            &\qquad + \min\{v(i, k-1, \ell+1), v(i, k-2, \ell+2)\}\\
            &\quad \geq v(i-1, k, \ell) - v(i-1, k-2, \ell+2) - v(i, k, \ell) + v(i, k-2, \ell+2)\\
            &\quad = \big(v(i-1, k, \ell) - v(i-1, k-1, \ell+1) - v(i, k, \ell) + v(i, k-1, \ell+1)\big)\\
            &\qquad + \big(v(i-1, k-1, \ell+1) - v(i-1, k-2, \ell+2) \nonumber \\
            &\qquad - v(i, k-1, \ell+1) + v(i, k-2, \ell+2)\big) \nonumber\\
            &\quad \geq 0,
        \end{align*}
        where the last inequality follows by adding the two inequalities \eqref{eq:larger_mu1_dec_inte2} and \eqref{eq:larger_mu1_dec_inte3}.
    \end{enumerate}
\end{enumerate}
It follows that the term in \eqref{eq:larger_mu1_dec_inte1} with coefficient $\frac{(k-1) \mu_1}{d(k-1,\ell+1)}$ is non-negative. An analogous argument applies to show the non-negativity of the term in \eqref{eq:larger_mu1_dec_inte1} with coefficient $ \frac{\min\{\ell,C_2\}\mu_2}{d(k,\ell)}$.
Noticing $\mu_1-\mathbb{I}\{\ell<C_2\}\mu_2 \geq \mu_1-\mu_2 \geq 0$ yields the non-negativity of the first term in \eqref{eq:larger_mu1_dec_inte1}, and hence the desired result \eqref{eq:larger_mu1_dec} holds.

It remains to prove the second statement. Given that the difference $v(i,k,\ell) - v(i,k-1,\ell+1)$ is monotone non-increasing in 
$i$, it suffices to show it becomes negative for all $i > N'_1(\ell)$ for some finite $N'_1(\ell)$. By the monotonicity result in \eqref{eq:larger_mu1_dec}, the following holds for all $i \geq 0$ that
\begin{align*}
     \lefteqn{v(i,k,\ell) - v(i,k-1,\ell+1)}&\\
     & \quad \leq v(0,k,\ell) - v(0,k-1,\ell+1)\\
     & \quad = \frac{h_1}{\mu_1} - \frac{\max\{\ell+1,C_2\}}{C_2}\frac{h_2}{\mu_2}.
\end{align*}
Since we only care about the sign of $D(i,k,\ell)$ when $i$ is large, without loss of generality assume $i \geq 1$. Consider the expression in \eqref{eq:diff}. Using a comparable argument as in \eqref{eq:larger_mu1_dec_inte4}, it holds that the difference in the term with coefficient $\frac{(k-1) \mu_1}{d(k-1,\ell+1)}$
\begin{align*}
    \lefteqn{\min\{v(i-1, k, \ell), v(i-1, k-1, \ell+1)\}}&\\
    &\qquad - \min\{v(i-1, k-1, \ell+1), v(i-1, k-2, \ell+2)\}\\
    &\quad \leq \max\Big\{v(i-1, k, \ell) - v(i-1, k-1, \ell+1), \\
    &\qquad \quad v(i-1, k-1, \ell+1) - v(i-1, k-2, \ell+2)\Big\} \\
    &\quad \leq \frac{h_1}{\mu_1} - \frac{\max\{\ell+1,C_2\}}{C_2}\frac{h_2}{\mu_2}.
\end{align*}
 Similarly, the difference in the term in \eqref{eq:diff} with coefficient $ \frac{\min\{\ell,C_2\}\mu_2}{d(k,\ell)}$is upper bounded by
\begin{align*}
    &\min\{v(i-1, k+1, \ell-1), v(i-1, k, \ell)\}- \min\{v(i-1, k, \ell), v(i-1, k-1, \ell+1)\}\\
    &\quad \leq \frac{h_1}{\mu_1} - \frac{\max\{\ell,C_2\}}{C_2}\frac{h_2}{\mu_2}.
\end{align*}
Applying the two upper bounds in the expression in \eqref{eq:diff} yields
\begin{align*}
    \lefteqn{v(i,k,\ell) - v(i,k-1,\ell+1)}&\\
    &\quad \leq \frac{-i\big(\mu_1-\mathbb{I}\{\ell<C_2\}\mu_2\big)h_0}{d(k,\ell)d(k-1,\ell+1)} + \frac{k h_1 + \ell h_2 }{d(k,\ell)} - \frac{(k-1)h_1+(\ell+1)h_2}{d(k-1,\ell+1)}\nonumber \\
    &\qquad +\frac{(k-1) \mu_1}{d(k-1,\ell+1)} \left(\frac{h_1}{\mu_1} - \frac{\max\{\ell+1,C_2\}}{C_2}\frac{h_2}{\mu_2} \right)\nonumber \\
    &\qquad + \frac{\min\{\ell,C_2\}\mu_2}{d(k,\ell)} \left(\frac{h_1}{\mu_1} - \frac{\max\{\ell,C_2\}}{C_2}\frac{h_2}{\mu_2}\right)\\
    &\quad = \frac{-i\big(\mu_1-\mathbb{I}\{\ell<C_2\}\mu_2\big)h_0}{d(k,\ell)d(k-1,\ell+1)} + \left(\frac{h_1}{\mu_1} - \frac{\max\{\ell+1,C_2\}}{C_2}\frac{h_2}{\mu_2} \right)\\
    &\quad = -c_1(\ell)i+c_2(\ell),
\end{align*}where $c_1(\ell) = \frac{\big(\mu_1-\mathbb{I}\{\ell<C_2\}\mu_2\big)h_0}{d(k,\ell)d(k-1,\ell+1)}$ and $c_2(\ell) = \frac{h_1}{\mu_1} - \frac{\max\{\ell+1,C_2\}}{C_2}\frac{h_2}{\mu_2}$, both are functions of $\ell$ (recall $k+\ell = C_1$) and independent of $i$.

Consider the case where $\ell \geq C_2$ so that $c_1(\ell) = \frac{\mu_1h_0}{d(k,\ell)d(k-1,\ell+1)} > 0$ and $c_2(\ell) = \frac{h_1}{\mu_1} - \frac{\ell+1}{C_2}\frac{h_2}{\mu_2}$ the result follows by letting $N'_1(\ell) = \frac{c_2(\ell)}{c_1(\ell)}$. If further assume $\mu_2 > \mu_1$, the result holds by letting $N'_1(\ell) = \frac{c_2(\ell)}{c_1(\ell)}$ again where $c_1(\ell) > 0$ since $\mu_1-\mathbb{I}\{\ell<C_2\}\mu_2 \geq \mu_1-\mu_2 > 0$. 
\end{proof}

\begin{proof} [Proof of Proposition \ref{prop:larger_mu2_always_action1}.]
    Proof by induction on $i$. When $i=0$, by recalling \eqref{eq:diff_bdry}, inequality \eqref{eq:larger_mu2_always_action1} becomes
    \begin{align*}
        v(0,k,\ell) - v(0,k-1,\ell+1) = \frac{h_1}{\mu_1} - \frac{h_2}{\mu_2} \geq 0.
    \end{align*}
    Assume the result holds at $i-1$ to prove it at $i \geq 1$. The expansion of the left-hand side of \eqref{eq:larger_mu2_always_action1} is given in \eqref{eq:diff} by Statement \ref{state:diff} of Lemma \ref{silemma:diff}. 
Consider the difference $\min\{v(i-1, k, \ell), v(i-1, k-1, \ell+1)\} - \min\{v(i-1, k-1, \ell+1), v(i-1, k-2, \ell+2)\}$ in the term with coefficient $\frac{(k-1) \mu_1}{d(k-1,\ell+1)}$ in \eqref{eq:diff}. Notice that $v(i-1, k, \ell) \geq v(i-1, k-1, \ell+1)$ by the inductive hypothesis since $\ell <C_2$. Replacing the second minimum with an upper bound $v(i-1, k-1, \ell+1)$ yields that 
\begin{align}
    \lefteqn{\min\{v(i-1, k, \ell), v(i-1, k-1, \ell+1)\}}& \nonumber \\
    &\quad - \min\{v(i-1, k-1, \ell+1), v(i-1, k-2, \ell+2)\} \geq 0. \label{eq:larger_mu2_always_action1_mu1_terms}
\end{align}
It follows similarly that ($\ell-1<C_2$) the difference with coefficient $\frac{\min\{\ell,C_2\}\mu_2}{d(k,\ell)}$ in \eqref{eq:diff}
\begin{align}
    \lefteqn{\min\{v(i-1, k+1, \ell-1), v(i-1, k, \ell)\}}& \nonumber\\
    &\quad - \min\{v(i-1, k, \ell), v(i-1, k-1, \ell+1)\} \geq 0. \label{eq:larger_mu2_always_action1_mu2_terms}
\end{align} 
Applying both bounds in \eqref{eq:larger_mu2_always_action1_mu1_terms} and \eqref{eq:larger_mu2_always_action1_mu2_terms}, and \eqref{eq:diff_const_expand} with $\ell<C_2$ back in \eqref{eq:diff} results in
\begin{align*}
    \lefteqn{v(i,k,\ell) - v(i,k-1,\ell+1)}&\\
    &\quad \geq \frac{C_1\mu_1\mu_2}{d(k,\ell)d(k-1,\ell+1)} \left(\frac{ih_0}{C_1}\left(\frac{1}{\mu_1} - \frac{1}{\mu_2}\right) + \frac{h_1}{\mu_1} - \frac{h_2}{\mu_2} \right)\\
    &\quad \geq 0,
\end{align*}where the last inequality follows since $\mu_2 \geq \mu_1$ and $\frac{h_1}{\mu_1} \geq \frac{h_2}{\mu_2}$.

Clearly, if $\frac{h_1}{\mu_1} > \frac{h_2}{\mu_2}$ holds strictly,  then by a similar argument, we have $D(i,k,\ell) >0$ for all $i \geq 0$.
\end{proof}

\begin{proof} [Proof of Proposition \ref{prop:larger_mu2_dec}.]
We prove \eqref{eq:larger_mu2_dec} by induction on $i$. 
When $i=0$, recalling $\ell\geq C_2$ and using the values \eqref{eq:diff_bdry} and \eqref{eq:diff} with $\ell \geq C_2$ at $i=1$ yield
\begin{align}
    \lefteqn{v(0,k,\ell) - v(0,k-1,\ell+1) - v(1,k,\ell) + v(1,k-1,\ell+1)}& \nonumber \\
    & \quad = \frac{h_1}{\mu_1} - \frac{(\ell+1) h_2}{C_2\mu_2} - \left(\frac{-\mu_1h_0}{d(k,\ell)d(k-1,\ell+1)} + \frac{k h_1 + \ell h_2 }{d(k,\ell)} - \frac{(k-1)h_1+(\ell+1)h_2}{d(k-1,\ell+1)}\right) \nonumber \\
    &\qquad -\frac{(k-1) \mu_1}{d(k-1,\ell+1)} \Big[\min\{v(0, k, \ell), v(0, k-1, \ell+1)\} \nonumber \\
    &\qquad \qquad - \min\{v(0, k-1, \ell+1), v(0, k-2, \ell+2)\} \Big] \nonumber \\
    &\qquad - \frac{C_2\mu_2}{d(k,\ell)} \Big[\min\{v(0, k+1, \ell-1), v(0, k, \ell)\} - \min\{v(0, k, \ell), v(0, k-1, \ell+1)\}\Big] \nonumber \\
    &\quad \geq \frac{h_1}{\mu_1} - \frac{(\ell+1) h_2}{C_2\mu_2} - \left(\frac{k h_1 + \ell h_2 }{d(k,\ell)} - \frac{(k-1)h_1+(\ell+1)h_2}{d(k-1,\ell+1)}\right) \nonumber \\
    &\qquad -\frac{(k-1) \mu_1}{d(k-1,\ell+1)} \Big[\min\{v(0, k, \ell), v(0, k-1, \ell+1)\} \nonumber \\
    &\qquad \qquad - \min\{v(0, k-1, \ell+1), v(0, k-2, \ell+2)\} \Big] \nonumber \\
    &\qquad - \frac{C_2\mu_2}{d(k,\ell)} \Big[\min\{v(0, k+1, \ell-1), v(0, k, \ell)\} \nonumber\\
    & \qquad \qquad - \min\{v(0, k, \ell), v(0, k-1, \ell+1)\}\Big]\label{eq:larger_mu2_dec_inte1}
\end{align}
where the inequality follows by dropping the non-negative term proportional to $h_0$. 
Recall the following bound from \eqref{eq:larger_mu1_dec_inte4} on the term with coefficient $-\frac{(k-1) \mu_1}{d(k-1,\ell+1)}$,
\begin{align}
    &\min\{v(0, k, \ell), v(0, k-1, \ell+1)\} - \min\{v(0, k-1, \ell+1), v(0, k-2, \ell+2)\} \nonumber\\
    &\quad \leq \frac{h_1}{\mu_1} - \frac{(\ell+1)h_2}{C_2\mu_2}. \label{eq:bound-boundary0}
\end{align}
Similarly, the difference in the term with coefficient $- \frac{C_2\mu_2}{d(k,\ell)}$, 
\begin{align}
    &\min\{v(0, k+1, \ell-1), v(0, k, \ell)\} - \min\{v(0, k, \ell), v(0, k-1, \ell+1)\}\leq\frac{h_1}{\mu_1} - \frac{\ell h_2}{C_2\mu_2}.\label{eq:bound-boundary1}
\end{align}
Applying the bounds \eqref{eq:bound-boundary0} and \eqref{eq:bound-boundary1} in \eqref{eq:larger_mu2_dec_inte1} yields
\begin{align*}
    \lefteqn{v(0,k,\ell) - v(0,k-1,\ell+1) - v(1,k,\ell) + v(1,k-1,\ell+1)}&\\
    &\quad \geq \frac{h_1}{\mu_1} - \frac{(\ell+1) h_2}{C_2\mu_2} - \left(\frac{k h_1 + \ell h_2 }{d(k,\ell)} - \frac{(k-1)h_1+(\ell+1)h_2}{d(k-1,\ell+1)}\right)\\
    &\qquad -\frac{(k-1) \mu_1}{d(k-1,\ell+1)} \left(\frac{h_1}{\mu_1} - \frac{(\ell+1)h_2}{C_2\mu_2}\right) - \frac{C_2\mu_2}{d(k,\ell)} \left(\frac{h_1}{\mu_1} - \frac{\ell h_2}{C_2\mu_2}\right)\\
    &\quad = 0.
\end{align*}
Assume the result in \eqref{eq:larger_mu2_dec} holds at $i-1$ and verify it at $i \geq 1$. Using \eqref{eq:diff} in Statement \ref{state:diff} with $\ell \geq C_2$ at $i$ and $i+1$ yields
\begin{align}
    \lefteqn{v(i,k,\ell) - v(i,k-1,\ell+1) - v(i+1,k,\ell) + v(i+1,k-1,\ell+1)}& \nonumber \\
    & \quad = \frac{\mu_1h_0}{d(k,\ell)d(k-1,\ell+1)}\nonumber \\
    &\qquad +\frac{(k-1) \mu_1}{d(k-1,\ell+1)} \Big[\min\{v(i-1, k, \ell), v(i-1, k-1, \ell+1)\}\nonumber \\
    &\qquad \quad - \min\{v(i-1, k-1, \ell+1), v(i-1, k-2, \ell+2)\}\nonumber \\
    &\qquad \quad - \min\{v(i, k, \ell), v(i, k-1, \ell+1)\}\nonumber \\
    &\qquad \quad + \min\{v(i, k-1, \ell+1), v(i, k-2, \ell+2)\}\Big] \nonumber \\ 
    &\qquad + \frac{C_2\mu_2}{d(k,\ell)} \Big[\min\{v(i-1, k+1, \ell-1), v(i-1, k, \ell)\}\nonumber \\
    &\qquad \quad - \min\{v(i-1, k, \ell), v(i-1, k-1, \ell+1)\}\nonumber \\
    &\qquad \quad - \min\{v(i, k+1, \ell-1), v(i, k, \ell)\}\nonumber \\
    &\qquad \quad + \min\{v(i, k, \ell), v(i, k-1, \ell+1)\}\Big]. 
    \label{eq:larger_mu2_dec_inte2}
\end{align}
One can prove the non-negativity of the term with coefficient $\frac{(k-1) \mu_1}{d(k-1,\ell+1)}$ 
by a comparable argument as that showing the result with similar terms in \eqref{eq:larger_mu1_dec_inte1} 
(recall $\ell+1 \geq \ell \geq C_2$ so that the inductive hypotheses apply).

Now consider the four minima in the term with coefficient $\frac{C_2\mu_2}{d(k,\ell)}$ in \eqref{eq:larger_mu2_dec_inte2}. If $\ell-1 \geq C_2$, the non-negativity follows similarly as before, since the inductive hypotheses apply in this case. It suffices to consider if $\ell-1<C_2$ (this is the case when $\ell = C_2$ since we assumed $\ell \geq C_2$). Notice that $v(i-1, k+1, \ell-1) \geq v(i-1, k, \ell)$ and $v(i, k+1, \ell-1) \geq v(i, k, \ell)$ by Proposition \ref{prop:larger_mu2_always_action1}. There are two cases to examine based on whether $v(i, k, \ell) \leq v(i, k-1, \ell+1)$.
\begin{enumerate}[label= \textbf{Case} \arabic*:, leftmargin=3\parindent]
    \item If $v(i, k, \ell) \leq  v(i, k-1, \ell+1)$, choosing an upper bound $v(i-1, k, \ell)$ in the second minimum yields
    \begin{align*}
    \lefteqn{\min\{v(i-1, k+1, \ell-1), v(i-1, k, \ell)\}}&\\
    &\qquad - \min\{v(i-1, k, \ell), v(i-1, k-1, \ell+1)\}\\
    &\qquad - \min\{v(i, k+1, \ell-1), v(i, k, \ell)\}\\
    &\qquad + \min\{v(i, k, \ell), v(i, k-1, \ell+1)\}\\
    &\quad \geq v(i-1, k, \ell) -v(i-1, k, \ell) - v(i, k, \ell) + v(i, k, \ell)\\
    &\quad = 0.
\end{align*}
\item If $v(i, k, \ell) >  v(i, k-1, \ell+1)$, choosing an upper bound $v(i-1, k-1, \ell+1)$ in the second minimum yields
    \begin{align*}
    \lefteqn{\min\{v(i-1, k+1, \ell-1), v(i-1, k, \ell)\}}&\\
    &\qquad - \min\{v(i-1, k, \ell), v(i-1, k-1, \ell+1)\}\\
    &\qquad - \min\{v(i, k+1, \ell-1), v(i, k, \ell)\}\\
    &\qquad + \min\{v(i, k, \ell), v(i, k-1, \ell+1)\}\\
    &\quad \geq v(i-1, k, \ell) -v(i-1, k-1, \ell+1) - v(i, k, \ell) + v(i, k-1, \ell+1)\\
    &\quad \geq 0,
\end{align*}
where the last inequality follows from the inductive hypothesis since $\ell \geq C_2$.
\end{enumerate}
The non-negativity of the terms with coefficients $\frac{(k-1) \mu_1}{d(k-1,\ell+1)}$ and $\frac{C_2\mu_2}{d(k,\ell)}$ in \eqref{eq:larger_mu2_dec_inte2} proves the Statement \ref{state:larger_mu2_dec}. 

It remains to show the existence of the finite threshold $N_2(\ell)$. To do so, it suffices to show that the difference $D(i,k,\ell)$ as expressed in \eqref{eq:diff} by Statement \ref{state:diff} becomes negative when $i$ is sufficiently large, i.e., for all $i> N'_2(\ell)$ for some finite $N'_2(\ell)$. Without loss of generality assume $i\geq 1$. Following a similar argument as that proving the second statement of Proposition \ref{prop:larger_mu1_dec},
the following holds for all $i \geq 0$ and $\ell\geq C_2$, 
\begin{align}
     \lefteqn{v(i,k,\ell) - v(i,k-1,\ell+1)}& \nonumber \\
     & \quad \leq v(0,k,\ell) - v(0,k-1,\ell+1) \nonumber \\
     & \quad = \frac{h_1}{\mu_1} - \frac{\ell+1}{C_2}\frac{h_2}{\mu_2}, \label{eq:diff_larger_mu2_upbd}
\end{align}
so that the difference in the term with coefficient $\frac{(k-1) \mu_1}{d(k-1,\ell+1)}$ in \eqref{eq:diff}
\begin{align}
    \lefteqn{\min\{v(i-1, k, \ell), v(i-1, k-1, \ell+1)\}}& \nonumber \\
    &\qquad - \min\{v(i-1, k-1, \ell+1), v(i-1, k-2, \ell+2)\} \nonumber \\
    &\quad \leq \max\Big\{v(i-1, k, \ell) - v(i-1, k-1, \ell+1), \nonumber \\
    &\qquad \quad v(i-1, k-1, \ell+1) - v(i-1, k-2, \ell+2)\Big\} \nonumber \\
    &\quad \leq \frac{h_1}{\mu_1} - \frac{\ell+1}{C_2}\frac{h_2}{\mu_2},\label{eq:diff_larger_mu2_upbd_mu1_terms}
\end{align}
where the last inequality follows using \eqref{eq:diff_larger_mu2_upbd}. 

The difference in the term with coefficient $\frac{C_2\mu_2}{d(k,\ell)}$ (recall $\ell \geq C_2$) in \eqref{eq:diff} needs to be considered in two cases.
\begin{enumerate}[label= \textbf{Case} \arabic*:, leftmargin=3\parindent]
    \item If $\ell - 1 \geq C_2$, referring to \eqref{eq:diff_larger_mu2_upbd} again yields
        \begin{align}
            &\min\{v(i-1, k+1, \ell-1), v(i-1, k, \ell)\}- \min\{v(i-1, k, \ell), v(i-1, k-1, \ell+1)\} \nonumber \\
            &\quad \leq \frac{h_1}{\mu_1} - \frac{\ell}{C_2}\frac{h_2}{\mu_2}. \label{eq:diff_larger_mu2_upbd_mu2_terms}
        \end{align}
    \item If $\ell - 1 < C_2$, i.e., $\ell = C_2$, choosing an upper bound $v(i-1, k, \ell)$ in the first minimum yields
    \begin{align*}
            &\min\{v(i-1, k+1, \ell-1), v(i-1, k, \ell)\}- \min\{v(i-1, k, \ell), v(i-1, k-1, \ell+1)\}\\
            & \quad \leq v(i-1, k, \ell) - \min\{v(i-1, k, \ell), v(i-1, k-1, \ell+1)\}\\
            &\quad = \max\{0,v(i-1, k, \ell)- v(i-1, k-1, \ell+1)\}\\
            &\quad \leq \max\Big\{0, \frac{h_1}{\mu_1} - \frac{\ell+1}{C_2}\frac{h_2}{\mu_2}\Big\}\\
            &\quad \leq \frac{h_1}{\mu_1} - \frac{h_2}{\mu_2},
        \end{align*}
        where the second inequality uses \eqref{eq:diff_larger_mu2_upbd} and the last inequality follows by the assumption that $\frac{h_1}{\mu_1} \geq \frac{h_2}{\mu_2}$ and $\ell = C_2$.
\end{enumerate}
Hence, the expression given in \eqref{eq:diff_larger_mu2_upbd_mu2_terms} holds for both cases.
Using \eqref{eq:diff_larger_mu2_upbd_mu1_terms} and \eqref{eq:diff_larger_mu2_upbd_mu2_terms}, the difference $D(i,k,\ell)$ as given in \eqref{eq:diff} can be upper bounded by
\begin{align*}
    \lefteqn{v(i,k,\ell) - v(i,k-1,\ell+1)}&\\
    &\quad \leq \frac{-i\mu_1h_0}{d(k,\ell)d(k-1,\ell+1)} + \left(\frac{k h_1 + \ell h_2 }{d(k,\ell)} - \frac{(k-1)h_1+(\ell+1)h_2}{d(k-1,\ell+1)}\right) \nonumber \\
    &\qquad +\frac{(k-1) \mu_1}{d(k-1,\ell+1)} \left(\frac{h_1}{\mu_1} - \frac{\ell+1}{C_2}\frac{h_2}{\mu_2} \right) + \frac{C_2\mu_2}{d(k,\ell)} \left(\frac{h_1}{\mu_1} - \frac{\ell}{C_2}\frac{h_2}{\mu_2} \right) \\
    &\quad = \frac{-i\mu_1h_0}{d(k,\ell)d(k-1,\ell+1)} + \left(\frac{h_1}{\mu_1} - \frac{\ell+1}{C_2}\frac{h_2}{\mu_2} \right).
\end{align*}
The result follows by defining $N'_2(\ell) = \frac{c_4(\ell)}{c_3(\ell)}$, where $c_3(\ell) = \frac{\mu_1h_0}{d(k,\ell)d(k-1,\ell+1)} > 0$ and $c_4(\ell) = \frac{h_1}{\mu_1} - \frac{\ell+1}{C_2}\frac{h_2}{\mu_2}$ both are functions of $\ell$ (recall $k+\ell = C_1$) and independent of $i$.
\end{proof}

\begin{proof} [Proof of Proposition \ref{prop:larger_mu2_always_action0}]
    Proof by induction on $i$. The base case when $i=0$ follows by \eqref{eq:diff_bdry}:
    \begin{align*}
        D(0,k,\ell) = v(0,k,\ell) - v(0,k-1,\ell+1) = \frac{h_1}{\mu_1} - \frac{(\ell+1) h_2}{C_2\mu_2} \leq 0.
    \end{align*}
    Assume the \eqref{eq:larger_mu2_always_action0} holds at $i-1$ and consider it at $i \geq 1$. Recall the expansion of the Bellman equations of the difference $D(i,k,\ell)$ is given in \eqref{eq:diff}. 
    Notice that 
    $v(i-1, k-1, \ell+1) \leq v(i-1, k-2, \ell+2)$ by inductive hypothesis since $\ell+1 \geq C_2$. Choosing an upper bound $v(i-1, k-1, \ell+1)$ within the first minimum in the term with coefficient $\frac{(k-1) \mu_1}{d(k-1,\ell+1)}$ in \eqref{eq:diff} yields that 
\begin{align}
    \lefteqn{\min\{v(i-1, k, \ell), v(i-1, k-1, \ell+1)\}}& \nonumber \\
    &\quad - \min\{v(i-1, k-1, \ell+1), v(i-1, k-2, \ell+2)\} \leq 0. \label{eq:larger_mu2_always_action0_mu1_terms}
\end{align}
The difference in the term with coefficient $\frac{C_2\mu_2}{d(k,\ell)}$ in \eqref{eq:diff} is also non-positive:
\begin{align}
    \lefteqn{\min\{v(i-1, k+1, \ell-1), v(i-1, k, \ell)\}}& \nonumber \\
    &\quad - \min\{v(i-1, k, \ell), v(i-1, k-1, \ell+1)\} \leq 0, \label{eq:larger_mu2_always_action0_mu2_terms}
\end{align}following a comparable argument by inductive hypothesis.
It remains to check the one-step cost terms as presented in \eqref{eq:diff_const_expand}. Noticing $\ell \geq C_2$ in \eqref{eq:diff_const_expand} yields
\begin{align}
    \lefteqn{\frac{C_2\mu_1\mu_2}{d(k,\ell)d(k-1,\ell+1)}\left(\frac{-ih_0}{C_2 \mu_2} + \frac{h_1-h_2}{\mu_1} - \frac{C_1h_2}{C_2 \mu_2} \right)}& \nonumber \\
    &\quad \leq \frac{C_2\mu_1\mu_2}{d(k,\ell)d(k-1,\ell+1)}\left(\frac{h_1-h_2}{\mu_1} - \frac{C_1h_2}{C_2 \mu_2} \right) \nonumber \\
    &\quad \leq \frac{C_2\mu_1\mu_2}{d(k+1,\ell-1)d(k,\ell)}\left(\frac{h_1}{\mu_1} - \frac{h_2}{\mu_2}\right) \nonumber \\
    &\quad \leq 0, \label{eq:larger_mu2_always_action0_const}
\end{align}where the second inequality uses $C_2 \leq C_1$ and the last inequality holds by the assumption that $\frac{h_1}{\mu_1} \leq \frac{h_2}{\mu_2}$.
Applying all three bounds in \eqref{eq:larger_mu2_always_action0_mu1_terms} - \eqref{eq:larger_mu2_always_action0_const} back in \eqref{eq:diff} proves \eqref{eq:larger_mu2_always_action0}.
\end{proof}
\begin{proof} [Proof of Proposition \ref{prop:larger_mu2_inc}]
Proof by induction on $i$. When $i=0$, using the boundary values \eqref{eq:diff_bdry} and applying \eqref{eq:diff_const_expand} with $\ell < C_2$ in \eqref{eq:diff} at $i=1$ yield
\begin{align}
    \lefteqn{v(0,k,\ell) - v(0,k-1,\ell+1) - v(1,k,\ell) + v(1,k-1,\ell+1)}& \nonumber \\
    & \quad = \frac{h_1}{\mu_1} - \frac{h_2}{\mu_2} - \frac{C_1\mu_1\mu_2}{d(k,\ell)d(k-1,\ell+1)} \left(\frac{h_0}{C_1}\left(\frac{1}{\mu_1} - \frac{1}{\mu_2}\right) + \frac{h_1}{\mu_1} - \frac{h_2}{\mu_2} \right) \nonumber \\
    &\qquad -\frac{(k-1) \mu_1}{d(k-1,\ell+1)} \Big[\min\{v(0, k, \ell), v(0, k-1, \ell+1)\} \nonumber \\
    &\qquad \qquad - \min\{v(0, k-1, \ell+1), v(0, k-2, \ell+2)\} \Big] \nonumber \\
    &\qquad - \frac{\ell\mu_2}{d(k,\ell)} \Big[\min\{v(0, k+1, \ell-1), v(0, k, \ell)\} - \min\{v(0, k, \ell), v(0, k-1, \ell+1)\}\Big] \nonumber \\
    &\quad \leq \left(1- \frac{C_1\mu_1\mu_2}{d(k,\ell)d(k-1,\ell+1)}\right) \left(\frac{h_1}{\mu_1} - \frac{h_2}{\mu_2}\right) \nonumber \\
    &\qquad -\frac{(k-1) \mu_1}{d(k-1,\ell+1)} \Big[\min\{v(0, k, \ell), v(0, k-1, \ell+1)\} \nonumber \\
    &\qquad \qquad - \min\{v(0, k-1, \ell+1), v(0, k-2, \ell+2)\} \Big] \label{eq:larger_mu2_inc_inte1} \\
    &\qquad - \frac{\ell \mu_2}{d(k,\ell)} \Big[\min\{v(0, k+1, \ell-1), v(0, k, \ell)\} - \min\{v(0, k, \ell), v(0, k-1, \ell+1)\}\Big]\nonumber
\end{align}since $\mu_2\geq\mu_1$. 
Choosing an upper bound $v(0, k-1, \ell+1)$ within the second minimum in the second term with coefficient $-\frac{(k-1) \mu_1}{d(k-1,\ell+1)}$ in \eqref{eq:larger_mu2_inc_inte1} yields
\begin{align}
    &\min\{v(0, k, \ell), v(0, k-1, \ell+1)\} - \min\{v(0, k-1, \ell+1), v(0, k-2, \ell+2)\} \nonumber \\
    &\quad \geq \min\{v(0, k, \ell), v(0, k-1, \ell+1)\} -  v(0, k-1, \ell+1) \nonumber \\
    &\quad = \min\{v(0, k, \ell) - v(0, k-1, \ell+1), 0\} \nonumber \\
    &\quad = \min\{\frac{h_1}{\mu_1} - \frac{h_2}{\mu_2}, 0 \} \nonumber \\
    &\quad = \frac{h_1}{\mu_1} - \frac{h_2}{\mu_2},\label{eq:larger_mu2_inc_inte1_mu1_terms}
\end{align}where the second-to-last equality follows by \eqref{eq:diff_bdry} since $\ell < C_2$, and the last equality uses the assumption that $\frac{h_1}{\mu_1} \leq \frac{h_2}{\mu_2}$.
Using the boundary values \eqref{eq:diff_bdry} again with $\ell-1<\ell <C_2$ in the difference in the last term with coefficient $- \frac{\ell \mu_2}{d(k,\ell)}$ in \eqref{eq:larger_mu2_inc_inte1} gives
\begin{align}
    &\min\{v(0, k+1, \ell-1), v(0, k, \ell)\} - \min\{v(0, k, \ell), v(0, k-1, \ell+1)\} \nonumber \\
    &\quad \geq \min\{v(0, k+1, \ell-1) - v(0, k, \ell), v(0, k, \ell) - v(0, k-1, \ell+1)\} \nonumber \\
    &\quad =  \frac{h_1}{\mu_1} - \frac{h_2}{\mu_2}. \label{eq:larger_mu2_inc_inte1_mu2_terms}
\end{align}
Substituting the two bounds \ref{eq:larger_mu2_inc_inte1_mu1_terms} and \eqref{eq:larger_mu2_inc_inte1_mu2_terms} back in \eqref{eq:larger_mu2_inc_inte1} yields
\begin{align*}
    \lefteqn{v(0,k,\ell) - v(0,k-1,\ell+1) - v(1,k,\ell) + v(1,k-1,\ell+1)}&\\
    &\quad \leq \left(1- \frac{C_1\mu_1\mu_2}{d(k,\ell)d(k-1,\ell+1)} - \frac{(k-1) \mu_1}{d(k-1,\ell+1)} - \frac{C_2\mu_2}{d(k,\ell)}\right) \left(\frac{h_1}{\mu_1} - \frac{h_2}{\mu_2}\right)\\
    &\quad = 0,
\end{align*}where the last step follows by recalling the definitions of $p_k,q_k$ and $r_k$ in \eqref{def:pk}, \eqref{def:qk} and \eqref{def:rk}, respectively, and using the identity in \eqref{eq:identity}.
Assume \eqref{eq:larger_mu2_inc} holds at $i-1$ to prove it at $i \geq 1$. Using \eqref{eq:diff} with $\ell < C_2$ at $i$ and $i+1$ yields
    \begin{align}
    \lefteqn{v(i,k,\ell) - v(i,k-1,\ell+1) - v(i+1,k,\ell) + v(i+1,k-1,\ell+1)}& \nonumber \\
    & \quad = \frac{-(\mu_2-\mu_1)h_0}{d(k,\ell)d(k-1,\ell+1)}\nonumber \\
    &\qquad +\frac{(k-1) \mu_1}{d(k-1,\ell+1)} \Big[\min\{v(i-1, k, \ell), v(i-1, k-1, \ell+1)\}\nonumber \\
    &\qquad \quad - \min\{v(i-1, k-1, \ell+1), v(i-1, k-2, \ell+2)\}\nonumber \\
    &\qquad \quad - \min\{v(i, k, \ell), v(i, k-1, \ell+1)\}\nonumber \\
    &\qquad \quad + \min\{v(i, k-1, \ell+1), v(i, k-2, \ell+2)\}\Big] \nonumber \\
    &\qquad + \frac{\ell\mu_2}{d(k,\ell)} \Big[\min\{v(i-1, k+1, \ell-1), v(i-1, k, \ell)\}\nonumber \\
    &\qquad \quad - \min\{v(i-1, k, \ell), v(i-1, k-1, \ell+1)\}\nonumber \\
    &\qquad \quad - \min\{v(i, k+1, \ell-1), v(i, k, \ell)\}\nonumber \\
    &\qquad \quad + \min\{v(i, k, \ell), v(i, k-1, \ell+1)\}\Big]. \label{eq:larger_mu2_inc_inte}
\end{align}
Consider first the four minima in the second term with coefficient $\frac{(k-1) \mu_1}{d(k-1,\ell+1)}$ in \eqref{eq:larger_mu2_inc_inte} based on whether $v(i-1, k+1, \ell-1) > v(i-1, k, \ell)$. 
\begin{enumerate}[label= \textbf{Case} \arabic*:, leftmargin=3\parindent]
    \item If $v(i-1, k, \ell) \leq v(i-1, k-1, \ell+1)$, replace the last minimum with an upper bound $v(i, k, \ell)$ and consider the following two subcases.
    \begin{enumerate}[label= \textbf{Subcase} \alph*:, leftmargin=3\parindent]
        \item If $v(i, k+1, \ell-1) > v(i, k, \ell)$, choosing an upper bound $v(i-1, k, \ell)$ in the first minimum yields
        \begin{align*}
            \lefteqn{\min\{v(i-1, k+1, \ell-1), v(i-1, k, \ell)\}}& \nonumber \\
            &\qquad - \min\{v(i-1, k, \ell), v(i-1, k-1, \ell+1)\}\nonumber \\
            &\qquad - \min\{v(i, k+1, \ell-1), v(i, k, \ell)\}\nonumber \\
            &\qquad + \min\{v(i, k, \ell), v(i, k-1, \ell+1)\}\\
            &\quad \leq v(i-1, k, \ell) - v(i-1, k, \ell) - v(i, k, \ell) + v(i, k, \ell)\\
            & \quad = 0.
        \end{align*}
        \item If $v(i, k+1, \ell-1) \leq v(i, k, \ell)$, choosing an upper bound $v(i-1, k+1, \ell-1)$ in the first minimum yields
        \begin{align}
            \lefteqn{\min\{v(i-1, k+1, \ell-1), v(i-1, k, \ell)\}}& \nonumber \\
            &\qquad - \min\{v(i-1, k, \ell), v(i-1, k-1, \ell+1)\}\nonumber \\
            &\qquad - \min\{v(i, k+1, \ell-1), v(i, k, \ell)\}\nonumber \\
            &\qquad + \min\{v(i, k, \ell), v(i, k-1, \ell+1)\} \nonumber \\
            &\quad \leq v(i-1, k+1, \ell-1) - v(i-1, k, \ell) - v(i, k+1, \ell-1) + v(i, k, \ell) \nonumber \\
            & \quad \leq 0, \label{eq:larger_mu2_inc_inte3}
        \end{align}where the last step holds by inductive hypothesis since $\ell - 1 < C_2$.
    \end{enumerate}
    \item If $v(i-1, k, \ell) > v(i-1, k-1, \ell+1)$, replace the last minimum with an upper bound $v(i, k-1, \ell+1)$ and again consider the following two subcases.
    \begin{enumerate}[label= \textbf{Subcase} \alph*:, leftmargin=3\parindent]
        \item If $v(i, k+1, \ell-1) > v(i, k, \ell)$, choosing an upper bound $v(i-1, k, \ell)$ in the first minimum yields
        \begin{align}
            \lefteqn{\min\{v(i-1, k+1, \ell-1), v(i-1, k, \ell)\}}& \nonumber \\
            &\qquad - \min\{v(i-1, k, \ell), v(i-1, k-1, \ell+1)\}\nonumber \\
            &\qquad - \min\{v(i, k+1, \ell-1), v(i, k, \ell)\}\nonumber \\
            &\qquad + \min\{v(i, k, \ell), v(i, k-1, \ell+1)\}\nonumber \\
            &\quad \leq v(i-1, k, \ell) - v(i-1, k-1, \ell+1) - v(i, k, \ell) + v(i, k-1, \ell+1)\nonumber \\
            & \quad \leq 0, \label{eq:larger_mu2_inc_inte4}
        \end{align}where the last step follows from inductive hypothesis since $\ell < C_2$.
        \item If $v(i, k+1, \ell-1) \leq v(i, k, \ell)$, choosing an upper bound $v(i-1, k+1, \ell-1)$ in the first minimum yields
        \begin{align*}
            \lefteqn{\min\{v(i-1, k+1, \ell-1), v(i-1, k, \ell)\}}& \nonumber \\
            &\qquad - \min\{v(i-1, k, \ell), v(i-1, k-1, \ell+1)\}\nonumber \\
            &\qquad - \min\{v(i, k+1, \ell-1), v(i, k, \ell)\}\nonumber \\
            &\qquad + \min\{v(i, k, \ell), v(i, k-1, \ell+1)\} \nonumber \\
            &\quad \leq v(i-1, k+1, \ell-1) - v(i-1, k-1, \ell+1) \nonumber \\
            &\quad \qquad - v(i, k+1, \ell-1) + v(i, k-1, \ell+1) \nonumber \\
            & \quad \leq 0,
        \end{align*}where the last step adds the two inequalities \eqref{eq:larger_mu2_inc_inte3} and \eqref{eq:larger_mu2_inc_inte4}.
    \end{enumerate}
\end{enumerate}
Consider now the four minima in the last term in \eqref{eq:larger_mu2_inc_inte} with coefficient $\frac{\ell\mu_2}{d(k,\ell)}$. If $\ell+1 < C_2$, an analogous reasoning as before applies to prove its non-positivity since all the inductive hypotheses apply. It remains to verify when $\ell+1 \geq C_2$ (in fact, we can only have $\ell = C_2-1$ since we assumed $\ell < C_2$). Proposition \ref{prop:larger_mu2_always_action0} guarantees that $v(i-1, k-1, \ell+1) \leq v(i-1, k-2, \ell+2)$ and $v(i, k-1, \ell+1) \leq v(i, k-2, \ell+2)$, which leaves two cases to analyze based on whether $v(i, k, \ell) \leq v(i, k-1, \ell+1)$.
\begin{enumerate}[label= \textbf{Case} \arabic*:, leftmargin=3\parindent]
    \item If $v(i, k, \ell) \leq v(i, k-1, \ell+1)$, choosing an upper bound $v(i-1, k, \ell)$ in the first minimum yields
    \begin{align*}
    \lefteqn{\min\{v(i-1, k, \ell), v(i-1, k-1, \ell+1)\}}&\nonumber \\
    &\qquad - \min\{v(i-1, k-1, \ell+1), v(i-1, k-2, \ell+2)\}\nonumber \\
    &\qquad - \min\{v(i, k, \ell), v(i, k-1, \ell+1)\}\nonumber \\
    &\qquad + \min\{v(i, k-1, \ell+1), v(i, k-2, \ell+2)\}\\
    &\quad \leq v(i-1, k, \ell) -v(i-1, k-1, \ell+1) - v(i, k, \ell) + v(i, k-1, \ell+1)\\
    &\quad \leq 0,
\end{align*}where the last inequality holds by the inductive hypothesis since $\ell < C_2$.
\item If $v(i, k, \ell) >  v(i, k-1, \ell+1)$, choosing an upper bound $v(i-1, k-1, \ell+1)$ in the first minimum yields
    \begin{align*}
    \lefteqn{\min\{v(i-1, k+1, \ell-1), v(i-1, k, \ell)\}}&\\
    &\qquad - \min\{v(i-1, k, \ell), v(i-1, k-1, \ell+1)\}\\
    &\qquad - \min\{v(i, k+1, \ell-1), v(i, k, \ell)\}\\
    &\qquad + \min\{v(i, k, \ell), v(i, k-1, \ell+1)\}\\
    &\quad \leq v(i-1, k-1, \ell+1) -v(i-1, k-1, \ell+1)\\
    &\qquad - v(i, k-1, \ell+1) + v(i, k-1, \ell+1)\\
    &\quad = 0.
\end{align*}
\end{enumerate}
The non-positivity of both two terms with coefficients $\frac{(k-1) \mu_1}{d(k-1,\ell+1)}$ and $\frac{\ell\mu_2}{d(k,\ell)}$ in \eqref{eq:larger_mu2_inc_inte}
together with the assumption that $\mu_2 \geq \mu_1$ proves the desired non-decreasing result. 

It remains to show the second statement regarding the existence of the finite threshold $N_3(\ell)$ with the additional assumption that $\mu_2>\mu_1$. It suffices to demonstrate that the difference in \eqref{eq:diff} becomes positive when $i$ is sufficiently large, i.e., for all $i> N'_3(\ell)$ for some finite $N'_3(\ell)$. Hence without loss of generality assume $i\geq 1$ for the consideration of large $i$ only.
The monotonicity ensures that the following holds for all $i \geq 0$ and $\ell< C_2$:
\begin{align}
     \lefteqn{v(i,k,\ell) - v(i,k-1,\ell+1)}& \nonumber \\
     & \quad \geq v(0,k,\ell) - v(0,k-1,\ell+1) \nonumber \\
     & \quad = \frac{h_1}{\mu_1} - \frac{h_2}{\mu_2}. \label{eq:diff_larger_mu2_lobd}
\end{align}
Using \eqref{eq:diff_larger_mu2_lobd}, the difference in \eqref{eq:diff} with coefficient $\frac{\ell\mu_2}{d(k,\ell)}$ is lower bounded by
\begin{align}
    &\min\{v(i-1, k+1, \ell-1), v(i-1, k, \ell)\}- \min\{v(i-1, k, \ell), v(i-1, k-1, \ell+1)\} \nonumber \\
    &\quad \geq \min\{v(i-1, k+1, \ell-1) - v(i-1, k, \ell), v(i-1, k, \ell) - v(i-1, k-1, \ell+1)\} \nonumber \\
    &\quad \geq \frac{h_1}{\mu_1} - \frac{h_2}{\mu_2}. \label{eq:diff_larger_mu2_lobd_mu2_terms}
\end{align}
Following a comparable reasoning as \eqref{eq:larger_mu2_inc_inte1_mu1_terms} with the substitution of \eqref{eq:diff_larger_mu2_lobd} (rather than the boundary values \eqref{eq:diff_bdry} in \eqref{eq:larger_mu2_inc_inte1_mu1_terms}), the difference in the term with coefficient $\frac{(k-1) \mu_1}{d(k-1,\ell+1)}$ in \eqref{eq:diff} is also lower bounded by
 \begin{align}
    \lefteqn{\min\{v(i-1, k, \ell), v(i-1, k-1, \ell+1)\}}& \nonumber \\
    &\qquad - \min\{v(i-1, k-1, \ell+1), v(i-1, k-2, \ell+2)\} \nonumber \\
    &\quad \geq \frac{h_1}{\mu_1} - \frac{h_2}{\mu_2}. \label{eq:diff_larger_mu2_lobd_mu1_terms}
\end{align}

Therefore, the difference in \eqref{eq:diff} (applying \eqref{eq:diff_const_expand} with $\ell <C_2$) can be further lower bounded by
\begin{align*}
    \lefteqn{v(i,k,\ell) - v(i,k-1,\ell+1)}&\\
    &\quad \geq \frac{C_1\mu_1\mu_2}{d(k,\ell)d(k-1,\ell+1)} \left(\frac{ih_0}{C_1}\left(\frac{1}{\mu_1} - \frac{1}{\mu_2}\right) + \frac{h_1}{\mu_1} - \frac{h_2}{\mu_2} \right) \nonumber \\
    &\qquad +\left(\frac{(k-1) \mu_1}{d(k-1,\ell+1)} + \frac{\ell \mu_2}{d(k,\ell)}\right) \left(\frac{h_1}{\mu_1} -\frac{h_2}{\mu_2} \right)\\
    &\quad = \frac{i(\mu_2-\mu_1)h_0}{d(k,\ell)d(k-1,\ell+1)} + \left(\frac{h_1}{\mu_1} - \frac{h_2}{\mu_2} \right),
\end{align*}where the equality relies on the definitions of $p_k,q_k$ and $r_k$ in \eqref{def:pk}, \eqref{def:qk} and \eqref{def:rk}, respectively, and the identity in \eqref{eq:identity}.
The result follows by defining $N'_3(\ell) = \frac{c_6(\ell)}{c_5(\ell)}$, where $c_5(\ell) = \frac{(\mu_2-\mu_1)h_0}{d(k,\ell)d(k-1,\ell+1)} > 0$ and $c_6(\ell) = -\left(\frac{h_1}{\mu_1} - \frac{h_2}{\mu_2} \right) \geq 0$ both are functions of $\ell$ (recall $k+\ell = C_1$) and independent of $i$.
\end{proof}

\begin{proof} [Proof of Corollary \ref{cor:always_action0}.]
    We discuss the proof in two cases based on the assumptions.
    \begin{enumerate}[label= \textbf{Case} \arabic*:, leftmargin=3\parindent]
        \item If $\mu_1 \geq \mu_2$, the result is a direct application of Proposition \ref{prop:larger_mu1_dec} and the boundary values as in \eqref{eq:diff_bdry}.
        \item If $\mu_1 \leq \mu_2$ and $\ell \geq C_2$, then the assumption $\frac{h_1}{\mu_1} \leq \frac{\max\{\ell+1,C_2\}}{C_2}\frac{h_2}{\mu_2}$ says $\frac{h_1}{\mu_1} \leq \frac{\ell+1}{C_2}\frac{h_2}{\mu_2}$, the result follows from Proposition \ref{prop:larger_mu2_dec} and the boundary values in \eqref{eq:diff_bdry} if $\frac{h_1}{\mu_1} \geq \frac{h_2}{\mu_2}$, and from Proposition \ref{prop:larger_mu2_always_action0} if $\frac{h_1}{\mu_1} \leq \frac{h_2}{\mu_2}$.
    \end{enumerate}
\end{proof}

\begin{proof} [Proof of Proposition \ref{prop:mono_ell}.]
    Proof by induction on $i$. Consider first $\mu_1\geq \mu_2$. Using \eqref{eq:diff_bdry}, inequality \eqref{eq:mono_ell} at $i = 0$ becomes
    \begin{align}
        & v(0,k+1,\ell-1) - v(0,k,\ell) - \Big[v(0,k,\ell) - v(0,k-1,\ell+1)\Big] \nonumber \\
        & \quad = \left(\frac{h_1}{\mu_1} - \frac{\max\{\ell,C_2\}}{C_2}\frac{h_2}{\mu_2} \right) - \left(\frac{h_1}{\mu_1} - \frac{\max\{\ell+1,C_2\}}{C_2}\frac{h_2}{\mu_2} \right) \nonumber \\
        & \quad \geq 0. \label{eq:mono_ell_larger_mu1_base}
    \end{align}
    Assume the result holds at $i-1$, consider it at $i\geq 1$ (hence $x \in \tilde{\X}_{diff}$) based on three cases $\ell \geq C_2$, $\ell<C_2$, and $\ell = C_2$.
    \begin{enumerate}[label= \textbf{Case} \arabic*:, leftmargin=3\parindent]
        \item If $\ell > C_2$, it holds that $d(k+1,\ell-1) > d(k,\ell) > d(k-1,\ell+1)$.
        Thinning the MDP with a higher rate $d(k+1,\ell-1) = (k+1)\mu_1 + C_2 \mu_2$ yields
            \begin{align}
            & v(i,k+1,\ell-1) - v(i,k,\ell) - \Big[v(i,k,\ell) - v(i,k-1,\ell+1)\Big] \nonumber \\
            &\quad =\frac{(k-1) \mu_1}{d(k+1,\ell-1)} \Big[\min\{v(i-1, k+1, \ell-1), v(i-1, k, \ell)\} \nonumber \\
            &\qquad \quad - \min\{v(i-1, k, \ell), v(i-1, k-1, \ell+1)\} \nonumber \\
            &\qquad \quad - \min\{v(i-1, k, \ell), v(i-1, k-1, \ell+1)\} \nonumber \\
            &\qquad \quad + \min\{v(i-1, k-1, \ell+1), v(i-1, k-2, \ell+2)\} \Big] \nonumber \\
            &\qquad + \frac{C_2\mu_2}{d(k+1,\ell-1)} \Big[\min\{v(i-1, k+2, \ell-2), v(i-1, k+1, \ell-1)\} \nonumber \\
            &\qquad \quad - \min\{v(i-1, k+1, \ell-1), v(i-1, k, \ell)\} \nonumber \\
            &\qquad \quad - \min\{v(i-1, k+1, \ell-1), v(i-1, k, \ell)\} \nonumber \\
            &\qquad \quad + \min\{v(i-1, k, \ell), v(i-1, k-1, \ell+1)\}\Big] \nonumber \\
            &\quad + \frac{2\mu_1}{d(k+1,\ell-1)}\Big[\min\{v(i-1, k+1, \ell-1), v(i-1, k, \ell)\} \nonumber \\
            &\qquad \quad - \min\{v(i-1, k, \ell), v(i-1, k-1, \ell+1)\} \nonumber \\
            &\qquad \quad -v(i,k,\ell) + v(i,k-1,\ell+1)\Big]. 
            \label{eq:mono_ell_larger_mu1_larger_ell}
        \end{align} 
Consider the four minima in the first term in \eqref{eq:mono_ell_larger_mu1_larger_ell} with coefficient $\frac{(k-1) \mu_1}{d(k+1,\ell-1)}$
based on whether $v(i-1, k+1, \ell-1) \leq v(i-1, k, \ell)$ or the reverse inequality holds.
\begin{enumerate}[label= \textbf{Case} \alph*:, leftmargin=3\parindent]
    \item If $v(i-1, k+1, \ell-1) > v(i-1, k, \ell)$, choose an upper bound $v(i-1,k-1,\ell+1)$ in the third minimum and consider the two sub-cases.
    \begin{enumerate}[label= \textbf{Subcase} (\roman*):, leftmargin=3\parindent]
        \item If $v(i-1, k-1, \ell+1) \leq v(i-1, k-2, \ell+2)$, replacing the second minimum with an upper bound $v(i-1,k,\ell)$ yields
        \begin{align*}
            \lefteqn{\min\{v(i-1, k+1, \ell-1), v(i-1, k, \ell)\}}& \\
            &\qquad - \min\{v(i-1, k, \ell), v(i-1, k-1, \ell+1)\}\\
            &\qquad - \min\{v(i-1, k, \ell), v(i-1, k-1, \ell+1)\} \\
            &\qquad + \min\{v(i-1, k-1, \ell+1), v(i-1, k-2, \ell+2)\} \\
            &\quad \geq v(i-1, k, \ell) - v(i-1,k,\ell)\\
            &\qquad - v(i-1,k-1,\ell+1) + v(i-1, k-1, \ell+1) \\
            &\quad = 0.
        \end{align*}
        \item If $v(i-1, k-1, \ell+1) > v(i-1, k-2, \ell+2)$, replacing the second minimum with an upper bound $v(i-1,k-1,\ell+1)$ yields
        \begin{align}
            \lefteqn{\min\{v(i-1, k+1, \ell-1), v(i-1, k, \ell)\}}& \nonumber \\
            &\qquad - \min\{v(i-1, k, \ell), v(i-1, k-1, \ell+1)\} \nonumber \\
            &\qquad - \min\{v(i-1, k, \ell), v(i-1, k-1, \ell+1)\} \nonumber \\
            &\qquad + \min\{v(i-1, k-1, \ell+1), v(i-1, k-2, \ell+2)\} \nonumber \\
            &\quad \geq v(i-1, k, \ell) - v(i-1,k-1,\ell+1) \nonumber \\
            &\qquad - v(i-1,k-1,\ell+1) + v(i-1, k-2, \ell+2) \nonumber \\
            &\quad \geq 0, \label{eq:mono_ell_inte1}
        \end{align}where the last inequality holds by inductive hypothesis.
    \end{enumerate}
    \item If $v(i-1, k+1, \ell-1) \leq v(i-1, k, \ell)$, choose an upper bound $v(i-1,k,\ell)$ in the third minimum and consider the two sub-cases.
    \begin{enumerate}
        \item If $v(i-1, k-1, \ell+1) \leq v(i-1, k-2, \ell+2)$, replacing the second minimum with an upper bound $v(i-1,k,\ell)$ yields
        \begin{align}
            \lefteqn{\min\{v(i-1, k+1, \ell-1), v(i-1, k, \ell)\}}& \nonumber \\
            &\qquad - \min\{v(i-1, k, \ell), v(i-1, k-1, \ell+1)\} \nonumber \\
            &\qquad - \min\{v(i-1, k, \ell), v(i-1, k-1, \ell+1)\} \nonumber \\
            &\qquad + \min\{v(i-1, k-1, \ell+1), v(i-1, k-2, \ell+2)\} \nonumber \\
            &\quad \geq v(i-1, k+1, \ell-1) - v(i-1,k,\ell) \nonumber \\
            &\qquad - v(i-1,k,\ell) + v(i-1, k-1, \ell+1) \nonumber \\
            &\quad \geq 0, \label{eq:mono_ell_inte2}
        \end{align}where the last inequality holds by inductive hypothesis.
        \item If $v(i-1, k-1, \ell+1) > v(i-1, k-2, \ell+2)$, replacing the second minimum with an upper bound $v(i-1,k-1,\ell+1)$ yields
        \begin{align*}
            \lefteqn{\min\{v(i-1, k+1, \ell-1), v(i-1, k, \ell)\}}& \\
            &\qquad - \min\{v(i-1, k, \ell), v(i-1, k-1, \ell+1)\}\\
            &\qquad - \min\{v(i-1, k, \ell), v(i-1, k-1, \ell+1)\} \\
            &\qquad + \min\{v(i-1, k-1, \ell+1), v(i-1, k-2, \ell+2)\} \\
            &\quad \geq v(i-1, k+1, \ell-1) - v(i-1,k-1,\ell+1)\\
            &\qquad - v(i-1,k,\ell) + v(i-1, k-2, \ell+2) \\
            &\quad \geq 0,
        \end{align*}where the last inequality follows by adding \eqref{eq:mono_ell_inte1} and \eqref{eq:mono_ell_inte2}.
    \end{enumerate}
\end{enumerate}
A comparable argument can be found to show the second term in \eqref{eq:mono_ell_larger_mu1_larger_ell} with coefficient $\frac{C_2\mu_2}{d(k+1,\ell-1)}$
is also non-negative.
Now consider the last term in \eqref{eq:mono_ell_larger_mu1_larger_ell} with coefficient $\frac{2\mu_1}{d(k+1,\ell-1)}$. 
There are two cases to examine based on whether $v(i-1, k+1, \ell-1) \leq v(i-1, k, \ell)$ or not.
\begin{enumerate}[label= \textbf{Case} \alph*:, leftmargin=3\parindent]
    \item If $v(i-1, k+1, \ell-1) > v(i-1, k, \ell)$, replacing the second minimum with an upper bound $(i-1,k-1,\ell+1)$ yields
    \begin{align*}
        \lefteqn{\min\{v(i-1, k+1, \ell-1), v(i-1, k, \ell)\}}& \\
        &\qquad - \min\{v(i-1, k, \ell), v(i-1, k-1, \ell+1)\}\\
        &\qquad - v(i,k,\ell) + v(i,k-1,\ell+1)\\
        & \quad \geq v(i-1, k, \ell) - v(i-1, k-1, \ell+1) - v(i,k,\ell) + v(i,k-1,\ell+1) \\
        & \quad \geq 0,
    \end{align*}where the last inequality follows by Proposition \ref{prop:larger_mu1_dec}.
    \item If $v(i-1, k+1, \ell-1) \leq v(i-1, k, \ell)$, replacing the second minimum with an upper bound $(i-1,k,\ell)$ yields
    \begin{align*}
        \lefteqn{\min\{v(i-1, k+1, \ell-1), v(i-1, k, \ell)\}}& \\
        &\qquad - \min\{v(i-1, k, \ell), v(i-1, k-1, \ell+1)\}\\
        &\qquad - v(i,k,\ell) + v(i,k-1,\ell+1)\\
        & \quad \geq v(i-1, k+1, \ell-1) - v(i-1, k, \ell) - v(i,k,\ell) + v(i,k-1,\ell+1)\\
        & \quad = \Big[v(i-1, k+1, \ell-1) - v(i-1, k, \ell)\\
        & \qquad \quad - v(i-1,k,\ell) + v(i-1,k-1,\ell+1)\Big]\\
        & \qquad + \Big[v(i-1,k,\ell) - v(i-1,k-1,\ell+1)\\
        & \qquad \quad - v(i,k,\ell) + v(i,k-1,\ell+1)\Big]\\
        & \quad \geq 0,
    \end{align*}where the last inequality holds by the inductive hypothesis and Proposition \ref{prop:larger_mu1_dec}.
    
    A substitution of the non-negativity of the three terms with coefficients $\frac{(k-1) \mu_1}{d(k+1,\ell-1)}$, $\frac{C_2\mu_2}{d(k+1,\ell-1)}$, and $\frac{2\mu_1}{d(k+1,\ell-1)}$ in \eqref{eq:mono_ell_larger_mu1_larger_ell}
    establishes the result for $\ell > C_2$.
    \end{enumerate}
    \item If $\ell = C_2$, it still holds that $d(k+1,\ell-1) > d(k,\ell) > d(k-1,\ell+1)$ since $\mu_1\geq\mu_2$. Thinning the MDP with a higher rate $d(k+1,\ell-1) = (k+1)\mu_1 + (C_2-1)\mu_2$ yields
    \begin{align}
        & v(i,k+1,\ell-1) - v(i,k,\ell) - \Big[v(i,k,\ell) - v(i,k-1,\ell+1)\Big] \nonumber \\
        &\quad =\frac{(k-1) \mu_1}{d(k+1,\ell-1)} \Big[\min\{v(i-1, k+1, \ell-1), v(i-1, k, \ell)\} \nonumber \\
        &\qquad \quad - \min\{v(i-1, k, \ell), v(i-1, k-1, \ell+1)\} \nonumber \\
        &\qquad \quad - \min\{v(i-1, k, \ell), v(i-1, k-1, \ell+1)\} \nonumber \\
        &\qquad \quad + \min\{v(i-1, k-1, \ell+1), v(i-1, k-2, \ell+2)\} \Big] \nonumber \\
        &\qquad + \frac{(C_2-1)\mu_2}{d(k+1,\ell-1)} \Big[\min\{v(i-1, k+2, \ell-2), v(i-1, k+1, \ell-1)\} \nonumber \\
        &\qquad \quad - \min\{v(i-1, k+1, \ell-1), v(i-1, k, \ell)\} \nonumber \\
        &\qquad \quad - \min\{v(i-1, k+1, \ell-1), v(i-1, k, \ell)\} \nonumber \\
        &\qquad \quad + \min\{v(i-1, k, \ell), v(i-1, k-1, \ell+1)\}\Big] \nonumber \\
        &\qquad + \frac{\mu_2}{d(k+1,\ell-1)} \Big[\min\{v(i-1, k+1, \ell-1), v(i-1, k, \ell)\} \nonumber \\
        &\qquad \quad - 2\min\{v(i-1, k+1, \ell-1), v(i-1, k, \ell)\} \nonumber \\
        &\qquad \quad + \min\{v(i-1, k, \ell), v(i-1, k-1, \ell+1)\}\Big] \nonumber \\
        &\qquad + \frac{\mu_1}{d(k+1,\ell-1)} \Big[\min\{v(i-1, k+1, \ell-1), v(i-1, k, \ell)\} \nonumber \\
        &\qquad \quad - 2\min\{v(i-1, k, \ell), v(i-1, k-1, \ell+1)\} + v(i,k-1,\ell+1) \Big] \nonumber \\
        &\qquad + \frac{(\mu_1-\mu_2)}{d(k+1,\ell-1)}\Big[\min\{v(i-1, k+1, \ell-1), v(i-1, k, \ell)\} \nonumber \\
        &\qquad \quad -2v(i,k,\ell) + v(i,k-1,\ell+1)\Big]. \label{eq:mono_ell_larger_mu1_equal}
    \end{align}
    The non-negativity of the first two terms with coefficients $\frac{(k-1) \mu_1}{d(k+1,\ell-1)}$ and $\frac{(C_2-1)\mu_2}{d(k+1,\ell-1)}$ in \eqref{eq:mono_ell_larger_mu1_equal} follows by a similar argument as the one used for the term with coefficient $\frac{(k-1) \mu_1}{d(k+1,\ell-1)}$ in the proof for the case when $\ell > C_2$. It remains to consider the last three terms in \eqref{eq:mono_ell_larger_mu1_equal}.
    We now decompose the coefficient in the second-to-last term: 
    \begin{align*}
        \frac{\mu_1}{d(k+1,\ell-1)} = \frac{\mu_2}{d(k+1,\ell-1)} + \frac{\mu_1 - \mu_2}{d(k+1,\ell-1)}.
    \end{align*}
    We then redistribute the second-to-last term with these new coefficients above, adding them to the third and last terms in \eqref{eq:mono_ell_larger_mu1_equal} with coefficients $\frac{\mu_2}{d(k+1,\ell-1)}$ and $\frac{\mu_1 - \mu_2}{d(k+1,\ell-1)}$, respectively, while preserving their common coefficient. This allows us to rearrange the sum of the last three terms as follows.
        \begin{align}
        \lefteqn{\frac{\mu_2}{d(k+1,\ell-1)} \Big[\min\{v(i-1, k+1, \ell-1), v(i-1, k, \ell)\}}&\nonumber \\
        &\qquad \quad - 2\min\{v(i-1, k+1, \ell-1), v(i-1, k, \ell)\} \nonumber \\
        &\qquad \quad + \min\{v(i-1, k, \ell), v(i-1, k-1, \ell+1)\}\Big] \nonumber \\
        &\qquad + \frac{\mu_1}{d(k+1,\ell-1)} \Big[\min\{v(i-1, k+1, \ell-1), v(i-1, k, \ell)\} \nonumber \\
        &\qquad \quad - 2\min\{v(i-1, k, \ell), v(i-1, k-1, \ell+1)\} + v(i,k-1,\ell+1) \Big] \nonumber \\
        &\qquad + \frac{(\mu_1-\mu_2)}{d(k+1,\ell-1)}\Big[\min\{v(i-1, k+1, \ell-1), v(i-1, k, \ell)\} \nonumber \\
        &\qquad \quad -2v(i,k,\ell) + v(i,k-1,\ell+1)\Big] \nonumber \\
        &\quad = \frac{\mu_2}{d(k+1,\ell-1)} \Big[v(i,k-1,\ell+1) - \min\{v(i-1, k, \ell), v(i-1, k-1, \ell+1)\}\Big] \nonumber \\
        &\qquad + \frac{2(\mu_1-\mu_2)}{d(k+1,\ell-1)}\Big[\min\{v(i-1, k+1, \ell-1), v(i-1, k, \ell)\} \label{eq:mono_ell_larger_mu1_equal_inte} \\
        &\qquad \quad - \min\{v(i-1, k, \ell), v(i-1, k-1, \ell+1)\}-v(i,k,\ell) + v(i,k-1,\ell+1)\Big] \nonumber
    \end{align}
    The non-negativity of the last term in \eqref{eq:mono_ell_larger_mu1_equal_inte}, with coefficient $\frac{2(\mu_1-\mu_2)}{d(k+1,\ell-1)}$, follows by the same reasoning used to establish the non-negativity of the term with coefficient $\frac{2\mu_1}{d(k+1,\ell-1)}$ in \eqref{eq:mono_ell_larger_mu1_larger_ell} for $\ell > C_2$. Thus, it remains to analyze the difference in the first term of \eqref{eq:mono_ell_larger_mu1_equal_inte}, which has coefficient $\frac{\mu_2}{d(k+1,\ell-1)}$. Selecting an upper bound $v(i-1, k-1, \ell+1)$ within the minimum in that first term gives
    \begin{align*}
        \lefteqn{v(i,k-1,\ell+1) - \min\{v(i-1, k, \ell), v(i-1, k-1, \ell+1)\}}&\\
        &\quad \geq v(i,k-1,\ell+1) - v(i-1, k-1, \ell+1)\\
        &\quad \geq 0,
    \end{align*}where the last inequality follows from Lemma \ref{lemma:mono_i}. Consequently, the sum of the last three terms in \eqref{eq:mono_ell_larger_mu1_equal}, which is equal to \eqref{eq:mono_ell_larger_mu1_equal_inte},
    is non-negative, thus establishing the result for $\ell=C_2$.
    \item If $\ell < C_2$, it is still the case that $d(k+1,\ell-1) > d(k,\ell) > d(k-1,\ell+1)$.
    Thinning the MDP with a higher rate $d(k+1,\ell-1) = (k+1)\mu_1 + (\ell-1)\mu_2$ yields
        \begin{align}
            & v(i,k+1,\ell-1) - v(i,k,\ell) - \Big[v(i,k,\ell) - v(i,k-1,\ell+1)\Big] \nonumber \\
            &\quad =\frac{(k-1) \mu_1}{d(k+1,\ell-1)} \Big[\min\{v(i-1, k+1, \ell-1), v(i-1, k, \ell)\} \nonumber \\
            &\qquad \quad - \min\{v(i-1, k, \ell), v(i-1, k-1, \ell+1)\} \nonumber \\
            &\qquad \quad - \min\{v(i-1, k, \ell), v(i-1, k-1, \ell+1)\} \nonumber \\
            &\qquad \quad + \min\{v(i-1, k-1, \ell+1), v(i-1, k-2, \ell+2)\} \Big] \nonumber \\
            &\qquad + \frac{(\ell-1)\mu_2}{d(k+1,\ell-1)} \Big[\min\{v(i-1, k+2, \ell-2), v(i-1, k+1, \ell-1)\} \nonumber \\
            &\qquad \quad - \min\{v(i-1, k+1, \ell-1), v(i-1, k, \ell)\} \nonumber \\
            &\qquad \quad - \min\{v(i-1, k+1, \ell-1), v(i-1, k, \ell)\} \nonumber \\
            &\qquad \quad + \min\{v(i-1, k, \ell), v(i-1, k-1, \ell+1)\}\Big] \nonumber \\
            &\qquad + \frac{2\mu_2}{d(k+1,\ell-1)} \Big[\min\{v(i-1, k+1, \ell-1), v(i-1, k, \ell)\} \nonumber \\
            &\qquad \quad - \min\{v(i-1, k, \ell), v(i-1, k-1, \ell+1)\} \nonumber \\
            &\qquad \quad - \min\{v(i-1, k+1, \ell-1), v(i-1, k, \ell)\} \nonumber \\
            &\qquad \quad + \min\{v(i-1, k, \ell), v(i-1, k-1, \ell+1)\}\Big] \nonumber \\
            &\qquad + \frac{2(\mu_1-\mu_2)}{d(k+1,\ell-1)}\Big[\min\{v(i-1, k+1, \ell-1), v(i-1, k, \ell)\} \nonumber \\
            &\qquad \quad - \min\{v(i-1, k, \ell), v(i-1, k-1, \ell+1)\} \nonumber \\
            &\qquad \quad -v(i,k,\ell) + v(i,k-1,\ell+1)\Big].
            \label{eq:mono_ell_larger_mu1_smaller_ell}
        \end{align}
    Equation \ref{eq:mono_ell_larger_mu1_smaller_ell} is non-negative because the second-to-last term, with coefficient $\frac{2\mu_2}{d(k+1,\ell-1)}$, is zero, and the remaining terms are non-negative by an argument analogous to \eqref{eq:mono_ell_larger_mu1_larger_ell}.
    \end{enumerate}
    
    The results are now established for $\mu_1 \geq \mu_2$. We proceed with the proof of \eqref{eq:mono_ell} for $\mu_2 \geq \mu_1$ using induction on $i$. When $i=0$, the result follows directly from \eqref{eq:mono_ell_larger_mu1_base} again. Assuming it holds at $i-1$, we prove it for $i\geq1$ so that $x \in \tilde{\X}_{diff}$, considering three cases: $\ell > C_2$, $\ell = C_2$, and $\ell <C_2$.
   \begin{enumerate}[label= \textbf{Case} \arabic*:, leftmargin=3\parindent]
        \item If $\ell > C_2$, then $\ell-1 \geq C_2$. Consider $D(i,k+1,\ell-1) - D(i,k,\ell)$. After comparable arithmetic that led to \eqref{eq:diff} and \eqref{eq:diff_const_expand}, 
        then grouping the terms that determine where to work after a Station 1 (and 2) service completion yields 
        \begin{align}
        & v(i,k+1,\ell-1) - v(i,k,\ell) - \Big[v(i,k,\ell) - v(i,k-1,\ell+1)\Big] \nonumber \\
        &\quad = \left(\frac{C_2\mu_1\mu_2}{d(k+1,\ell-1)d(k,\ell)}\left(\frac{-ih_0}{C_2 \mu_2} + \frac{h_1-h_2}{\mu_1} - \frac{C_1h_2}{C_2 \mu_2} \right)\right) \nonumber \\
        &\qquad \quad -\left(\frac{C_2\mu_1\mu_2}{d(k,\ell)d(k-1,\ell+1)}\left(\frac{-ih_0}{C_2 \mu_2} + \frac{h_1-h_2}{\mu_1} - \frac{C_1h_2}{C_2 \mu_2} \right) \right) \nonumber \\
        &\qquad +\frac{(k-1) \mu_1}{d(k-1,\ell+1)} \Big[\min\{v(i-1, k+1, \ell-1), v(i-1, k, \ell)\} \nonumber \\
        &\qquad \quad - \min\{v(i-1, k, \ell), v(i-1, k-1, \ell+1)\} \nonumber \\
        &\qquad \quad - \min\{v(i-1, k, \ell), v(i-1, k-1, \ell+1)\} \nonumber \\
        &\qquad \quad + \min\{v(i-1, k-1, \ell+1), v(i-1, k-2, \ell+2)\} \Big] \nonumber \\
        &\qquad + \frac{C_2\mu_2}{d(k+1,\ell-1)} \Big[\min\{v(i-1, k+2, \ell-2), v(i-1, k+1, \ell-1)\} \nonumber \\
        &\qquad \quad - \min\{v(i-1, k+1, \ell-1), v(i-1, k, \ell)\}\Big] \nonumber \\
        & \qquad \quad - \min\{v(i-1, k+1, \ell-1), v(i-1, k, \ell)\} \nonumber \\
        &\qquad \quad + \min\{v(i-1, k, \ell), v(i-1, k-1, \ell+1)\}\Big] \nonumber \\
        &\qquad \quad + \frac{2C_2\mu_1^2\mu_2}{d(k+1,\ell-1)d(k,\ell)d(k-1,\ell+1)}
        \Big[ \nonumber \\
        &\qquad \quad \min\{v(i-1, k+1, \ell-1), v(i-1, k, \ell)\} \nonumber \\
        &\qquad \quad - \min\{v(i-1, k, \ell), v(i-1, k-1, \ell+1)\}\Big] 
        \label{eq:mono_ell_larger_mu2_larger_ell}
    \end{align}
    Following analogously to the proof of the non-negativity of the term with coefficient $\frac{(k-1) \mu_1}{d(k+1,\ell-1)}$ in \eqref{eq:mono_ell_larger_mu1_larger_ell} for $\ell > C_2$ where $\mu_1 \geq \mu_2$, both terms in \eqref{eq:mono_ell_larger_mu2_larger_ell} with coefficients $\frac{(k-1) \mu_1}{d(k-1,\ell+1)}$ and $\frac{C_2\mu_2}{d(k+1,\ell-1)}$, are non-negative. Consider the sum of the remaining terms in \eqref{eq:mono_ell_larger_mu2_larger_ell}.
    \begin{align}
        \lefteqn{\left(\frac{C_2\mu_1\mu_2}{d(k+1,\ell-1)d(k,\ell)}\left(\frac{-ih_0}{C_2 \mu_2} + \frac{h_1-h_2}{\mu_1} - \frac{C_1h_2}{C_2 \mu_2} \right)\right)}& \nonumber \\
        &\qquad \quad -\left(\frac{C_2\mu_1\mu_2}{d(k,\ell)d(k-1,\ell+1)}\left(\frac{-ih_0}{C_2 \mu_2} + \frac{h_1-h_2}{\mu_1} - \frac{C_1h_2}{C_2 \mu_2} \right) \right) \nonumber \\
        &\qquad + \frac{2C_2\mu_1^2\mu_2}{d(k+1,\ell-1)d(k,\ell)d(k-1,\ell+1)}
        \Big[ \nonumber \\
        &\qquad \quad \min\{v(i-1, k+1, \ell-1), v(i-1, k, \ell)\} \nonumber \\
        &\qquad \quad - \min\{v(i-1, k, \ell), v(i-1, k-1, \ell+1)\}\Big] \nonumber \\
        &\quad = \frac{2C_2\mu_1^2\mu_2}{d(k+1,\ell-1)d(k,\ell)d(k-1,\ell+1)}
        \Big[ \frac{ih_0}{C_2\mu_2} - \left( \frac{h_1-h_2}{\mu_1}- \frac{C_1h_2}{C_2\mu_2}\right)\nonumber \\
        &\qquad \quad + \min\{v(i-1, k+1, \ell-1), v(i-1, k, \ell)\} \nonumber \\
        &\qquad \quad - \min\{v(i-1, k, \ell), v(i-1, k-1, \ell+1)\}\Big].
        \label{eq:mono_ell_larger_mu2_larger_ell_inte}
    \end{align}
     For the difference $\min\{v(i-1, k+1, \ell-1), v(i-1, k, \ell)\} - \min\{v(i-1, k, \ell), v(i-1, k-1, \ell+1)\}$ in \eqref{eq:mono_ell_larger_mu2_larger_ell_inte}, we have
    \begin{align*}
        \lefteqn{\min\{v(i-1, k+1, \ell-1), v(i-1, k, \ell)\} - \min\{v(i-1, k, \ell), v(i-1, k-1, \ell+1)\}}&\\
        &\quad \geq \min\Big\{v(i-1, k+1, \ell-1) - v(i-1, k, \ell), v(i-1, k, \ell) - v(i-1, k-1, \ell+1)\Big\}\\
        &\quad \geq  \min\left\{-\frac{(i-1)h_0}{C_2\mu_2} + \left( \frac{h_1-h_2}{\mu_1}- \frac{C_1h_2}{C_2\mu_2}\right), -\frac{(i-1)h_0}{C_2\mu_2} + \left( \frac{h_1-h_2}{\mu_1}- \frac{C_1h_2}{C_2\mu_2}\right)\right\}\\
        &\quad \geq -\frac{ih_0}{C_2\mu_2} + \left( \frac{h_1-h_2}{\mu_1}- \frac{C_1h_2}{C_2\mu_2}\right),
    \end{align*}where the second inequality follows from Statement \ref{state:diff_lb} in Lemma \ref{silemma:diff}. Therefore, \eqref{eq:mono_ell_larger_mu2_larger_ell_inte} is non-negative, and so is \eqref{eq:mono_ell_larger_mu2_larger_ell}.
        \item If $\ell = C_2$, thinning the MDPs with a higher rate $d(k,\ell) = k\mu_1+ C_2\mu_2$ in both the differences $v(i,k+1,\ell-1) - v(i,k,\ell)$ and $v(i,k,\ell) - v(i,k-1,\ell+1)$ respectively yields (recall $\ell = C_2$):
    \begin{align*}
        \lefteqn{v(i,k+1,\ell-1) - v(i,k,\ell)}&\\
        &\quad = \frac{h_1-h_2}{d(k,\ell)} + \frac{k \mu_1}{d(k,\ell)} \Big[\min\{v(i-1, k+1, \ell-1), v(i-1, k, \ell)\}\\
        &\qquad \quad - \min\{v(i-1, k, \ell), v(i-1, k-1, \ell+1)\} \Big]\\
        &\qquad + \frac{(C_2-1)\mu_2}{d(k+1,\ell-1)} \Big[\min\{v(i-1, k+2, \ell-2), v(i-1, k+1, \ell-1)\}\\
        &\qquad \quad - \min\{v(i-1, k+1, \ell-1), v(i-1, k, \ell)\}\Big]\\
        &\qquad + \frac{\mu_1}{d(k,\ell)}\Big[\min\{v(i-1,k+1,\ell-1),v(i-1,k,\ell)\}\\
        &\qquad \quad - \min\{v(i-1,k+1,\ell-1),v(i-1,k,\ell)\}\Big]\\
        &\qquad + \frac{\mu_2-\mu_1}{d(k,\ell)}\Big[v(i,k+1,\ell-1) - \min\{v(i-1,k+1,\ell-1),v(i-1,k,\ell)\}\Big],
    \end{align*}and
    \begin{align*}
        \lefteqn{v(i,k,\ell) - v(i,k-1,\ell+1)}&\\
        &\quad = \frac{h_1-h_2}{d(k,\ell)} + \frac{(k-1)\mu_1}{d(k,\ell)}\Big[\min\{v(i-1,k,\ell),v(i-1,k-1,\ell+1)\}\\
        &\qquad \quad - \min\{v(i-1,k-1,\ell+1),v(i-1,k-2,\ell+2)\}\Big]\\
        &\qquad + \frac{C_2\mu_2}{d(k,\ell)}\Big[\min\{v(i-1,k+1,\ell-1),v(i-1,k,\ell)\}\\
        &\qquad \quad - \min\{v(i-1,k,\ell),v(i-1,k-1,\ell+1)\}\Big]\\
        &\qquad + \frac{\mu_1}{d(k,\ell)} \Big[\min\{v(i-1,k,\ell),v(i-1,k-1,\ell+1)\} - v(i,k-1,\ell+1)\Big].
    \end{align*}Hence \eqref{eq:mono_ell} equals to
        \begin{align}
        &v(i,k+1,\ell-1) - v(i,k,\ell) - \Big[v(i,k,\ell) - v(i,k-1,\ell+1)\Big] \nonumber \\
        &\quad = \frac{(k-1)\mu_1}{d(k,\ell)}\Big[\min\{v(i-1,k+1,\ell-1),v(i-1,k,\ell)\} \nonumber \\
        &\qquad \quad - \min\{v(i-1,k,\ell),v(i-1,k-1,\ell+1)\} \nonumber \\
        &\qquad \quad -\min\{v(i-1,k,\ell),v(i-1,k-1,\ell+1)\} \nonumber \\
        &\qquad \quad +\min\{v(i-1,k-1,\ell+1),v(i-1,k-2,\ell+2)\}\Big] \nonumber \\
        &\qquad + \frac{(C_2-1)\mu_2}{d(k+1,\ell-1)}\Big[\min\{v(i-1,k+2,\ell-2),v(i-1,k+1,\ell-1)\} \nonumber \\
        &\qquad \quad -\min\{v(i-1,k+1,\ell-1),v(i-1,k,\ell)\} \nonumber \\
        &\qquad \quad - \min\{v(i-1,k+1,\ell-1),v(i-1,k,\ell)\} \nonumber \\
        &\qquad \quad + \min\{v(i-1,k,\ell),v(i-1,k-1,\ell+1)\}\Big] \nonumber \\
        &\qquad + \frac{\mu_2-\mu_1}{d(k,\ell)}\Big[v(i,k+1,\ell-1) - \min\{v(i-1,k+1,\ell-1),v(i-1,k,\ell)\} \nonumber \\
        &\qquad \quad - \min\{v(i-1,k+1,\ell-1)-v(i-1,k,\ell)\} \nonumber \\
        &\qquad \quad + \min\{v(i-1,k,\ell),v(i-1,k-1,\ell+1)\}\Big] \label{eq:mono_ell_larger_mu2_equal} \\
        &\qquad - \frac{\mu_1}{d(k,\ell)}\Big[\min\{v(i-1,k,\ell),v(i-1,k-1,\ell+1)\} - v(i,k-1,\ell+1)\Big] \nonumber.
    \end{align}
    The non-negativity of the first two terms in \eqref{eq:mono_ell_larger_mu2_equal} with coefficients $\frac{(k-1)\mu_1}{d(k,\ell)}$ and $\frac{(C_2-1)\mu_2}{d(k+1,\ell-1)}$ follows similarly as that of the term with coefficient $\frac{(k-1)\mu_1}{d(k+1,\ell-1)}$ in \eqref{eq:mono_ell_larger_mu1_larger_ell} for $\ell > C_2$ and $\mu_1 \geq \mu_2$. Consider the third term in \eqref{eq:mono_ell_larger_mu2_equal} with coefficient $\frac{\mu_2-\mu_1}{d(k,\ell)}$ in two cases based on whether $v(i-1,k,\ell) \leq v(i-1,k-1,\ell+1)$ or not.
    \begin{enumerate}[label= \textbf{Case} \alph*:, leftmargin=3\parindent]
        \item If $v(i-1,k,\ell) \leq v(i-1,k-1,\ell+1)$, choosing upper bounds $v(i-1,k+1,\ell-1)$ and $v(i-1,k,\ell)$ in the second and the third minima, respectively, yields
        \begin{align*}
        \lefteqn{v(i,k+1,\ell-1) - \min\{v(i-1,k+1,\ell-1),v(i-1,k,\ell)\}}& \nonumber \\
        &\qquad - \min\{v(i-1,k+1,\ell-1)-v(i-1,k,\ell)\} \nonumber \\
        &\qquad + \min\{v(i-1,k,\ell),v(i-1,k-1,\ell+1)\}\\
        &\quad \geq v(i,k+1,\ell-1) - v(i-1,k+1,\ell-1) \nonumber \\
        &\qquad - v(i-1,k,\ell) + v(i-1,k,\ell)\\
        &\quad = v(i,k+1,\ell-1) - v(i-1,k+1,\ell-1) \nonumber \\
        &\quad \geq 0,
    \end{align*}where the last inequality follows from Lemma \ref{lemma:mono_i}.
    \item If $v(i-1,k,\ell) > v(i-1,k-1,\ell+1)$, choosing the same upper bound $v(i-1,k,\ell)$ in both the second and the third minima yields
        \begin{align}
        \lefteqn{v(i,k+1,\ell-1) - \min\{v(i-1,k+1,\ell-1),v(i-1,k,\ell)\}}& \nonumber \\
        &\qquad - \min\{v(i-1,k+1,\ell-1)-v(i-1,k,\ell)\} \nonumber \\
        &\qquad + \min\{v(i-1,k,\ell),v(i-1,k-1,\ell+1)\} \nonumber \\
        &\quad \geq v(i,k+1,\ell-1) - v(i-1,k,\ell)& \nonumber \\
        &\qquad - v(i-1,k,\ell) + v(i-1,k-1,\ell+1) \nonumber \\
        &\quad = \Big[v(i,k+1,\ell-1) - v(i-1,k+1,\ell-1)\Big] \nonumber \\
        &\qquad + \Big[v(i-1,k+1,\ell-1)- v(i-1,k,\ell)& \nonumber \\
        &\qquad - v(i-1,k,\ell) + v(i-1,k-1,\ell+1)\Big] \nonumber \\
        &\quad \geq 0, \label{eq:mono_ell_larger_mu2_equal_inte1}
    \end{align}where the last step holds by Lemma \ref{lemma:mono_i} and the inductive hypothesis.
    \end{enumerate}
    It remains to analyze the last term in \eqref{eq:mono_ell_larger_mu2_equal} with coefficient $- \frac{\mu_1}{d(k,\ell)}$. Selecting an upper bound $v(i-1,k-1,\ell+1)$ within the first minimum gives
    \begin{align}
        \lefteqn{\min\{v(i-1,k,\ell),v(i-1,k-1,\ell+1)\} - v(i,k-1,\ell+1)}& \nonumber \\
        &\quad \leq v(i-1,k-1,\ell+1) - v(i,k-1,\ell+1) \nonumber \\
        &\quad \leq 0,
    \end{align}
    where the last step applies Lemma \ref{lemma:mono_i} again, showing that the final term in \eqref{eq:mono_ell_larger_mu2_equal} is also non-negative. Together, the non-negativity of all four terms in \eqref{eq:mono_ell_larger_mu2_equal} ensures that the entire expression is non-negative.
    \item If $\ell < C_2$, again consider the expansion of the difference as expressed in \eqref{eq:diff} applied at $(i,k+1,\ell-1)$ and $(i,k,\ell)$, respectively, followed by the substitution of \eqref{eq:diff_const_expand}. By grouping the terms representing where to work after a Station 1 (and 2) service completion, we have 
    \begin{align}
        & v(i,k+1,\ell-1) - v(i,k,\ell) - \Big[v(i,k,\ell) - v(i,k-1,\ell+1)\Big] \nonumber \\
        &\quad = \frac{C_1\mu_1\mu_2}{d(k+1,\ell-1)d(k,\ell)} \left(\frac{ih_0}{C_1}\left(\frac{1}{\mu_1} - \frac{1}{\mu_2}\right) + \frac{h_1}{\mu_1} - \frac{h_2}{\mu_2} \right) \nonumber \\
        &\qquad - \frac{C_1\mu_1\mu_2}{d(k,\ell)d(k-1,\ell+1)} \left(\frac{ih_0}{C_1}\left(\frac{1}{\mu_1} - \frac{1}{\mu_2}\right) + \frac{h_1}{\mu_1} - \frac{h_2}{\mu_2} \right) \nonumber \\
        &\qquad +\frac{(k-1) \mu_1}{d(k-1,\ell+1)} \Big[\min\{v(i-1, k+1, \ell-1), v(i-1, k, \ell)\} \nonumber \\
        &\qquad \quad - \min\{v(i-1, k, \ell), v(i-1, k-1, \ell+1)\} \nonumber \\
        &\qquad \quad - \min\{v(i-1, k, \ell), v(i-1, k-1, \ell+1)\} \nonumber \\
        &\qquad \quad + \min\{v(i-1, k-1, \ell+1), v(i-1, k-2, \ell+2)\} \Big] \nonumber \\
        &\qquad + \frac{(\ell-1)\mu_2}{d(k+1,\ell-1)} \Big[\min\{v(i-1, k+2, \ell-2), v(i-1, k+1, \ell-1)\} \nonumber \\
        &\qquad \quad - \min\{v(i-1, k+1, \ell-1), v(i-1, k, \ell)\}\Big] \nonumber \\
        & \qquad \quad - \min\{v(i-1, k+1, \ell-1), v(i-1, k, \ell)\} \nonumber \\
        &\qquad \quad + \min\{v(i-1, k, \ell), v(i-1, k-1, \ell+1)\}\Big] \nonumber \\
        &\qquad \quad - \frac{2C_1\mu_1\mu_2(\mu_2-\mu_1)}{d(k+1,\ell-1)d(k,\ell)d(k-1,\ell+1)}
        \Big[ \nonumber \\
        &\qquad \qquad \min\{v(i-1, k+1, \ell-1), v(i-1, k, \ell)\} \nonumber \\
        &\qquad \quad - \min\{v(i-1, k, \ell), v(i-1, k-1, \ell+1)\}\Big]. \label{eq:mono_ell_larger_mu2_smaller_ell}
    \end{align}
    Notice that the terms in \eqref{eq:mono_ell_larger_mu2_smaller_ell} with coefficients $\frac{(k-1) \mu_1}{d(k-1,\ell+1)}$ and $\frac{(\ell-1)\mu_2}{d(k+1,\ell-1)}$ are non-negative by a similar reasoning as that of the term with coefficient $\frac{(k-1)\mu_1}{d(k+1,\ell-1)}$ in \eqref{eq:mono_ell_larger_mu1_larger_ell} for $\ell > C_2$ and $\mu_1 \geq \mu_2$.
    Consider the sum of the remaining terms in \eqref{eq:mono_ell_larger_mu2_smaller_ell}.
    \begin{align}
        \lefteqn{\frac{C_1\mu_1\mu_2}{d(k+1,\ell-1)d(k,\ell)} \left(\frac{ih_0}{C_1}\left(\frac{1}{\mu_1} - \frac{1}{\mu_2}\right) + \frac{h_1}{\mu_1} - \frac{h_2}{\mu_2} \right)}& \nonumber \\
        &\qquad - \frac{C_1\mu_1\mu_2}{d(k,\ell)d(k-1,\ell+1)} \left(\frac{ih_0}{C_1}\left(\frac{1}{\mu_1} - \frac{1}{\mu_2}\right) + \frac{h_1}{\mu_1} - \frac{h_2}{\mu_2} \right) \nonumber \\
        &\qquad - \frac{2C_1\mu_1\mu_2(\mu_2-\mu_1)}{d(k+1,\ell-1)d(k,\ell)d(k-1,\ell+1)}
        \Big[ \nonumber \\
        &\qquad \quad \min\{v(i-1, k+1, \ell-1), v(i-1, k, \ell)\} \nonumber \\
        &\qquad - \min\{v(i-1, k, \ell), v(i-1, k-1, \ell+1)\}\Big] \nonumber\\
        & \quad = \frac{2C_1\mu_1\mu_2(\mu_2-\mu_1)}{d(k+1,\ell-1)d(k,\ell)d(k-1,\ell+1)} \Big[\left(\frac{ih_0}{C_1}\left(\frac{1}{\mu_1} - \frac{1}{\mu_2}\right) + \frac{h_1}{\mu_1} - \frac{h_2}{\mu_2} \right) \nonumber \\
        &\qquad -\Big(\min\{v(i-1, k+1, \ell-1), v(i-1, k, \ell)\} \nonumber \\
        &\qquad - \min\{v(i-1, k, \ell), v(i-1, k-1, \ell+1)\}\Big)\Big].\label{eq:diff_ub_larger_mu2_inte}
    \end{align}
    Note, for the difference $\min\{v(i-1, k+1, \ell-1), v(i-1, k, \ell)\} - \min\{v(i-1, k, \ell), v(i-1, k-1, \ell+1)\}$ in \eqref{eq:diff_ub_larger_mu2_inte} we have 
    \begin{align*}
        \lefteqn{\min\{v(i-1, k+1, \ell-1), v(i-1, k, \ell)\} - \min\{v(i-1, k, \ell), v(i-1, k-1, \ell+1)\}}&\\
        &\quad \leq \max\Big\{v(i-1, k+1, \ell-1) - v(i-1, k, \ell), v(i-1, k, \ell) - v(i-1, k-1, \ell+1)\Big\}\\
        &\quad \leq \max\left\{\frac{(i-1)h_0}{C_1}\left(\frac{1}{\mu_1} - \frac{1}{\mu_2}\right) + \frac{h_1}{\mu_1} - \frac{h_2}{\mu_2}, \frac{(i-1)h_0}{C_1}\left(\frac{1}{\mu_1} - \frac{1}{\mu_2}\right) + \frac{h_1}{\mu_1} - \frac{h_2}{\mu_2} \right\}\\
        &\quad \leq \frac{ih_0}{C_1}\left(\frac{1}{\mu_1} - \frac{1}{\mu_2}\right) + \frac{h_1}{\mu_1} - \frac{h_2}{\mu_2},
    \end{align*}
    where the second inequality follows from Statement \ref{state:diff_ub_larger_mu2} in Lemma \ref{silemma:diff}, given that $\mu_1 \geq \mu_1$ and $\ell - 1 < \ell < C_2$, and the last inequality holds by $\mu_2 \geq \mu_1$. Moreover, the coefficient of \eqref{eq:diff_ub_larger_mu2_inte} is non-negative by another application of $\mu_2 \geq \mu_1$, ensuring that \eqref{eq:diff_ub_larger_mu2_inte} is non-negative, and consequently, so is \eqref{eq:mono_ell_larger_mu2_smaller_ell}.
\end{enumerate}
To conclude, inequality \eqref{eq:mono_ell} holds for both cases, $\mu_1 \geq \mu_2$ and $\mu_2\geq\mu_1$.
\end{proof}

This subsection presents the proofs for the supporting results for the case where $\frac{h_1}{\mu_1} \leq \frac{h_2}{\mu_2}$, as outlined in Subsection \ref{sec:support-results}.

\begin{proof} [Proof of Corollary \ref{cor:enough_Type_II}.]
    Recall that $x \in \X_{diff}$ implies that $k = C_1 - \ell \geq 1$ and $x \in \X$. Since $C_2 \geq C_1$, it follows that $\ell \leq C_1-k \leq C_1 - 1 \leq C_2 -1$, and therefore $\ell < C_2$ necessarily holds. Consequently, the Bellman equations for the case where $C_2 \geq C_1$ coincide with those for $C_2 < C_1$ with $\ell < C_2$. A similar argument applies to proving the corresponding statements.
\end{proof}

\subsection{Preliminaries for heuristic design} \label{sec:prelim-heuristics}
In addition to Statement \ref{state:recursive_pos} in Lemma \ref{lemma:diff} and the entire Lemma \ref{lemma:heur} outlined in Section \ref{sec:prelim}, we introduce Lemma \ref{silemma:heur} as another essential preliminary result for heuristic analysis, with all proofs provided in this subsection.

Recall from Definitions \ref{def:probs_k}, a little arithmetic yields the following identity:
\begin{align}
    p_k+q_k+r_k = 1. \label{eq:identity}
\end{align}
Additionally, Definition \ref{def:probs_quantities_k} directly leads to
\begin{align}
    b' \leq b_k \leq b, \quad  c' \leq c_k \leq c. \label{eq:cmp_bk_b_ck_c}
\end{align}

As needed, we reindex these quantities from $k$ to $\ell$ to facilitate analysis:

\begin{definition} \label{def:probs_quantities_ell}
    For $\ell = 0, 1, \ldots, C_1-1$ and $k+\ell = C_1$, let
    \begin{align*}
        \tilde{p}_\ell:=p_k, \quad \tilde{q}_\ell:=q_k, \quad \tilde{r}_\ell:=r_k, \quad \tilde{c}_\ell:=c_k, \quad \tilde{b}_\ell:=b_k,
    \end{align*}where $p_k, q_k, r_k, c_k$ and $b_k$ are defined in Definition \ref{def:probs_quantities_k}. 
\end{definition}

\noindent Lemma \ref{silemma:heur} below, which complements Lemma \ref{lemma:heur}, bridges the results concerning the comparison of $H(i,k,\ell)$ and $D(i,k,\ell)$.
Here, Statement \ref{state:recursive_neg} of Lemma \ref{silemma:heur} reformulates the expression for $D(i,k,\ell)$ when $D(i-1,k,\ell) \leq 0$, using the notations in Definition \ref{def:probs_quantities_ell}. 
Statements \ref{state:heur_mono_y} and \ref{state:y_fraction} in Lemma \ref{silemma:heur} provide additional properties of $y_k$, where $y_k$ is defined in \eqref{def:sequence_y}, and Statements \ref{state:heur_mono_z} and \ref{state:z_fraction} in the same lemma provide further properties of $z_\ell$ with $z_\ell$ specified in \eqref{def:sequence_z}.

\begin{lemma}\label{silemma:heur}
    The following results hold.
    \begin{enumerate}
        \item \label{state:recursive_neg}
        Consider $x =(i,k,\ell)\in \tilde{\X}_{diff}$, where $i \geq 1$. Then $D(i-1,k,\ell) \leq 0$ implies
        \begin{align}
            D(i,k,\ell) = \tilde{p}_\ell(i\tilde{c}_\ell+\tilde{b}_\ell) + \tilde{q}_\ell D(i-1,k,\ell) + \tilde{r}_\ell \min\{D(i-1,k+1,\ell-1), 0\}, \label{eq:recursive_neg}
        \end{align}where $\tilde{p}_\ell, \tilde{q}_\ell, \tilde{r}_\ell, \tilde{c}_\ell$ and $\tilde{b}_\ell$ are specified in Definition \eqref{def:probs_quantities_ell}.
        \item Fix any $k = 1, 2, \ldots, C_1$. The following results concerning $y_k$ (defined in \eqref{def:sequence_y}) hold:
        \begin{enumerate} [ref=\theenumi(\alph*)]
            \item \label{state:heur_mono_y}
            For $k=1, 2, \ldots, C_1-1$,
            \begin{align*}
                (y_k+1) - y_{k+1}&=\left\{
                \begin{array}{lcl}
                    m &  \text{if } &\ell-1 =C_1-(k+1)< C_2,\\
                    0 &  \text{if } &\ell-1 =C_1-(k+1)\geq C_2.
                \end{array}
                \right.
            \end{align*} 
            Consequently, $(y_k+1) - y_{k+1} \geq 0$ by the definition of $m$ in \eqref{def:ratio_of_rates}.
            \item \label{state:y_fraction}
            $y_k \geq \frac{r_k}{p_k}$. In particular, if $\ell = C_1-k \geq C_2$, $y_k = \frac{r_k}{p_k}$ holds with equality, where $r_k$ and $p_k$ are provided in \eqref{def:rk} and \eqref{def:pk}, respectively.
        \end{enumerate}    
        \item Fix any $\ell = 0, 1, \ldots, C_2-1$. The following results hold for $z_\ell$ (defined in \eqref{def:sequence_z}):
        \begin{enumerate} [ref=\theenumi(\alph*)]
            \item \label{state:heur_mono_z}
            $z_{\ell}+1 - z_{\ell+1} > 0$, where $\ell=0, 1, \ldots, C_2-2$ (if one exists). 
            \item \label{state:z_fraction}
            $z_{\ell} \geq \frac{\tilde{q}_\ell}{\tilde{p}_\ell}$, where $\tilde{q}_\ell$ and $\tilde{p}_\ell$ are given in Definition \ref{def:probs_quantities_ell}. 
        \end{enumerate}
    \end{enumerate}
\end{lemma}

\noindent We begin by proving properties of $D(i,k,\ell)$ as established in Statement \ref{state:recursive_pos} of Lemma \ref{lemma:diff} and Statement \ref{state:recursive_neg} of Lemma \ref{silemma:heur}.

\begin{proof} [Proof of Statement \ref{state:recursive_pos} in Lemma \ref{lemma:diff}]
    Observe that for $i \geq 1$, if $D(i-1,k,\ell) = v(i-1, k, \ell) - v(i-1, k-1, \ell+1) \geq 0$, the difference in the second term with coefficient $q_k$ in \eqref{eq:alge_diff} becomes
    \begin{align}
        \lefteqn{\min\{v(i-1, k, \ell), v(i-1, k-1, \ell+1)\}}& \nonumber \\
        &\qquad - \min\{v(i-1, k-1, \ell+1), v(i-1, k-2, \ell+2)\} \nonumber \\
        &\quad = v(i-1, k-1, \ell+1)\} - \min\{v(i-1, k-1, \ell+1), v(i-1, k-2, \ell+2)\} \nonumber \\
        &\quad = \max\{0,D(i-1,k-1,\ell+1)\}. \label{eq:recursive_qk_terms}
    \end{align}
    Consider now the difference in the third term with coefficient $r_k$ in \eqref{eq:alge_diff}. Proposition \ref{prop:mono_ell} ensures that if $D(i-1,k,\ell) \geq 0$, $v(i-1, k+1, \ell-1) \geq v(i-1, k, \ell)$. Therefore,
        \begin{align}
            \lefteqn{\min\{v(i-1, k+1, \ell-1), v(i-1, k, \ell)\}}& \nonumber \\
            &\qquad - \min\{v(i-1, k, \ell), v(i-1, k-1, \ell+1)\} \nonumber \\
            &\quad = v(i-1, k, \ell) - v(i-1, k-1, \ell+1) \nonumber \\
            &\quad =D(i-1,k,\ell). \label{eq:recursive_rk_terms}
        \end{align}
    Under the assumption $D(i-1,k,\ell) \geq 0$, applying both bounds in \eqref{eq:recursive_qk_terms} and \eqref{eq:recursive_rk_terms} in \eqref{eq:alge_diff} yields the recursive equation \eqref{eq:recursive_pos} in Statement \ref{state:recursive_pos} in Lemma \ref{lemma:diff}.
\end{proof}

\begin{proof} [Proof of Statement \ref{state:recursive_neg} in Lemma \ref{silemma:heur}.]
    Using the notations from Definition \ref{def:probs_quantities_ell} in the equation in \eqref{eq:diff} (together with \eqref{eq:diff_const_expand}), yields that for $x \in \tilde{\X}_{diff}$ and $i \geq 1$,
        \begin{align}
        \lefteqn{D(i,k,\ell) = v(i,k,\ell) - v(i,k-1,\ell+1)}& \nonumber \\
        &\quad = \tilde{p}_\ell(i\tilde{c}_\ell+\tilde{b}_\ell) \nonumber \\
        &\qquad + \tilde{q}_\ell \Big[\min\{v(i-1, k, \ell), v(i-1, k-1, \ell+1)\} \nonumber \\
        &\qquad \quad - \min\{v(i-1, k-1, \ell+1), v(i-1, k-2, \ell+2)\} \Big] \nonumber \\
        &\qquad + \tilde{r}_\ell \Big[\min\{v(i-1, k+1, \ell-1), v(i-1, k, \ell)\} \nonumber \\
        &\qquad \quad - \min\{v(i-1, k, \ell), v(i-1, k-1, \ell+1)\}\Big]. \label{eq:alge_diff_ell} 
    \end{align}
    Notice that by Proposition \ref{prop:mono_ell}, if $i \geq 1$, $D(i-1,k,\ell) \leq 0$ implies $D(i-1,k-1,\ell+1)\leq 0$. The difference in the second term in \eqref{eq:alge_diff_ell} with coefficient $\tilde{q}_\ell$ becomes
    \begin{align}
        \lefteqn{\min\{v(i-1, k, \ell), v(i-1, k-1, \ell+1)\}}& \nonumber \\
        &\qquad - \min\{v(i-1, k-1, \ell+1), v(i-1, k-2, \ell+2)\} \nonumber \\
        &\quad = v(i-1, k, \ell)\} - v(i-1, k-1, \ell+1) \nonumber \\
        &\quad = D(i-1,k,\ell). \label{eq:alge_diff_q_ell_terms}
    \end{align}
And the difference in the third term with coefficient $\tilde{r}_\ell$ in \eqref{eq:alge_diff_ell} can be expressed as follows:
    \begin{align}
        \lefteqn{\min\{v(i-1, k+1, \ell-1), v(i-1, k, \ell)\}}& \nonumber \\
        &\qquad - \min\{v(i-1, k, \ell), v(i-1, k-1, \ell+1)\} \nonumber \\
        &\quad = \min\{v(i-1, k+1, \ell-1), v(i-1, k, \ell)\} - v(i-1, k, \ell) \nonumber \\
        &\quad = \min\{D(i-1,k+1,\ell-1), 0\}. \label{eq:alge_diff_r_ell_terms}
    \end{align}
    Applying both bounds in \eqref{eq:alge_diff_q_ell_terms} and \eqref{eq:alge_diff_r_ell_terms} in \eqref{eq:alge_diff_ell} gives the recursive equation \eqref{eq:recursive_neg} in Statement \ref{state:recursive_neg} in Lemma \ref{silemma:heur}.    
\end{proof} 

\noindent Next, we verify the previously stated assertions on the properties of $y_k$ and $z_\ell$. For $y_k$, see Statements \ref{state:heur_mono_y} and \ref{state:y_fraction} in Lemma \ref{silemma:heur} and Statement \ref{state:y_ub} in Lemma \ref{lemma:heur}.
For $z_\ell$, refer to Statements \ref{state:heur_mono_z} and \ref{state:z_fraction} in Lemma \ref{silemma:heur} and Statement \ref{state:z_ub} in Lemma \ref{lemma:heur}.

\begin{proof} [Proof of Statement \ref{state:heur_mono_y} in Lemma \ref{silemma:heur}]
    For any fixed $k = 1, 2, \ldots, C_1-2$ (with $\ell = C_1 - k$), directly applying the expression of $y_k$ from \eqref{def:sequence_y} yield
    \begin{align*}
        \lefteqn{(y_k+1) - y_{k+1}}\\
        &\quad = \big((k-1)+\min\{\ell,C_2\}m + 1\big) - \big(k+\min\{\ell-1,C_2\}m \big)\\
        &\quad = \Big(\min\{\ell,C_2\} - \min\{\ell-1,C_2\}\Big)m.
    \end{align*}
    The result follows by discussing whether $\ell-1 \geq C_2$.
\end{proof}

\begin{proof} [Proof of Statement \ref{state:y_fraction} in Lemma \ref{silemma:heur}.]
    Recall the expressions of $r_k$ and $p_k$ from \eqref{def:rk} and \eqref{def:pk} in Definition \ref{def:probs_k}, respectively. For any $k=1,2, \ldots,C_1$, we have
    \begin{align*}
        & \frac{r_k}{p_k}= \left\{
        \begin{array}{lcl}
             \frac{\ell}{C_1}\big((k-1)+(\ell+1)m\big) & \text{if} & \ell = C_1 - k < C_2,\\
              (k-1)+C_2 m & \text{if} & \ell = C_1 - k \geq C2_,
        \end{array}
        \right.
    \end{align*}
    where $m =\frac{\mu_2}{\mu_1} >0.$ 
    Recall the expression of $y_k$ from \eqref{def:sequence_y} and the results are established by considering the following two cases based on whether $\ell < C_2$.
    \begin{enumerate}[label= \textbf{Case} \arabic*:, leftmargin=3\parindent]
        \item If $\ell < C_2$, $\frac{r_k}{p_k} \leq k-1+\ell m = y_k$ since $\ell < \ell+1 \leq C_1$ and $m>0$.
        \item If $\ell \geq C_2$, $\frac{r_k}{p_k} = k-1+C_2 m = y_k$.
    \end{enumerate}
\end{proof}

\begin{proof} [Proof of Statement \ref{state:y_ub} in Lemma \ref{lemma:heur}]
    Referring to the expression of $y_k$ from \eqref{def:sequence_y} results in
    \begin{align*}
        y_k =(k-1)+\min\{\ell,C_2\} m \leq (k-1)+\ell m \leq (k-1)+\ell, 
    \end{align*}where the first inequality replaces the minimum with an upper bound (recall $m= \frac{\mu_2}{\mu_1}>0$), and the second inequality applies $m  \leq 1$ since $\mu_1 \geq \mu_2$. The result thus follows by noticing $k+\ell = C_1$.
\end{proof}

\begin{proof} [Proof of Statement \ref{state:heur_mono_z} in Lemma \ref{silemma:heur}.]
    For any $\ell = 0, 1, \ldots, C_2-2$ (if one exists, i.e., $C_2 \geq 2$) with $k = C_1 - \ell$, directly applying the definition of $z_\ell$ in \ref{def:sequence_z} gives
    \begin{align*}
        \lefteqn{(z_{\ell}+1) - z_{\ell+1}}\\
        &\quad = \left(\frac{k-1}{m} + \ell + 1\right) - \left(\frac{k-2}{m} + (\ell + 1)\right)\\
        &\quad = \frac{k-1}{m} - \frac{k-2}{m}\\
        &\quad > 0.
    \end{align*}
\end{proof}

\begin{proof} [Proof of Statement \ref{state:z_fraction} in Lemma \ref{silemma:heur}.]
    By Definition \ref{def:probs_quantities_ell}, for $\ell < C_2$, we have
    \begin{align*}
        & \frac{\tilde{q}_\ell}{\tilde{p}_\ell}=\frac{q_k}{p_k}= \frac{k-1}{C_1 m}(k+\ell m) \leq \frac{k-1}{m} + \ell,
    \end{align*}where the second equality applies expressions of $q_k$ and $p_k$ in \eqref{def:qk} and \eqref{def:pk}, respectively, and the inequality holds since $k-1 < k \leq C_1$ and $m>0$ (recall \ref{def:ratio_of_rates}).
    The result now follows from the expression of $z_\ell$ in \eqref{def:sequence_z}.
\end{proof}

\begin{proof} [Proof of Statement \ref{state:z_ub} in Lemma \ref{lemma:heur}]
    Using the expression for $z_\ell$ where $\ell < C_2$ in  \eqref{def:sequence_z} yields
    \begin{align*}
        z_\ell = \frac{k-1}{m} + \ell \leq (k-1)+\ell = C_1-1, 
    \end{align*}where the inequality applies $m =\frac{\mu_2}{\mu_1}\geq 1$ ($m$ is defined in \eqref{def:ratio_of_rates}) since $\mu_2 \geq \mu_1$, and the last step follows by noticing $k+\ell = C_1$.
\end{proof}

\noindent Finally, we prove the remaining statements. These either compare $D(i,k,\ell)$ with an affine function closely related to $H(i,k,\ell)$ in particular parameter and state spaces or investigate properties of $i_D(k)$ and $i_H(k)$ (see Statements \ref{state:heur_linear_ub}–\ref{state:asymp} and Statement \ref{state:heur_linear_lb} in Lemma \ref{lemma:diff}).

\begin{proof} [Proof of Statement \ref{state:heur_linear_ub} in Lemma \ref{lemma:heur}]
    Recall the expression of $c$ given in \eqref{def:b_and_c}. It follows that $c=\frac{h_0}{C_1}\left(\frac{1}{\mu_1} - \frac{1}{\mu_2}\right)<0$ since $\mu_1 > \mu_2$. The proof proceeds by induction on $i$. When $i = 0$, by the boundary values \eqref{eq:diff_bdry}, 
    \begin{align*}
        D(0,k,\ell) = \frac{h_1}{\mu_1} - \frac{\max\{\ell+1,C_2\}h_2}{C_2\mu_2} \leq \frac{h_1}{\mu_2} - \frac{h_2}{\mu_2} = b \leq \big(0-y_{k}\big)c+b,
    \end{align*}where the last equality follows from the definition of $b$ in \eqref{def:b_and_c} and the last inequality holds since $c < 0$. Now suppose the result holds at $(i-1,k,\ell)$ for all $k$ and $\ell$ so that $(i-1,k,\ell) \in \tilde{X}_{diff}$, and consider it at $i$ for $i \geq 1$.
    If $D(i,k,\ell) \geq 0$, then by Proposition \ref{prop:larger_mu1_dec}, we also have $D(i-1,k,\ell) \geq 0$ given that $\mu_1 > \mu_2$. Consequently, the recursive equation \eqref{eq:recursive_pos} holds for $D(i,k,\ell)$ by Statement \ref{state:recursive_pos} in Lemma \ref{lemma:diff}. Additionally, the inductive hypothesis at $(i-1,k,\ell)$ implies 
    \begin{align}
        D(i-1,k,\ell) \leq (i-1-y_{k})c+b. \label{eq:induction_i-1_k_ell}
    \end{align}
    Next, observe that 
    \begin{align}
        \max\{0,D(i-1,k-1,\ell+1)\} \leq \max\{0,(i-1-y_{k-1})c+b\}. \label{eq:max_zero_H}
    \end{align}
    To establish this, we consider the following two cases.
    \begin{enumerate}[label= \textbf{Case} \arabic*:, leftmargin=3\parindent]
        \item If $D(i-1,k-1,\ell+1) \geq 0$, then the inductive hypothesis at $(i-1,k-1,\ell+1)$ implies that
        \begin{align*}
            \max\{0,D(i-1,k-1,\ell+1)\} & = D(i-1,k-1,\ell+1) \\
            \leq (i-1-y_{k-1})c+b & \leq \max\{0,(i-1-y_{k-1})c+b\}.
        \end{align*}
        \item If $D(i-1,k-1,\ell+1) < 0$,
        \begin{align*}
            \max\{0,D(i-1,k-1,\ell+1)\} = 0 \leq \max\{0,(i-1-y_{k-1})c+b\}.
        \end{align*}
    \end{enumerate}
    Substituting \eqref{eq:induction_i-1_k_ell} and \eqref{eq:max_zero_H} into the recursive equation \eqref{eq:recursive_pos} yields
    \begin{align}
        \lefteqn{D(i,k,\ell)= p_k(ic_k+b_k) + q_k \max\{0,D(i-1,k-1,\ell+1)\} + r_k D(i-1,k,\ell)}& \nonumber \\
        & \quad \leq p_k(ic_k+b_k) + q_k \max\{0,\left(i-1-y_{k-1}\right)c+b\} + r_k \big((i-1-y_{k})c+b\big) \nonumber \\
        & \quad \leq p_k(ic+b) + q_k\max\{0,(i-1-y_{k-1})c+b\} + r_k\big((i-1-y_{k})c+b\big), \label{eq:ub_max}
    \end{align}where the last inequality holds by \eqref{eq:cmp_bk_b_ck_c}. We consider the following two cases based on whether or not $i \leq i_{k-1}$, where $i_{k-1} := \left\lfloor\frac{b}{-c}+y_{k-1}\right\rfloor+1$.
    \begin{enumerate}[label= \textbf{Case} \arabic*:, leftmargin=3\parindent]
        \item If $i \leq i_{k-1}$, i.e., $(i-1-y_{k-1})c+b \geq 0$ (recall that $c<0$), then \eqref{eq:ub_max} becomes
        \begin{align}
        \lefteqn{D(i,k,\ell)\leq p_k(ic+b) + q_k \big((i-1-y_{k-1})c+b\big) + r_k \big((i-1-y_k)c+b\big)}& \nonumber \\
        & \quad = (p_k+q_k+r_k)(ic+b) - \big(q_k(y_{k-1}+1) + r_k(y_k+1)\big)c \nonumber \\
        &\quad = (ic+b) - \big(q_k(y_{k-1}+1) + r_k(y_k+1)\big)c, \label{eq:heur_linear_ub_case1_inte}
        \end{align}where the last equality follows by \eqref{eq:identity}.
        Consider the expression in the second term with coefficient $-c$ in \eqref{eq:heur_linear_ub_case1_inte}. Applying expressions of $q_k, r_k$ and $y_k$'s from \eqref{def:qk}, \eqref{def:rk} and \eqref{def:sequence_y}, respectively, results in 
        \begin{align}
            \lefteqn{q_k(y_{k-1}+1) + r_k(y_k+1)}& \nonumber \\ 
            &\quad = \frac{k-1}{(k-1)+\min\{\ell+1,C_2\}m} \big(k-2+\min\{\ell+1,C_2\}m + 1\big) \nonumber \\
            &\quad \qquad + \frac{\min\{\ell,C_2\} m}{k + \min\{\ell,C_2\}m} \big((k-1)+\min\{\ell,C_2\} m + 1\big) \nonumber \\
            &\quad = k-1+ \min\{\ell,C_2\} m \nonumber \\
            &\quad = y_k,
            \label{eq:y_rearrange}
        \end{align}where the last step relies on the definition of $y_k$ in \eqref{def:sequence_y} again.
        Applying \eqref{eq:y_rearrange} in \eqref{eq:heur_linear_ub_case1_inte} yields that if $i \leq i_{k-1}$ and $D(i,k,\ell) \geq 0$,
        \begin{align}
            D(i,k,\ell) \leq (i-y_k)c+b. \label{eq:heur_linear_ub_case1}
        \end{align}
        It remains to consider the case when $i \geq i_{k-1}+1$.
        \item Consider $i \geq i_{k-1}+1$, i.e., $(i-1-y_{k-1})c+b < 0$. We first prove $D(i_{k-1}+1,k,\ell) < 0$ by contradiction. Suppose, contrary to the claim, that $D(i_{k-1}+1,k,\ell) \geq 0$, then $D(i_{k-1},k,\ell) \geq 0$ by Proposition \ref{prop:larger_mu1_dec}. As a result, $D(i_{k-1},k,\ell) \leq (i_{k-1}-y_k)c+b$ by \eqref{eq:heur_linear_ub_case1} in the first case. Therefore, applying \eqref{eq:ub_max} at $i_{k-1}+1$ and noticing that $\big((i_{k-1} + 1)-1-y_{k-1}\big)c+b < 0$ yield
        \begin{align}
            \lefteqn{
            D(i_{k-1}+1,k,\ell) \leq p_k\big((i_{k-1}+1)c+b\big) + r_k\big((i_{k-1}-y_k)c+b\big)}& \nonumber \\
            &\quad = (p_k+r_k)\big((i_{k-1}+1)c+b\big) - r_k(1+y_k)c \nonumber \\
            &\quad = (p_k+r_k)\left(\left(i_{k-1}+1 - \frac{r_k}{p_k+r_k}(1+y_k)\right)c+b\right) \nonumber \\
            &\quad < (p_k+r_k)\left(\left(-\frac{b}{c}+y_{k-1}+1 - \frac{r_k}{p_k+r_k}(1+y_k)\right) \cdot c+b\right) \nonumber \\
            & \quad = (p_k+r_k)\left(y_{k-1}+1 - \frac{r_k}{p_k+r_k}(1+y_k)\right) \cdot c, \label{eq:case2_neg}
        \end{align}where the strict inequality holds since $i_{k-1} = \left\lfloor\frac{b}{-c}+y_{k-1}\right\rfloor+1 >\frac{b}{-c}+y_{k-1}$ and $c<0$.
        Observe that the terms in the second parenthesis in \eqref{eq:case2_neg} is non-negative for all $k$:
        \begin{align*}
            \lefteqn{y_{k-1}+1 - \frac{r_k}{p_k+r_k}(1+y_k)}&\\
            &\quad \geq y_k - \frac{r_k}{p_k+r_k}(1+y_k)\\
            &\quad = \frac{p_k}{p_k+r_k}\left(y_k-\frac{r_k}{p_k}\right)\\
            &\quad \geq 0,
        \end{align*}where the first inequality follows from Statement \ref{state:heur_mono_y} in Lemma \ref{silemma:heur}, and the last inequality holds by Statement \ref{state:y_fraction} in Lemma \ref{silemma:heur}.
        Together with $c<0$ in \eqref{eq:case2_neg}, it hence follows that $D(i_{k-1}+1,k,\ell) < 0$. This yields a contradiction to our hypothesis that $D(i_{k-1}+1,k,\ell) \geq 0$, verifying that $D(i_{k-1}+1,k,\ell) < 0$ is a must.
        
        Finally, we know $D(i,k,\ell) \leq D(i_{k-1}+1,k,\ell) <0$ for all $i \geq i_{k-1}+1$ by Proposition \ref{prop:larger_mu1_dec}. This in turn says that $D(i,k,\ell) \geq 0$ implies $i \leq i_{k-1}$.
    \end{enumerate}
    To conclude from both cases, $D(i,k,\ell) \geq 0$ implies that $i \leq i_{k-1}$ and $D(i,k,\ell)\leq (i-y_k)c+b$. 
    The proof now completes itself.
\end{proof}

\begin{proof} [Proof of Statement \ref{state:queue_in_collab_bd} in Lemma \ref{lemma:heur}.]
    Suppose first $h_0 \geq h_2$ and proceed by induction on $i$. For the base case $i = 0$, substituting the boundary values from \eqref{eq:diff_bdry} along with the expressions for $c',b'$ from \eqref{def:b'_and_c'} yields that for $\ell \geq C_2$,
    \begin{align}
        \lefteqn{D(0,k,\ell) - \big((0-y_k)c'+b'\big)}& \nonumber \\
        &\quad = \left(\frac{h_1}{\mu_1} - \frac{(\ell+1)h_2}{C_2\mu_2}\right) -\frac{h_0}{C_2\mu_2} y_k - \left(\frac{h_1-h_2}{\mu_1} - \frac{C_1h_2}{C_2\mu_2} \right) \nonumber \\
        &\quad = \frac{h_2}{\mu_1} + \frac{(k-1)h_2}{C_2\mu_2} - \frac{h_0}{C_2\mu_2} y_k \nonumber \\
        &\quad = \frac{h_2}{C_2 \mu_2} (C_2 m +k-1)  - \frac{h_0}{C_2\mu_2} y_k \nonumber \\
        &\quad = \frac{h_2-h_0}{C_2\mu_2}y_k, \label{eq:base_queue_in_collab}
    \end{align}where the last step follows from the fact that $y_k = k-1+C_2 m$, as stated in \eqref{def:sequence_y}. The result for $i=0$ then holds by applying $h_0 \geq h_2$ in \eqref{eq:base_queue_in_collab}. 
    
    Now, suppose that Statement \ref{state:queue_in_collab_costly_queue_ub} in Lemma \ref{lemma:heur} holds at $i-1$ for all $k$ and $\ell \geq C_2$ so that $(i-1,k,\ell) \in \tilde{\X}_{diff}$, and consider the case at $i$, for $i \geq 1$. The remainder of the proof follows a similar structure to that of Statement \ref{state:heur_linear_ub} in Lemma \ref{lemma:heur}, but with two key differences. First, all $c$ and $b$ are replaced with $c'$ and $b'$, respectively. In fact, in \eqref{eq:recursive_pos}, $b_k = b'$ and $c_k=c'<0$ when $\ell \geq C_2$ by their definitions in \eqref{def:bk} and \eqref{def:ck}, respectively.  Additionally, while the proof of Statement \ref{state:heur_linear_ub} in Lemma \ref{lemma:heur} relied solely on Proposition \ref{prop:larger_mu1_dec}, here we apply Proposition \ref{prop:larger_mu1_dec} when $\mu_1 \geq \mu_2$ and Proposition \ref{prop:larger_mu2_dec} when $\mu_1 < \mu_2$.

    Consider now $h_0 \leq h_2$ for the proof of Statement \ref{state:queue_in_collab_costly_collab_lb} in Lemma \ref{lemma:heur}, which is also carried out by induction on $i$. The base case $i = 0$ follows by replacing $h_0 \leq h_2$ in \eqref{eq:base_queue_in_collab}.
    Suppose the result holds at $i-1$ for all $k$ and $\ell \geq C_2$ so that $(i-1,k,\ell) \in \tilde{\X}_{diff}$, and consider it at $i$, where $i \geq 1$. We discuss based on whether or not $D(i-1,k,\ell) \geq 0$.
    \begin{enumerate}[label= \textbf{Case} \arabic*:, leftmargin=3\parindent]
        \item If $D(i-1,k,\ell) \geq 0$, then \eqref{eq:recursive_pos} holds by Statement \ref{state:recursive_pos} in Lemma \ref{lemma:diff}. Noting $c_k = c'$ and $b_k = b'$ when $\ell \geq C_2$ in \eqref{eq:recursive_pos} yields
        \begin{align}
        \lefteqn{D(i,k,\ell) = p_k(ic_k+b_k) + q_k \max\{0,D(i-1,k-1,\ell+1)\} + r_k D(i-1,k,\ell)}& \nonumber \\
        &\quad = p_k(ic'+b') + q_k \max\{0,D(i-1,k-1,\ell+1)\} + r_k D(i-1,k,\ell). \label{eq:queue_in_collab_costly_collab_lb_case1_inte}
        \end{align}
       Replacing the maximum in \eqref{eq:queue_in_collab_costly_collab_lb_case1_inte} with a lower bound $D(i-1,k-1,\ell+1)$ yields 
        \begin{align}
            \lefteqn{D(i,k,\ell) \geq p_k(ic'+b') + q_k D(i-1,k-1,\ell+1) + r_k D(i-1,k,\ell)}& \nonumber \\
            &\quad \geq p_k(ic'+b') + q_k \big((i-1-y_{k-1})c'+b'\big) + r_k \big((i-1-y_k)c'+b'\big) \nonumber \\
            &\quad = (p_k+q_k+r_k)(ic'+b') - \big(q_k(y_{k-1}+1) + r_k(y_k+1)\big)c' \nonumber \\
            &\quad = (i-y_k)c'+b', \label{eq:queue_in_collab_costly_collab_lb_case1}
        \end{align}where the second inequality uses inductive hypotheses at $i-1$ (twice) and the last step follows from the identities \eqref{eq:identity} and \eqref{eq:y_rearrange} (in the proof of Statement \ref{state:heur_linear_ub} in Lemma \ref{lemma:heur}).
        \item If $D(i-1,k,\ell) < 0$, then \eqref{eq:recursive_neg} holds by Statement \ref{state:recursive_neg} in Lemma \ref{silemma:heur}. Recall definitions of $\tilde{p}_\ell, \tilde{q}_\ell, \tilde{r}_\ell, \tilde{c}_\ell$ and $\tilde{b}_\ell$ in Definition \ref{def:probs_quantities_ell}, for $\ell \geq C_2$, equation \eqref{eq:recursive_neg} becomes
        \begin{align}
            \lefteqn{D(i,k,\ell) = \tilde{p}_\ell(i\tilde{c}_\ell+\tilde{b}_\ell) + \tilde{q}_\ell D(i-1,k,\ell) + \tilde{r}_\ell \min\{D(i-1,k+1,\ell-1), 0\}}& \nonumber \\
            &\quad = p_k(ic'+b') + q_k D(i-1,k,\ell) + r_k \min\{D(i-1,k+1,\ell-1), 0\}. \label{eq:ub_costly_queue_D_neg}
        \end{align}
        Noticing two inequalities $D(i-1,k,\ell) \geq D(i-1,k-1,\ell+1)$ and $D(i-1,k+1,\ell-1) \geq D(i-1,k,\ell)$ by Proposition \ref{prop:mono_ell} and applying them in \eqref{eq:ub_costly_queue_D_neg} gives
        \begin{align*}
            \lefteqn{D(i,k,\ell)\geq p_k(ic'+b') + q_k D(i-1,k-1,\ell+1) + r_k \min\{D(i-1,k,\ell), 0\}}&\\
            &\quad = p_k(ic'+b') + q_k D(i-1,k-1,\ell+1) + r_k D(i-1,k,\ell),
        \end{align*}where the last inequality relies on the assumption that $D(i-1,k,\ell) < 0$. The result follows by a comparable argument as \eqref{eq:queue_in_collab_costly_collab_lb_case1}.
    \end{enumerate}
    Both cases together confirm Statement \ref{state:queue_in_collab_costly_collab_lb} in Lemma \ref{lemma:heur}, which states that if $\ell \geq C_2$ and $h_0 \leq h_2$, then $D(i,k,\ell) \geq (i-y_k)c'+b'$.
\end{proof}

\begin{proof} [Proof of Statement \ref{state:cond_k-1} in Lemma \ref{lemma:heur}.]
    To prove this claim, notice that $y_k = y_{k-1}+1 \geq y_{k-1}$ for all $k = 2,3, \ldots, C_1-C_2$ (if one exists) by Statement \ref{state:heur_mono_y} in Lemma \ref{silemma:heur}, resulting in
    \begin{align}
        R_2(k) = -\frac{b'}{c'} + y_k = -\frac{b'}{c'} + y_{k-1}+1 = R_2(k-1)+1, \label{eq:zero_of_H}
    \end{align}
    where the first and last equalities apply the definition of $R_2(k)$ in \eqref{eq:thresh_queue_in_collab}.
    Now recall the expression of $i_H(k)$ in \eqref{def:iH_lowcost} for $\ell \geq C_2$ and $\frac{h_1}{\mu_1} > \frac{\ell+1}{C_2}\frac{h_2}{\mu_2}$. If $i_H(k) \geq 1$, then $\lceil R_2(k) \rceil \geq 1$, and hence $R_2(k) > 0$, which by \eqref{eq:zero_of_H}, implies $R_2(k-1) > -1$ and therefore $\lceil R_2(k-1) \rceil+1 \geq 1$. Hence, $i_H(k) = \lceil R_2(k) \rceil$ and $i_H(k-1) = \lceil R_2(k-1) \rceil $. Comparing with \eqref{eq:zero_of_H}, this leads to Statement \ref{state:sequence_iH} in Lemma \ref{lemma:heur}.
    
    Moreover, subtracting \eqref{eq:zero_of_H} from the equation in Statement \ref{state:sequence_iH} of this lemma yields
    \begin{align}
        i_H(k) - R_2(k) = i_H(k-1) - R_2(k-1).
        \label{eq:cond_highcost_collab_rhs}
    \end{align}This means the right-hand side of inequality \eqref{cond:highcost_collab} in Condition \ref{cond:sufficient_equal_thresh} remains the same for $k$ and $k-1$.
    Now subtracting $k$ from both sides of the equation in Statement \ref{state:sequence_iH} in Lemma \ref{lemma:heur} yields
    \begin{align}
        i_H(k) - k = i_H(k-1) - (k-1). 
        \label{eq:cond_highcost_collab_equiv1}
    \end{align}
    Additionally, referring to the expressions in Definition \ref{def:probs_k}, for $\ell \geq C_2$, we have
    \begin{align}
        \frac{r_k}{1-q_k} = \frac{k-1+C_2 m}{k+C_2 m}
    \end{align}is an increasing function in $k$. Consequently,
    \begin{align}
        \frac{r_{k-1}}{1-q_{k-1}} \leq \frac{r_k}{1-q_k}. \label{eq:cond_highcost_collab_LHS}
    \end{align}
    From $y_{k-1} \leq y_k$, equation \eqref{eq:cond_highcost_collab_equiv1}, and inequality \eqref{eq:cond_highcost_collab_LHS}, it follows that
    \begin{align}
        \frac{r_{k-1}}{1-q_{k-1}}y_{k-1} r^{i_H(k-1) - (k-1)} \leq \frac{r_k}{1-q_k}y_k r^{i_H(k) - k}. \label{eq:cond_highcost_collab_lhs}
    \end{align}
    Combining \eqref{eq:cond_highcost_collab_lhs} with \eqref{eq:cond_highcost_collab_rhs}, we conclude that if $k$ satisfies \eqref{cond:highcost_collab} in Condition \ref{cond:sufficient_equal_thresh}, then $k-1$ also satisfies it. A comparable argument applies to establish the corresponding result for \eqref{cond:highcost_queue} in Condition \ref{cond:sufficient_equal_thresh}.
\end{proof}

\begin{proof} [Proof of Statement \ref{state:asymp} in Lemma \ref{lemma:heur}.] 
    Before proceeding with the proof, we collect several observations that apply symmetrically to $h_0 >h_2$ and $h_0<h_2$:
    \begin{enumerate}
        \item 
        If $h_0 > h_2$, then inequality \eqref{cond:highcost_queue} in Condition \ref{cond:sufficient_equal_thresh} is also equivalent to
    \begin{align}
        (i_H(k)-1-y_k\big)c'+b' + \frac{h_2-h_0}{C_2 \mu_2}y_{k}r^{i_H(k)-k} > 0. \label{eq:positive_D_iH(k)-1}
    \end{align}
    Similarly, if $h_0 < h_2$, inequality  \eqref{cond:highcost_collab} in Condition \ref{cond:sufficient_equal_thresh} is the same as
    \begin{align}
        \big(i_H(k)-y_k\big)c'+b' + \frac{r_k}{1-q_k}\frac{h_2-h_0}{C_2 \mu_2}y_{k}r^{i_H(k)-k} \leq 0. \label{eq:negative_D_iH(k)}
    \end{align}
    The equivalence of \eqref{eq:positive_D_iH(k)-1} and \eqref{eq:negative_D_iH(k)} to \eqref{cond:highcost_queue} and \eqref{cond:highcost_collab} in Condition \ref{cond:sufficient_equal_thresh}, respectively, follows from some arithmetic using the definition of $R_2(k)$ (see \eqref{eq:thresh_queue_in_collab}).
    \item \label{obs:pos_iD}
    We have $i_D(k) \geq 1$ by Statement \ref{state:always_neg_queue_in_collab} in Theorem \ref{thm:heur_lowcost}, since $\frac{h_1}{\mu_1} > \frac{\ell+1}{C_2} \frac{h_2}{\mu_2}$, which implies $D(0,k,\ell) >0$ (recall the definition of $i_D(k)$ in \eqref{def:iD_lowcost}).
    \item For all $i$ such that $1\leq i \leq i_D(k)$ --- or equivalently $D(i-1,k,\ell) > 0$ (recall the definition of $i_D(k)$ in \eqref{def:iD_lowcost}) --- the following equation holds:
    \begin{align}
        D(i,k,\ell) = p_k(ic'+b') + q_k \max\{0,D(i-1,k-1,\ell+1)\} + r_k D(i-1,k,\ell). \label{eq:recursive_queue_in_collab}
    \end{align}
    This follows by applying $c_k = c'$ and $b_k = b'$ when $\ell \geq C_2$ (see their definitions in \eqref{def:ck} and \eqref{def:bk}, respectively) in \eqref{eq:recursive_pos}. 
    \item The following equation holds for all $i$ such that $0 \leq i \leq i_D(1)$:
    \begin{align}
        D(i,1,C_1-1) = \left(i - y_1\right)c' + b' + \frac{h_2-h_0}{C_2 \mu_2}y_1 r^i. \label{eq:base_asymp}
    \end{align}
    Consider the proof of \eqref{eq:base_asymp} for $i=0$ first. From \eqref{eq:base_queue_in_collab} in the proof of Statement \ref{state:queue_in_collab_bd} in Lemma \ref{lemma:heur}, which holds under the assumption $\ell \geq C_2$ only, we know that (applied in particular at $k=1$)
    \begin{align*}
        D(0,k,\ell) = \big((0-y_k)c'+b'\big) + \frac{h_2-h_0}{C_2 \mu_2}y_1.
    \end{align*}
    For any $i$ such that $1\leq i \leq i_D(1)$, applying $q_1 = 0$ and $p_1+r_1 = 1$ (recall their expressions in Definition \ref{def:probs_k}) in \eqref{eq:recursive_queue_in_collab} leads to 
    \begin{align}
         D(i,1,C_1-1) = p_1 (i c' + b') + r_1 D(i-1,1,C_1-1). \label{eq:recursive1}
    \end{align} Further expanding this by reapplying \eqref{eq:recursive1} at $i-1,i-2, \ldots, 1$ yields 
    \begin{align*}
        & D(i,1,C_1-1) = p_1 (i c' + b') + r_1 \Big(p_1 \big((i-1) c' + b'\big) + r_1 D(i-2,1,C_1-1)\Big)\\
        &\quad = \left(i + r_1 (i-1)\right) p_1 c' + \left(1 + r_1\right)p_1 b' + (r_1)^2 D(i-2,1,C_1-1)\\
        &\quad = \ldots\\
        &\quad = \left(i + r_1 (i-1) + \ldots + (r_1)^{i-1}\right) p_1 c' + \left(1 + r_1 + \ldots + (r_1)^{i-1}\right)p_1 b' + (r_1)^i D(0,1,C_1-1)\\
        &\quad = \left(\frac{i}{1-r_1} - r_1\frac{1-(r_1)^i}{(1-r_1)^2}\right) p_1 c' + \frac{1-(r_1)^i}{1-r_1}p_1 b' + (r_1)^i D(0,1,C_1-1).
    \end{align*}A substitution with $p_1+r_1=1$ yields
    \begin{align}
         \lefteqn{D(i,1,C_1-1) = \left(i - r_1\frac{1-(r_1)^i}{1-r_1}\right)c' + \big(1-(r_1)^i\big) b' + (r_1)^i D(0,1,C_1-1)}&\nonumber\\
         &\quad = \left(i - \frac{r_1}{1-r_1}\right)c' + b' + (r_1)^i \left(D(0,1,C_1-1)- \left(-\frac{r_1}{1-r_1}c' + b'\right)\right) \nonumber \\
         &\quad = \left(i - y_1\right)c' + b' + (r_1)^i \big(D(0,1,C_1-1)- (-y_1c' + b')\big) \nonumber \\
         &\quad = \left(i - y_1\right)c' + b' + \frac{h_2-h_0}{C_2 \mu_2}y_1(r_1)^i, \nonumber
    \end{align}where the second-to-last equality follows from the identity that $y_1 = C_2 m = \frac{r_1}{1-r_1}$ (using the expressions of $y_1$ in \eqref{def:sequence_y} and $r_1$ in \eqref{def:rk}), and the last step applies \eqref{eq:base_queue_in_collab} at $k=1$ again. The result follows by noticing $r = r_1$ (see the definition of $r$ in Condition \ref{cond:sufficient_equal_thresh}).
    \end{enumerate}
    
    With the supporting facts above, the proof of Statement \ref{state:asymp} in Lemma \ref{lemma:heur} is now in reach. Consider first $h_0 > h_2$. 
    Note that Statement \ref{state:thresh_queue_in_collab} in Theorem \ref{thm:heur_lowcost} confirms that for any $k$ such that $\ell = C_1-k \geq C_2$, it holds that
    \begin{align}
        i_H(k) \geq i_D(k) \geq 1, \label{eq:thresh_queue_in_collab_costly_queue}
    \end{align}where the last inequality holds by Observation \ref{obs:pos_iD}.
    The proof is carried out by induction on $k$. For the base case where $k=1$, Statement \ref{state:exp_bound} in Lemma \ref{lemma:heur} follows immediately from \eqref{eq:base_asymp} by defining $\varphi_1(i) = \frac{h_2-h_0}{C_2 \mu_2}y_1r^i$. 
    It remains to show Statement \ref{state:thresh_reverse} to complete the proof for $k=1$. 
    Suppose towards a contradiction that $i_D(1) \leq i_H(1)-1$ assuming that $k=1$ satisfies \eqref{cond:highcost_queue} in Condition \ref{cond:sufficient_equal_thresh} (or \eqref{eq:positive_D_iH(k)-1}). Define the following function of $i$
    \begin{align*}
        \phi(i) := \left(i - y_1\right)c' + b' + \frac{h_2 - h_0}{C_2 \mu_2} y_1 r^i, \quad \forall i \in \mathbb{Z_+}.
    \end{align*}
    By \eqref{eq:base_asymp}, we have $D(i,1,C_1-1) = \phi(i)$ for all $0 \leq i \leq i_D(1)$, where $i_D(1) \geq 1$. Consequently, by the definition of $i_D(1)$,
    \begin{align}
        \phi\big(i_D(1)\big) = D\big(i_D(1),1,C_1-1\big) \leq 0. \label{eq:contradiction}
    \end{align}
    Next, consider the difference
    \begin{align*}
        \phi(i+1) - \phi(i) = c' + \frac{h_2-h_0}{C_2 \mu_2}y_1r^i(r-1).
    \end{align*}
    Since $h_0 > h_2$ and $r< 1$ as given in Definition \ref{cond:sufficient_equal_thresh}, this expression is decreasing in $i$. It follows that 
    \begin{align*}
        \phi(i+1) - \phi(i) \leq \phi(1) - \phi(0) = D(1,1,C_1-1) - D(0,1,C_1-1) \leq 0,
    \end{align*}where the last inequality follows from Proposition \ref{prop:larger_mu1_dec} if $\mu_1\geq \mu_2$ and from Proposition \ref{prop:larger_mu2_dec} if $\mu_2 > \mu_1$. Therefore, $\phi(i)$ is non-increasing in $i$. Together with the assumption that $i_D(1) \leq i_H(1)-1$, where $i_H(1)-1 \geq 0$ (recall $i_H(1) \geq 1$ by \eqref{eq:thresh_queue_in_collab_costly_queue}), this leads to 
    \begin{align*}
        \phi\big(i_D(1)\big) \geq \phi\big(i_H(1)-1\big)
        = \big(i_H(1)-1 - y_1\big)c' + b' + \frac{h_2-h_0}{C_2 \mu_2}y_1 r^{i_H(1)-1} > 0,
    \end{align*}where the last step follows from \eqref{eq:positive_D_iH(k)-1}. Comparing with \eqref{eq:contradiction}, we reach a contradiction, proving that $i_D(1) \geq i_H(1)$. The proof for the base case $k=1$ now completes itself.

    Assuming Statements \ref{state:exp_bound} and \ref{state:thresh_reverse} in Lemma \ref{lemma:heur} hold at $k-1$, we verify them at $k$, where $k \geq 2$ such that $\ell = C_1-k \geq C_2$.
    We proceed with showing the following equation:
    \begin{align}
        i_D(k-1) = i_H(k-1). \label{eq:equal_thresh_induction}
    \end{align}
    This is true if $i_D(k-1) = 0$ by Statement \ref{state:always_neg_queue_in_collab} in Theorem \ref{thm:heur_lowcost}. If $i_D(k-1) \geq 1$, the inductive hypothesis of Statement \ref{state:thresh_reverse} in Lemma \ref{lemma:heur} leads to $i_H(k-1) \leq i_D(k-1)$ and the reverse direction is by \eqref{eq:thresh_queue_in_collab_costly_queue}.
    The proof proceeds with an inner induction on $i$, where $0 \leq i \leq i_D(k)-1$, to prove Statement \ref{state:exp_bound}. 
    When $i = 0$ at $k$, the result holds equation \eqref{eq:base_queue_in_collab} again:
    \begin{align*}
        \varphi_k(0) = D(0,k,\ell) - H_k(0) = \frac{h_2-h_0}{C_2 \mu_2}y_k \geq \frac{h_2-h_0}{C_2 \mu_2}y_k r^{-k+1}, 
    \end{align*}where the inequality follows since $h_0 > h_2$ and $k \geq 1$.  
    Suppose Statement \ref{state:exp_bound} in Lemma \ref{lemma:heur} holds at $i-1$ (for $k$) and consider it at $i$ (for $k$), where $1 \leq i \leq i_D(k)-1$.
    Consider the following two cases based on whether or not $i \leq i_D(k-1) = i_H(k-1)$.
    \begin{enumerate}[label= \textbf{Case} \arabic*:, leftmargin=3\parindent]
        \item If $i \leq i_H(k-1) = i_D(k-1)$ (and $1 \leq i \leq i_D(k)-1$), equation \eqref{eq:recursive_queue_in_collab} becomes:
        \begin{align}
            & D(i,k,\ell) = p_k (i c' + b') + q_k D(i-1,k-1,\ell+1) + r_k D(i-1,k,\ell) \nonumber \\
            &\quad = p_k (i c' + b') + q_k \Big(\big(i-1 - y_{k-1}\big)c' + b' +  \varphi_{k-1}(i-1)\Big) \nonumber \\
            &\quad \qquad + r_k \Big(\big(i-1 - y_{k}\big)c' + b' +  \varphi_{k}(i-1)\Big) \nonumber \\
            &\quad = (p_k+q_k+r_k)(ic' +b') - \big(q_k(1+y_{k-1}) + r_k(1+y_k)\big)c' \nonumber \\
            &\quad \qquad  + q_k\varphi_{k-1}(i-1) + r_k\varphi_{k}(i-1), \label{eq:asymp_induction}
        \end{align}where $\varphi_{k-1}(i-1) \geq \frac{h_2-h_0}{C_2 \mu_2}y_{k-1} r^{(i-1)-(k-1)+1}$ and $\varphi_k(i-1) \geq \frac{h_2-h_0}{C_2 \mu_2}y_k r^{(i-1)-k+1}$ in the second equality. This applies inductive hypotheses of Statement \ref{state:exp_bound} in Lemma \ref{lemma:heur} twice at $k-1$ and at $i-1$ for $k$, respectively.
        Using identities \eqref{eq:identity} and \eqref{eq:y_rearrange} in \eqref{eq:asymp_induction} yields 
        \begin{align}
            & D(i,k,\ell)= (ic' +b') - y_kc' + q_k\cdot\varphi_{k-1}(i-1) + r_k\cdot\varphi_{k}(i-1)\nonumber \\
            &\quad = (i-y_k)c'+b' + \varphi_{k}(i), \label{eq:recursive_asymp}
        \end{align}by letting $\varphi_{k}(i):=q_k\cdot\varphi_{k-1}(i-1) + r_k\cdot\varphi_{k}(i-1)$.
        In addition,
        \begin{align}
            \lefteqn{\varphi_k(i) = q_k\cdot\varphi_{k-1}(i-1) + r_k\cdot\varphi_{k}(i-1)}& \nonumber \\
            & \quad\geq  q_k \frac{h_2-h_0}{C_2 \mu_2}y_{k-1} r^{(i-1)-(k-1)+1} + r_k \frac{h_2-h_0}{C_2 \mu_2}y_{k} r^{(i-1)-k+1} \nonumber \\ 
            & \quad= \frac{h_2-h_0}{C_2 \mu_2}  y_k r^{i-k+1} \left(q_k \frac{y_{k-1}}{y_k} +\frac{r_k }{r}\right)
            \label{eq:phi_induction_costly_queue}
        \end{align}
        Consider the term in the parenthesis. Using expressions of $y_k$ in \eqref{def:sequence_y}, as well as $q_k$ and $r_k$ in \eqref{def:qk} and \eqref{def:rk}, respectively, from Definition \ref{def:probs_k} for $\ell \geq C_2$, we obtain
        \begin{align}
            \lefteqn{q_k \frac{y_{k-1}}{y_k} +\frac{r_k}{r}}& \nonumber\\
            &\quad = \frac{k-1}{(k-1)+ C_2 m} \frac{k-2+C_2 m}{k-1+C_2 m} + \frac{1+C_2 m}{k + C_2 m} \nonumber\\
            &\quad \leq \frac{k-1}{(k-1)+ C_2 m} \frac{k-1+C_2 m}{k+C_2 m} + \frac{1+C_2 m}{k + C_2 m} \nonumber\\
            &\quad = 1, \label{eq:asymp_phi_induction}
        \end{align}where the inequality relies on the fact that $\frac{k-1+C_2 m}{k+C_2 m}$ is increasing in $k$.
        Applying $h_0 > h_2$ and \eqref{eq:asymp_phi_induction} in \eqref{eq:phi_induction_costly_queue} and yields 
        \begin{align}
            \varphi_k(i) \geq \frac{h_2-h_0}{C_2 \mu_2}  y_k r^{i-k+1}. \label{eq:asymp_case1}
        \end{align}
        \item Consider now $i \geq i_H(k-1)+1 = i_D(k-1)+1$, that is, $D(i-1,k-1,\ell+1) \leq 0$. 
        Observe that Statement \ref{state:sequence_iH} in Lemma \ref{lemma:heur}, along with \eqref{eq:thresh_queue_in_collab_costly_queue}, implies
        \begin{align}
            i_H(k-1)+1 = i_H(k) \geq i_D(k), \label{eq:iH_inc}
        \end{align}leading to $i \geq i_D(k)$ in this case.
    \end{enumerate}
    Combining both cases, Statement \ref{state:exp_bound} in Lemma \ref{lemma:heur} is now established for $h_0 > h_2$. It remains to show $i_H(k) \leq i_D(k)$ to complete the inductive proof of Statement \ref{state:thresh_reverse}.
    The proof reduces to verifying that $D\big(i_H(k)-1,k,\ell\big) > 0$, or $D\big(i_H(k-1),k,\ell\big) > 0$ (recall $i_H(k-1)+1 = i_H(k)$ by \eqref{eq:iH_inc}). If $i_H(k) = 1$, or equivalently, $i_H(k-1) = 0$, $D(0,k,\ell) > 0$ holds by Observation \ref{obs:pos_iD}. Therefore, it suffices to consider the case where $i_H(k-1) \geq 1$.
    Notice that Proposition \ref{prop:mono_ell} implies the following 
    \begin{align}
        \lefteqn{D\big(i_H(k-1)-1,k,\ell\big)}& \nonumber \\
        &\quad \geq D\big(i_H(k-1)-1,k-1,\ell+1\big) \nonumber \\
        &\quad = D\big(i_D(k-1)-1,k-1,\ell+1\big) \nonumber \\
        &\quad > 0, \label{eq:pos_terms_in_recursion}
    \end{align}where the equality substitutes \eqref{eq:equal_thresh_induction} and the last step applies the definition of $i_D(k-1)$ (see \eqref{def:iD_lowcost}). It follows from \eqref{eq:pos_terms_in_recursion} that $D\big(i_H(k-1),k,\ell\big)$ satisfies \eqref{eq:recursive_queue_in_collab}:
    \begin{align*}
        \lefteqn{D\big(i_H(k-1),k,\ell\big)}& \\
        &\quad = p_k \big(i_H(k-1) c' + b'\big) + q_k D\big(i_H(k-1)-1,k-1,\ell+1\big) + r_k D\big(i_H(k-1)-1,k,\ell\big)\\
        &\quad = p_k \big(i_H(k-1) c' + b'\big) + q_k \Big(\big(i_H(k-1)-1 - y_{k-1}\big)c' + b' +  \varphi_{k-1}\big(i_H(k-1)-1\big)\Big) \nonumber \\
        &\quad \qquad + r_k \Big(\big(i_H(k-1)-1 - y_{k}\big)c' + b' +  \varphi_{k}\big(i_H(k-1)-1\big)\Big),
    \end{align*}where in the last step, we have $ \varphi_{k-1}\big(i_H(k-1)-1\big) \geq \frac{h_2-h_0}{C_2 \mu_2}y_{k-1} r^{\big(i_H(k-1)-1\big)-(k-1)+1}$ and $\varphi_{k}\big(i_H(k-1)-1\big) \geq \frac{h_2-h_0}{C_2 \mu_2}y_k r^{\big(i_H(k-1)-1\big)-k+1}$, which applies Statement \ref{state:exp_bound} in Lemma \ref{lemma:heur} since $i_H(k-1) -1 \leq i_D(k)-1<i_D(k)$ by  \eqref{eq:pos_terms_in_recursion}. 
     
    Now, following the same arithmetic steps used to derive \eqref{eq:recursive_asymp} and \eqref{eq:asymp_case1}, substituting $i = i_H(k-1)$ yields 
     \begin{align}
         D\big(i_H(k-1),k,\ell\big) = \big(i_H(k-1) - y_{k}\big)c' + b' +  \varphi_{k}\big(i_H(k-1)\big), \label{eq:D_iH(k-1)}
     \end{align}where 
    \begin{align}
        & \varphi_{k}\big(i_H(k-1)\big) := q_k\cdot\varphi_{k-1}\big(i_H(k-1)-1\big) + r_k\cdot\varphi_{k}\big(i_H(k-1)-1\big)\nonumber \\
        &\quad \geq \frac{h_2-h_0}{C_2 \mu_2}  y_k r^{i_H(k-1)-k+1}. \label{eq:phi_iH(k-1)}
    \end{align}
    Note that while \eqref{eq:recursive_asymp} and \eqref{eq:phi_induction_costly_queue} apply the inductive hypotheses of Statement \ref{state:exp_bound} in Lemma \ref{lemma:heur}, equations \eqref{eq:D_iH(k-1)} and \eqref{eq:phi_iH(k-1)} instead utilize its established result. Finally, substituting \eqref{eq:phi_iH(k-1)} in \eqref{eq:D_iH(k-1)} leads to
    \begin{align*}
        & D\big(i_H(k-1),k,\ell\big) \geq \big(i_H(k-1) - y_{k}\big)c' + b' +  \frac{h_2-h_0}{C_2 \mu_2}  y_k r^{i_H(k-1)-k+1}\\
        &\quad = \big(i_H(k)-1 - y_{k}\big)c' + b' +  \frac{h_2-h_0}{C_2 \mu_2}  y_k r^{i_H(k)-k}\\
        &\quad > 0,
    \end{align*}where the equality follows from Statement \ref{state:sequence_iH} in lemma \ref{lemma:heur} and the last step applies \eqref{eq:positive_D_iH(k)-1}. This confirms Statement~\ref{state:thresh_reverse} in the lemma for the case $h_0 > h_2$, completing the desired result.
    
    Now, consider the second case where $h_0 < h_2$. The proofs of both Statements \ref{state:exp_bound} and \ref{state:thresh_reverse} proceed analogously to the case $h_0 >h_2$, following induction on $k$.
    For the base case $k=1$, defining $\varphi_1(i) = \frac{h_2-h_0}{C_2 \mu_2}y_1r^i$ once again establishes Statement \ref{state:exp_bound} via \eqref{eq:base_asymp}.
    It remains to show \eqref{cond:highcost_collab} in Condition \ref{cond:sufficient_equal_thresh} (or \eqref{eq:negative_D_iH(k)}) implies $i_H(1) \geq i_D(1)$.
    In fact, we obtain that $i_H(1) \leq i_D(1)$ from Statement \ref{state:thresh_queue_in_collab} in Theorem \ref{thm:heur_lowcost}, since $h_0 < h_2$. Therefore, $D\big(i_H(1),1,C_1-1\big)$ satisfies \eqref{eq:base_asymp}:
    \begin{align*}
        \lefteqn{D\big(i_H(1),1,C_1-1\big) = \big(i_H(1) - y_1\big)c' + b' + \frac{h_2-h_0}{C_2 \mu_2}y_1 r^{i_H(1)}}&\\
        &\quad = \big(i_H(1) - y_1\big)c' + b' + \frac{h_2-h_0}{C_2 \mu_2}y_1r^{i_H(1)-1}\frac{r_1}{1-q_1},
    \end{align*}where the second equality follows since $q_1 = 0$ and $r_1=r$. The result follows by noticing that $i_H(1) \geq i_D(1)$ holds if and only if $D\big(i_H(1),1,C_1-1\big) \leq 0$, and is further equivalent to \eqref{eq:negative_D_iH(k)} (applied at $k=1$). The proof of both statements for $k=1$ is now complete.
    
    Next, suppose the results hold at $k-1$ and consider them at $k$, where $k \geq 2$ such that $\ell = C_1-k \geq C_2$. We prove $i_H(k) \geq 1$ by contradiction. Suppose, contrary to the claim, that $i_H(k) = 0$. Substituting this into \eqref{eq:negative_D_iH(k)} yields
    \begin{align}
        0 & \geq  (0-y_k)c'+b' + \frac{r_k}{1-q_k}\frac{h_2-h_0}{C_2 \mu_2}y_{k}r^{0-k} \nonumber \\
        &\quad = -y_kc'+b' + \frac{r_k}{r(1-q_k)}\frac{h_2-h_0}{C_2 \mu_2}y_{k}r^{1-k} \nonumber \\
        &\quad \geq -y_kc'+b' + \frac{h_2-h_0}{C_2 \mu_2}y_{k}, \label{eq:iH_contra1}
    \end{align}where the final inequality utilizes $h_0 <h_2$, $k \geq 1$ and the following fact:
    \begin{align*}
        \frac{r_k}{r(1-q_k)} = \frac{k-1+C_2 m}{k+C_2 m} \frac{1+C_2 m}{C_2 m} \geq 1.
    \end{align*}
    However, Observation \ref{obs:pos_iD} implies
    $0 < D\big(0,k,\ell\big) = -y_kc'+b' + \frac{h_2-h_0}{C_2 \mu_2}y_{k}$, where the equality follows by \eqref{eq:base_queue_in_collab}. A contradiction emerges when comparing with \eqref{eq:iH_contra1}, which confirms that $i_H(k)\geq 1$. Consequently, Statement \ref{state:cond_k-1} in Lemma \ref{lemma:heur} follows.
    
    An inner induction on $i$ is adopted to prove Statement \ref{state:exp_bound}, where $0 \leq i \leq i_D(k)-1$ ($i_D(k) \geq 1$ by Observation \ref{obs:pos_iD}). The base case when $i = 0$ at $k$ holds by \eqref{eq:base_queue_in_collab} once again:
    \begin{align*}
        \varphi_k(0) = D(0,k,\ell) - \big((0-y_k)c'+b'\big) = \frac{h_2-h_0}{C_2 \mu_2}y_k \leq \frac{h_2-h_0}{C_2 \mu_2}y_k r^{-k+1}, 
    \end{align*}where the inequality follows since $h_0 < h_2$ and $k \geq 1$.  
    Assuming that Statement \ref{state:exp_bound} holds at $i-1$ for $k$, we now analyze it at $i \geq 1$ for $k$. Observe that \eqref{eq:equal_thresh_induction} remains valid by a symmetric argument. 
    Combining \eqref{eq:equal_thresh_induction} with Statement \ref{state:sequence_iH} in Lemma \ref{lemma:heur} yields 
    \begin{align}
        i_D(k-1) = i_H(k-1) = i_H(k)-1, \label{eq:iH_inc2}
    \end{align}
    leading to the following two cases to consider.
    \begin{enumerate}[label= \textbf{Case} \arabic*:, leftmargin=3\parindent]
        \item When $1\leq i \leq i_H(k)-1=i_H(k-1)$ (and $1 \leq i \leq i_D(k)-1$), following the comparable arithmetic as \eqref{eq:asymp_induction}--\eqref{eq:asymp_phi_induction}, we obtain that equation \eqref{eq:recursive_asymp} in Lemma \ref{lemma:heur} still holds. Additionally, 
        \begin{align*}
            \lefteqn{\varphi_k(i) = q_k\cdot\varphi_{k-1}(i-1) + r_k\cdot\varphi_{k}(i-1)}& \nonumber \\
            & \quad\leq  q_k \frac{h_2-h_0}{C_2 \mu_2}y_{k-1} r^{(i-1)-(k-1)+1} + r_k \frac{h_2-h_0}{C_2 \mu_2}y_{k} r^{(i-1)-k+1} \nonumber \\ 
            & \quad= \frac{h_2-h_0}{C_2 \mu_2}  y_k r^{i-k+1} \left(q_k \frac{y_{k-1}}{y_k} +\frac{r_k }{r}\right)\nonumber\\
            &\quad \leq \frac{h_2-h_0}{C_2 \mu_2}  y_k r^{i-k+1}, 
        \end{align*}where the first inequality follows from inductive hypotheses twice of Statement \ref{state:exp_bound} at $k-1$ and at $i-1$ for $k$, and the last step holds by $h_0 < h_2$ and \eqref{eq:asymp_phi_induction}.

        In particular, since $i_H(k) -1 \leq i_D(k) -1 $ by Statement \ref{state:thresh_queue_in_collab} in Theorem \ref{thm:heur_lowcost}, we obtain\begin{align}
            D\big(i_H(k)-1,k,\ell\big) = \big(i_H(k)-1 - y_{k}\big)c' + b' +  \varphi_{k}\big(i_H(k)-1\big),  \label{eq:D_iH(k)-1}
        \end{align}where
        \begin{align}
            \varphi_{k}\big(i_H(k)-1\big) \leq \frac{h_2-h_0}{C_2 \mu_2}  y_k r^{\big(i_H(k)-1\big)-k+1}. \label{eq:phi_iH(k)-1}
        \end{align}
        Additionally, $D\big(i_H(k),k,\ell\big)$ satisfies \eqref{eq:recursive_queue_in_collab}.
        \item Consider $i \geq i_H(k)$. We prove that if $k$ satisfies \eqref{cond:highcost_collab} in Condition \ref{cond:sufficient_equal_thresh}, then $D\big(i_H(k),k,\ell\big) \leq 0$, so that $i_H(k) \geq i_D(k)$. 
        Applying both $D\big(i_H(k)-1,k-1,\ell+1\big) \leq 0$ (by \eqref{eq:iH_inc2}) and \eqref{eq:D_iH(k)-1} in \eqref{eq:recursive_queue_in_collab} gives
        \begin{align*}
            \lefteqn{D\big(i_H(k),k,\ell\big)= p_k \big(i_H(k) c' + b'\big) + r_k D\big(i_H(k)-1,k,\ell\big)}&\\
            &\quad =  p_k \big(i_H(k) c' + b'\big) + r_k \Big(\big(i_H(k)-1 - y_{k}\big)c' + b' +  \varphi_{k}\big(i_H(k)-1\big)\Big)\\
            &\quad = (p_k+r_k)\left(\left(i_H(k)-\frac{r_k}{p_k+r_k}(1+y_k)\right) \cdot c'+b' + \frac{r_k}{p_k+r_k} \varphi_{k}\big(i_H(k)-1\big)\right)\\
            &\quad = (p_k+r_k)\left(\big(i_H(k)-y_k\big) \cdot c'+b' + \frac{r_k}{p_k+r_k} \varphi_{k}\big(i_H(k)-1\big)\right).
        \end{align*}Here, the final step holds because 
        \begin{align*}
            \frac{r_k}{p_k+r_k}(1+y_k) = \frac{r_k}{p_k+r_k}\left(1+\frac{r_k}{p_k}\right) = \frac{r_k}{p_k} = y_k,
        \end{align*}where the first and last equalities follow from Statement \ref{state:y_fraction} in Lemma \ref{lemma:heur} for the case $\ell \geq C_2$.
        Finally, $D\big(i_H(k),k,\ell\big) \leq 0$ because the following term is non-positive:
        \begin{align*}
            \lefteqn{\left(i_H(k)-y_k\right) \cdot c'+b' + \frac{r_k}{p_k+r_k} \varphi_{k}\big(i_H(k)-1\big)}&\\
            &\quad = \big(i_H(k)-y_k\big) \cdot c'+b' + \frac{r_k}{1-q_k} \varphi_{k}\big(i_H(k)-1\big)\\
            &\quad \leq \big(i_H(k)-y_k\big) \cdot c'+b' + \frac{r_k}{1-q_k} \frac{h_2-h_0}{C_2 \mu_2}y_k r^{i_H(k)-k}\\
            &\quad \leq 0,
        \end{align*}where the equality applies the identity \eqref{eq:identity}, the second-to-last inequality uses \eqref{eq:phi_iH(k)-1} and the final step is a consequence of \eqref{eq:negative_D_iH(k)}. This confirms that $i_H(k) \geq i_D(k)$, thereby establishing Statement \ref{state:thresh_reverse} in Lemma \ref{lemma:heur}.
    \end{enumerate}
    Combining both cases, Statement \ref{state:exp_bound} in Lemma \ref{lemma:heur} is thus verified, completing the proof.
\end{proof}

\begin{proof} [Proof of Statement \ref{state:heur_linear_lb} in Lemma \ref{lemma:heur}]
    The proof follows analogously to that of Statement \ref{state:heur_linear_ub} in Lemma \ref{lemma:heur} by induction on $i$. 
    In this case, the definition of $c$ given in \eqref{def:b_and_c} yields that $c=\frac{h_0}{C_1}\left(\frac{1}{\mu_1} - \frac{1}{\mu_2}\right)>0$ because $\mu_1 < \mu_2$. 
    For the base case $i = 0$, the boundary values \eqref{eq:diff_bdry} for $\ell < C_2$ imply that
    \begin{align*}
        D(0,k,\ell) = \frac{h_1}{\mu_1} - \frac{h_2}{\mu_2} = b \geq \big(0-z_\ell\big)c+b,
    \end{align*}where the second equality follows from the definition of $b$ in \eqref{def:b_and_c} and the last inequality holds by noticing that $c > 0$. Now suppose the result holds at $i-1$ for all $k$ and $\ell<C_2$ such that $(i-1,k,\ell) \in \tilde{\X}_{diff}$, and consider it at $i$, where $i \geq 1$. For $i \geq 1$, Proposition \ref{prop:larger_mu2_inc} ensures that $D(i,k,\ell) \leq 0$ implies $D(i-1,k,\ell) \leq 0$, so $D(i,k,\ell)$ satisfies the recursive equation \eqref{eq:recursive_neg} by Statement \ref{state:recursive_neg} in Lemma \ref{silemma:heur}. Additionally, by the inductive hypothesis at $(i-1,k,\ell)$, we have
    \begin{align}
        D(i-1,k,\ell) \geq (i-1-z_\ell)c+b. \label{eq:induction_lb_i-1_k_ell}
    \end{align}
    We then show that 
    \begin{align}
        \min\{D(i-1,k+1,\ell-1),0\} \geq \min\{(i-1-z_{\ell-1})c+b,0\}  \label{eq:min_zero_H}
    \end{align}
    by considering the following two cases. 
    \begin{enumerate}[label= \textbf{Case} \arabic*:, leftmargin=3\parindent]
        \item If $D(i-1,k+1,\ell-1) \leq 0$, then the inductive hypothesis at $(i-1,k+1,\ell-1)$ implies that
        \begin{align*}
            \lefteqn{\min\{D(i-1,k+1,\ell-1),0\} = D(i-1,k+1,\ell-1)}\\
            &\quad \geq (i-1-z_{\ell-1})c+b \geq \min\{(i-1-z_{\ell-1})c+b,0\}.
        \end{align*}
        \item If $D(i-1,k+1,\ell-1) > 0$,
        \begin{align*}
            \min\{D(i-1,k+1,\ell-1),0\} = 0 \geq \min\{(i-1-z_{\ell-1})c+b,0\}.
        \end{align*}
    \end{enumerate}
    Substituting \eqref{eq:induction_lb_i-1_k_ell} and \eqref{eq:min_zero_H} into the recursive equation \eqref{eq:recursive_neg} for $\ell<C_2$ yields
    \begin{align}
        \lefteqn{D(i,k,\ell) = \tilde{p}_\ell(i\tilde{c}_\ell+\tilde{b}_\ell) + \tilde{q}_\ell D(i-1,k,\ell) + \tilde{r}_\ell \min\{D(i-1,k+1,\ell-1), 0\}}& \nonumber \\
        & \quad \geq \tilde{p}_\ell(i\tilde{c}_\ell+\tilde{b}_\ell) + \tilde{q}_\ell \big((i-1-z_\ell)c+b\big) + \tilde{r}_\ell \min\{(i-1-z_{\ell-1})c+b, 0\} \nonumber \\
        & \quad = \tilde{p}_\ell(ic+b) + \tilde{q}_\ell \big((i-1-z_\ell)c+b\big) + \tilde{r}_\ell \min\{(i-1-z_{\ell-1})c+b, 0\}, \label{eq:lb_min}
    \end{align}where the last step holds since $c_\ell = c$ and $b_\ell = b$ for $\ell<C_2$ (recall definitions of $c$ and $b$ in Definition \ref{def:probs_quantities_ell}). There are two cases to consider based on whether or not $i \leq \tilde{i}_{\ell-1}$, where $\tilde{i}_{\ell-1} := \left\lfloor\frac{-b}{c} + z_{\ell-1}\right\rfloor+1$.
    \begin{enumerate}[label= \textbf{Case} \arabic*:, leftmargin=3\parindent]
        \item If $i \leq \tilde{i}_{\ell-1}$ ,i.e., $(i-1-z_{\ell-1})c+b \leq 0$ (recall $c>0$), then \eqref{eq:lb_min} becomes
        \begin{align}
        \lefteqn{D(i,k,\ell) \geq \tilde{p}_\ell(ic+b) + \tilde{q}_\ell \big((i-1-z_\ell)c+b\big) + \tilde{r}_\ell \big((i-1-z_{\ell-1})c+b\big)}& \nonumber \\
        & \quad = (\tilde{p}_\ell+\tilde{q}_\ell+\tilde{r}_\ell)(ic+b) - \big(\tilde{q}_\ell(z_\ell+1) + \tilde{r}_\ell(z_{\ell-1}+1)\big)c \nonumber \\
        &\quad = (ic+b) - \big(\tilde{q}_\ell(z_\ell+1) + \tilde{r}_\ell(z_{\ell-1}+1)\big)c, \label{eq:heur_linear_lb_case1_inte}
        \end{align}where the last equality is by the definition that $\tilde{p}_\ell+\tilde{q}_\ell+\tilde{r}_\ell=1$ (see Definition \ref{def:probs_quantities_ell}). Consider the expression in the second term with coefficient $-c$ in \eqref{eq:heur_linear_lb_case1_inte}.
        \begin{align}
            \lefteqn{\tilde{q}_\ell(z_\ell+1) + \tilde{r}_\ell(z_{\ell-1}+1)}& \nonumber \\
            &\quad = \frac{k-1}{(k-1)+(\ell+1)m} \left( \frac{k-1}{m}+\ell+1\right) + \frac{\ell m}{k + \ell m} \left( \frac{k}{m}+(\ell-1)+1\right)\nonumber \\
            &\quad = \frac{k-1}{m}+\ell \nonumber \\
            &\quad = z_\ell, \label{eq:z_rearrange}
        \end{align}where the first equality applies the definitions of $\tilde{q}_\ell$ and $\tilde{r}_\ell$ in Definition \ref{def:probs_quantities_ell} and $z_\ell$ in \eqref{def:sequence_z}, and the last equality applies \eqref{def:sequence_z} again. Applying \eqref{eq:y_rearrange} in \eqref{eq:heur_linear_lb_case1_inte} yields that for $\ell < C_2$, if $D(i,k,\ell) \leq 0$ and $i \leq \tilde{i}_{\ell-1}$, then
        \begin{align}
            D(i,k,\ell) \geq (i-z_\ell)c+b. \label{eq:heur_linear_lb_case1}
        \end{align}
        \item If $i \geq \tilde{i}_{\ell-1} + 1$, that is, $(i-1-z_{\ell-1})c+b > 0$, we now prove by contradiction that $D(\tilde{i}_{\ell-1}+1,k,\ell) > 0$. If not, we know $D(\tilde{i}_{\ell-1},k,\ell) \leq D(\tilde{i}_{\ell-1}+1,k,\ell) \leq 0$ by Proposition \ref{prop:larger_mu2_inc}, leading to $D(\tilde{i}_{\ell-1},k,\ell) \geq (\tilde{i}_{\ell-1}-z_\ell)c+b$ by \eqref{eq:heur_linear_lb_case1} in the first case. Consequently, applying \eqref{eq:lb_min} at $(\tilde{i}_{\ell-1}+1,k,\ell)$ and noticing that $\big((\tilde{i}_{\ell-1} + 1)-1-z_{\ell-1}\big)c+b > 0$ result in
        \begin{align}
            \lefteqn{
            D(\tilde{i}_{\ell-1}+1,k,\ell) \geq \tilde{p}_\ell\big((\tilde{i}_{\ell-1}+1)c+b\big) + \tilde{q}_\ell \big((\tilde{i}_{\ell-1}-z_\ell)c+b\big)}& \nonumber \\
            &\quad = (\tilde{p}_\ell+\tilde{q}_\ell)\big((\tilde{i}_{\ell-1}+1)c+b\big) - \tilde{q}_\ell(1+z_\ell)c \nonumber \\
            &\quad = (\tilde{p}_\ell+\tilde{q}_\ell)\left(\left(\tilde{i}_{\ell-1}+1 - \frac{\tilde{q}_\ell}{\tilde{p}_\ell+\tilde{q}_\ell}(1+z_\ell)\right)c+b\right) \nonumber \\
            &\quad > (\tilde{p}_\ell+\tilde{q}_\ell)\left(\left(-\frac{b}{c}+z_{\ell-1}+1 - \frac{\tilde{q}_\ell}{\tilde{p}_\ell+\tilde{q}_\ell}(1+z_\ell)\right) \cdot c+b\right) \nonumber \\
            & \quad = (\tilde{p}_\ell+\tilde{q}_\ell)\left(z_{\ell-1}+1 - \frac{\tilde{q}_\ell}{\tilde{p}_\ell+\tilde{q}_\ell}(1+z_\ell)\right) \cdot c, \label{eq:case2_pos}
        \end{align}where the strict inequality holds by noticing $\tilde{i}_{\ell-1} = \left\lfloor\frac{-b}{c}+z_{\ell-1}\right\rfloor+1 >\frac{b}{-c}+z_{\ell-1}$ and $c>0$.
        Moreover, the term in the second parenthesis in \eqref{eq:case2_pos} is non-negative:
        \begin{align*}
            \lefteqn{z_{\ell-1}+1 - \frac{\tilde{q}_\ell}{\tilde{p}_\ell+\tilde{q}_\ell}(1+z_\ell)}&\\
            &\quad \geq z_\ell - \frac{\tilde{q}_\ell}{\tilde{p}_\ell+\tilde{q}_\ell}(1+z_\ell)\\
            &\quad = \frac{\tilde{p}_\ell}{\tilde{p}_\ell+\tilde{q}_\ell}\left(z_\ell-\frac{\tilde{q}_\ell}{\tilde{p}_\ell}\right)\\
            &\quad \geq 0,
        \end{align*}where the first inequality applies Statement \ref{state:heur_mono_z} in Lemma \ref{silemma:heur}, and the last inequality uses Statement \ref{state:z_fraction} in Lemma \ref{lemma:heur}.
        Along with $c>0$ in \eqref{eq:case2_pos}, it holds that $D(\tilde{i}_{\ell-1}+1,k,\ell) > 0$, which contradicts our hypothesis that $D(\tilde{i}_{\ell-1}+1,k,\ell) \leq 0$, thus verifying that $D(\tilde{i}_{\ell-1}+1,k,\ell) > 0$.
        
        A final observation is that $D(i,k,\ell) \geq D(\tilde{i}_{\ell-1}+1,k,\ell) > 0$ for all $i \geq \tilde{i}_{\ell-1}+1$ by Proposition \ref{prop:larger_mu2_inc}, and therefore, if $D(i,k,\ell) \leq 0$, then $i \leq \tilde{i}_{\ell-1}$ is a must.       
    \end{enumerate}
    Concluding both cases, for $\ell < C_2$, $D(i,k,\ell) \leq 0$ implies that $i \leq \tilde{i}_{\ell-1}$ and $D(i,k,\ell)\geq (i-z_\ell)c+b$. The result is now established by induction.
\end{proof}

\end{document}